\documentclass[11pt]{article}

\usepackage{amsmath}
\usepackage{amssymb}
\usepackage{latexsym}
\usepackage{amsmath, amsfonts,amssymb, amsthm, euscript,makeidx,color,mathrsfs}
\usepackage[english]{babel}
\usepackage{rotating}
\usepackage{tikz}
\usetikzlibrary{positioning} 
\usepackage{tcolorbox} 
\usepackage{subcaption}
\usepackage{tabularx}
\usepackage{arydshln,leftidx,mathtools}

\newtheorem{assumption}{Assumption}
\def\qed{ \ \vrule width.2cm height.2cm depth0cm\smallskip}

\newcommand{\la}{\langle}
\newcommand{\ra}{\rangle}

\newcommand{\ba}{\begin{array}}
\newcommand{\ea}{\end{array}}
\newcommand{\be}{\begin{equation}}
\newcommand{\ee}{\end{equation}}
\newcommand{\bea}{\begin{eqnarray}}
\newcommand{\eea}{\end{eqnarray}}
\newcommand{\beaa}{\begin{eqnarray*}}
\newcommand{\eeaa}{\end{eqnarray*}}

\def\dbE{\mathbb{E}}
\def\dbF{\mathbb{F}}

\def\dbP{\mathbb{P}}

\def\dbR{\mathbb{R}}
\def\dbS{\mathbb{S}}
\def\dbT{\mathbb{T}}
\def\dbV{\mathbb{V}}
\def\dbX{\mathbb{X}}

\def\a{\alpha}
\def\b{\beta}
\def\g{\gamma}
\def\d{\delta}
\def\e{\varepsilon}

\def\k{\kappa}
\def\l{\lambda}

\def\si{\sigma}

\def\f{\varphi}
\def\th{\theta}
\def\o{\omega}

\def\G{\Gamma}
\def\D{\Delta}
\def\Th{\Theta}
\def\L{\Lambda}

\def\O{\Omega}
\def\cF{{\cal F}}

\def\cI{{\cal I}}
\def\cJ{{\cal J}}
\def\cK{{\cal K}}

\def\cN{{\cal N}}

\def\cS{{\cal S}}

\def\cV{{\cal V}}

\def\no{\noindent}

\def\ms{\medskip}
\def\bs{\bigskip}
\def\q{\quad}
\def\qq{\qquad}

\def\pa{\partial}
\def\cd{\cdot}
\def\cds{\cdots}

\def\td{\nabla}

\def\bx{{\bf x}}

\def\qed{ \hfill \vrule width.25cm height.25cm depth0cm\smallskip}

\newcommand{\basa}{\begin{assumption}}
\newcommand{\easa}{\end{assumption}}

\newcommand{\bas}{\begin{assum}}
\newcommand{\eas}{\end{assum}}

\def\limP2{\,\mathop{\buildrel \Pi_2\over\longrightarrow\,}}

\def\pa{\partial}

 \def\cd{\cdot}
\def\cds{\cdots}

\def\esup{\mathop{\rm ess\;sup}}

\def\dis{\displaystyle}

\def\bx{{\bf x}}

\def\1{{\bf 1}}

\def\bla{{\big \la}}
\def\bra{{\big \ra}}

\def\:{\!:\!}
\def\reff#1{{\rm(\ref{#1})}}
\def \proof{{\noindent \bf Proof\quad}}

at 9pt

\begin{document}

\newtheorem{thm}{Theorem}[section]
\newtheorem{lem}[thm]{Lemma}
\newtheorem{cor}[thm]{Corollary}
\newtheorem{prop}[thm]{Proposition}
\newtheorem{rem}[thm]{Remark}
\newtheorem{eg}[thm]{Example}
\newtheorem{defn}[thm]{Definition}
\newtheorem{assum}[thm]{Assumption}

\renewcommand {\theequation}{\arabic{section}.\arabic{equation}}
\def\thesection{\arabic{section}}

\title{ Cubature Method for Stochastic Volterra Integral Equations}
\author{ Qi Feng\thanks{\noindent Department of
Mathematics, University of Michigan, Ann Arbor, 48109-1043; email: qif@umich.edu. This author is supported in part by NSF grant DMS-2306769.}
~ and ~ Jianfeng Zhang\thanks{\noindent Department of Mathematics, University of Southern California, Los Angeles, 90089; email: jianfenz@usc.edu. This author is supported in part by NSF grant DMS-1908665 and DMS-2205972.} \thanks{\noindent The authors would like to thank two anonymous referees for their constructive comments which have helped to improve the paper greatly.}
}
\maketitle

{
\abstract{In this paper, we introduce the cubature formula for Stochastic Volterra Integral Equations. We first derive the stochastic Taylor expansion in this setting, by utilizing a functional It\^{o} formula, and provide its tail estimates. We then introduce the cubature measure for such equations, and construct it explicitly in some special cases, including a long memory stochastic volatility model. We shall provide the error estimate rigorously. Our numerical examples show that the cubature method is much more efficient than the Euler scheme,  provided certain conditions are satisfied.}

\bs\bs

\no{\bf MSC2020.} 60H20, 65C30, 91G60

\bs
\no{\bf Keywords}. Stochastic Volterra integral equations, cubature formula, stochastic Taylor expansion,  fractional stochastic volatility model,  rough volatility model

}

\section{Introduction}
\label{sect-Introdcution}
\setcounter{equation}{0}
Consider a  stock price in a Brownian  setting under risk neutral measure $\dbP$:
\bea
\label{S}
dS_t = S_t \si_t dB_t,
\eea
In the Black-Scholes model,  the volatility process $\si_t \equiv \si_0$ is a constant. There is a large literature on stochastic volatility models where $\si$ is also a diffusion process, see e.g. Fouque-Papanicolaou-Sircar-Solna \cite{FPSS}. Strongly supported by empirical studies, the fractional stochastic volatility models and rough volatility models have received very strong attention in recent years, where $\si$ satisfies the following {Stochastic Volterra Integral Equation (SVIE)}:
\bea\label{si}
\si_t=\si_0+\int_0^tK(t,r)V_0(\si_r)dr + \int_0^t K(t,r)V_1(\si_r)d\tilde B_r.
\eea
Here $\tilde B$ is another Brownian motion possibly correlated with $B$, $V_i$'s are appropriate deterministic functions, and the deterministic  two time variable function $K$ has a Hurst parameter $H>0$ in the sense that $K(t,r) \sim (t-r)^{H-{1\over 2}}$ and $\pa_t K(t,r) \sim (t-r)^{H-{3\over 2}}$ when $t-r>0$ is small. Such a model was first proposed by Comte-Renault \cite{CR} for $H>{1\over 2}$, to model the long memory property of the volatility process. Another notable work is Gatheral-Jaisson-Rosenbaum \cite{GJR}, which finds market evidence that volatility's high-frequency behavior could be modeled as a rough path with $H<{1\over 2}$. We remark that one special case of \reff{si} is the fractional Brownian motion, where $V_0\equiv 0$ and $V_1\equiv 1$, see, e.g. Nualart \cite{Nualart}.

Our goal of this paper is to understand and more importantly to  numerically compute the option price in this market:  assuming  zero interest rate for simplicity,
\bea
\label{Y0}
\dbE[G(S_T)].
\eea 
Note that the volatility process $\si$ in \reff{si} is in general neither a Markov process nor a semimartingale\footnote{When $H>1$, $X$ is actually a semimartingale, see e.g. \cite{VZ17}. However, it is still highly non-Markovian, so the numerical challenge remains in this case.}. Consequently, $S_t$ is highly non-Markovian, in the sense that one cannot Markovianize it by adding finitely many extra states, and correspondingly the option price is characterized as a path dependent PDE (PPDE, for short), see Viens-Zhang \cite{VZ17}. This imposes significant challenges, both theoretically and numerically. Indeed, compared to the huge literature on numerical methods for PDEs, there are very few works on efficient numerical methods for such PPDEs. Besides the standard Euler scheme Zhang \cite{ZhangXicheng}, we refer to Wen-Zhang \cite{WZ} for an improved rectangular method; Jacquier-Oumgari \cite{jacquier2019deep} and Ruan-Zhang \cite{ruan2020} 
on numerical methods for high dimensional (nonlinear) PPDEs driven by SVIEs;
Richard-Tan-Yang \cite{RTY2021, RTY2023} on discrete-time simulation schemes, including the Euler and Milstein schemes, and the corresponding Multi-Level Monte-Carlo method;  Ma-Yang-Cui \cite{MYC} by using Markov chain approximation;
 and  Alfonsi-Kebaier \cite{AlfonsiKebaier2021},  Bayer-Breneis \cite{BayerBreneis2023}, and Harms \cite{Harms} by using Laplace transform for singular kernel functions. In recent years, there has also been a growing interest on the convergence analysis and error estimates for stochastic Volterra integral equations, see, e.g. Bayer-Fukasawa-Nakahara \cite{BayerFukasawaNakahara}, Bayer-Hall-Tempone \cite{BayerHallTempone}, Bonesini-Jacquier-Pannier \cite{BJP2023}, Friz-Salkeld-Wagenhofer \cite{FrizSalkeldWagenhofer}, Fukasawa-Ugai \cite{FT2021},  Gassiat \cite{Gassiat},  Li-Huang-Hu \cite{LHH2022},  and Nualart-Saikia \cite{Nualart2022}.
 In this paper, we propose the cubature method for the above option price \reff{Y0}. This is a deterministic method, and our numerical examples show that, under certain conditions, it is much more efficient than the simulation methods such as the Euler scheme.

The cubature method was first introduced by the seminal works Lyons-Victoir \cite{LV03} and Litterer-Lyons \cite{LL04} for diffusion processes, see also   Gyurk{\'o}-Lyons \cite{GL11},  Litterer-Lyons \cite{LL12}, Ninomiya-Shinozaki \cite{NinomiyaShinozaki}, and Ninomiya-Victoir \cite{NinomiyaVictoir}  for its extensive numerical implementations. The method builds upon the stochastic Taylor expansion for smooth $G$: 
\bea
\label{STE}
G(S_T) = I_N + R_N,
\eea
see, e.g. Kloeden-Platen \cite{KP}, where $I_N$ is a linear combination of multiple integrals against the Brownian motion $B$ (typically in Stratonovich form), called the signatures of $B$, and $R_N$ is the remainder term. The main idea is to introduce a discrete measure $Q$ to match the expectations of the signatures: recalling that $\dbE = \dbE^\dbP$ is the expectation under $\dbP$,
\bea
\label{PQ}
\dbE[I_N] = \dbE^Q[I_N],
\eea
Then we will have an approximation $\dbE[G(S_T)] \approx \dbE^Q[G(S_T)]$. Since $Q$ is discrete and it is easy to compute the exact value of $\dbE^Q[G(S_T)]$ (without involving simulations), the algorithm is very efficient, provided sufficient technical conditions to make the approximation error small enough.

In this paper we shall consider the general SVIE, see \reff{volterra sde} below, and our goal is to approximate $\dbE[G(X_T)]$. To introduce the cubature method for $X$,  our first step is to derive the stochastic Taylor expansion in this setting. Note that \reff{STE} relies heavily on the It\^{o} formula, but the solution $X$ to the SVIE is not a semimartinagle, which prohibits us from applying the It\^{o} formula directly.  To overcome this difficulty, we utilize an auxiliary two time variable process $\Th^s_t$ introduced by {Wang \cite{Wang}} and Viens-Zhang \cite{VZ17}, see \reff{volterra sde theta} below. This process satisfies $\Th^t_t = X_t$ and enjoys the desired semimartingale property:  for fixed $s$, the process $t\in [0, s] \mapsto \Th^s_t$ is a semimartingale. In particular, \cite{VZ17} established a functional It\^{o} formula, which enables us to derive the desired stochastic Taylor expansion, with more involved signatures for the SVIEs than the diffusion case. We then introduce a discrete cubature measure $Q$ for $X$, in the pirit of \reff{STE} and \reff{PQ}, and prove the following error estimate: for some constant $C_N$ which depends on the regularity of the coefficients, 
\bea
\label{err1}
\Big|\dbE[G(X_T)] - \dbE^Q[G(X_T)]\Big| \le C_N T^{N+1\over 2}.
\eea

The above result is desirable only when $T$ is small. For general $T$, we follow the idea of \cite[Theorem 3.3]{LV03} and utilize the flow property of the path dependent value function established in \cite{VZ17}. To be precise, we consider a uniform partition of $[0, T]$: $0=T_0<\cds<T_M=T$, and construct a cubature measure $Q_m$ on each subinterval $[T_m, T_{m+1}]$. Let $Q$ be the independent composition of $\{Q_m\}_{0\le m<M}$, we then have the following estimate:
 \bea
\label{err2}
\Big|\dbE[G(X_T)] - \dbE^Q[G(X_T)]\Big| \le C_N {T^{N+1\over 2}\over M^{N-1\over 2}},
\eea
where $C_N$ is  independent of $M$.  The above estimate clearly converges to $0$ as $M\to \infty$.

We remark that, while our stochastic Taylor expansion can be developed for any kernel $K$ with Hurst parameter $H>0$, the cubature formula becomes much more subtle when $H<{1\over 2}$. In this paper we restrict to the case $H\ge {1\over 2}$, and leave the  case $H<{1\over 2}$ to future study. For applications, we refer to Comte-Renault \cite{CR} for the long memory model with ${1\over 2}<H<1$,  Gulisashvili-Viens-Zhang \cite{GVZ} for the integrated variance model with $1<H<2$, and El Omari \cite{ElOmari}  for the mixed fractional Brownian motion model with more general $H$. We also refer to Beran \cite[Section 4.2]{Beran} for applications in hydrology,  Loussot-Harba-Jacquet-Benhamou-Lespesailles-Julien \cite{LHJBLJ} for applications in image processing, Gupta-Singh-Karlekar \cite{GSK} for applications in signal classification, Blu-Unser \cite{BU} for fractional spline estimators,  and Perrin-Harba-Berzin Joseph-Iribarren-Bonami \cite{PHBIB} for the theory of higher order fractional Brownian motions with general $H>1$. 
However, we should point out that our result does not cover the rough volatility models in Gatheral-Jaisson-Rosenbaum \cite{GJR} with $0<H<{1\over 2}$.  Moreover, the estimates \reff{err1} and  \reff{err2} require the coefficients to be sufficiently smooth, as we will specify in the paper. In particular, the constant $C_N$ will depend on such regularity. 

The efficiency of our cubature method comes down to the construction of the cubature measure $Q$ in \reff{err2}, which will involve $(2W)^M$ deterministic paths for some constant $W$. When the dimension of $X$ is large, or when $N$ is large, the $W$ will be large; and when $T$ is large, in light of \reff{err2} we will require $M$ to be large. 
We refer to \S \ref{sect-conclude} for more precise comments on the efficiency issue. When all the conditions are satisfied so that $(2W)^M$ is at a reasonable level, our numerical examples show that the cubature method is much more efficient than the Euler scheme.  

We should remark that the above efficiency issue was already present for the cubature method in the standard Brownian setting. There have been great efforts in the literature to overcome this difficulty and to apply the idea of the cubature method to more general models, see, e.g., Crisan-Manolarakis  \cite{CM12, CM14},  Crisan-McMurray \cite{crisan2019cubature}, Raynal-Trillos \cite{de2015cubature},  Filipovi\'c-Larsson-Pulido \cite{filipovic2020markov}, Foster-Lyons-Oberhauser \cite{FLO},  Foster-Reis-Strange \cite{FosterReisStrange}, and  Hayakawa-Oberhauser-Lyons \cite{HayakawaOberhauserLyons}.
It will be very interesting to explore if these ideas can help to improve the cubature method in the Volterra framework. We would also like to mention the very interesting connection between the signature, the kernel method, and machine learning, see Chevyrev-Oberhauser \cite{CO18}, Kidger-Bonnier-Arribas-Salvi-Lyons \cite{kidger2019deep}, Kiraly-Oberhauser  \cite{KO16}, Liao-Lyons-Yang-Ni \cite{liao2019learning},  and the references therein.  

Finally, we note that, while sharing many properties, the SVIE \reff{si} is different from the following SDE driven by a fractional Brownian motion $B^H_t:= \int_0^t K(t,r) dr + \int_0^t K(t,r) d\tilde B_r$:
\bea\label{si'}
\si'_t=\si_0+\int_0^tV_0(\si'_r)dr + \int_0^t V_1(\si'_r)d B^H_r.
\eea
We refer to Baudoin-Coutin \cite{BC07} and Passeggeri \cite{Passeggeri16}  for some works on signatures for fractional Brownian motions, and {Harang-Tindel \cite{harang2019volterra, harang2021volterra_2}
 on signatures defined for `` Volterra path"}. We shall remark that, unlike our signature which is directly for the solution $\si_t$ to the SVIE \reff{si} (instead of for the driving Brownian motion $\tilde B$), these signatures are for the driving fractional Brownian motion $B^H$ {or ``Volterra path"} which has much simpler structure. In particular, their signatures do not lead to the desired stochastic Taylor expansion which is crucial for the  cubature method. 
 
The rest of the paper is organized as follows.  In \S \ref{sect-Taylor} we derive the stochastic Taylor expansions for the general SVIEs and prove the tail estimate. In \S \ref{sect-cubature} we introduce the cubature formula  when $T$ is small, and in \S \ref{sect-cubatureM} we modify the cubature formula when $[0, T]$ is decomposed into $M$ parts. We construct the cubature measure $Q$ explicitly for a one dimensional SVIE in \S \ref{sect-1d} , and for the $2$-dimensional fractional stochastic  volatility model in \S \ref{sect-rough}. In \S\ref{sect-numerical} we present various numerical examples and compare its efficiency with the Euler scheme. Finally, we present some technical proofs in Appendix.

\section{Stochastic Taylor expansions} 
\label{sect-Taylor}
\setcounter{equation}{0}
Throughout this paper, let $(\O, \cF, \dbF, \dbP)$ be a filtered probability space, $B^0_t:=t$, and $B=(B^1, \cds, B^d)$ a $d$-dimensional Brownian motion. Let $T>0$ be a fixed terminal time. We consider $d_1$-dimensional state process $X= (X^1,\cds, X^{d_1})$ solving the following SVIE under Stratonovich integration $\circ$ : given $x=(x_1,\cds, x_{d_1})\in \dbR^{d_1}$,
\bea
\label{volterra sde}
X^i_t = x_i + \sum_{j=0}^d \int_0^t K_i(t,r) V^i_j(X_r) \circ dB^j_r,\q i=1,\cds, d_1,
\eea
and we are interested in the efficient numerical computation of 
\bea
\label{Y0EG}
Y_0:=\dbE\big[G(X_T)\big].
\eea

Throughout the paper, the following hypotheses will always be enforced: for some $N\ge 1$ which will be specified in contexts,

\ms
\no {\bf (H0)} Each $K_i: \{(t, r): 0\le r\le t\le T\} \to [0, \infty)$ is infinitely smooth on $\{r<t\}$, and either $K_i \equiv 1$ or $K_i$ has Hurst parameter $H_i>{1\over 2}$, that is, $K_i(t,r) \sim (t-r)^{H_i-{1\over 2}}$ and $\pa_t K_i(t,r) \sim  (t-r)^{H_i-{3\over 2}}$ when $t-r>0$ is small.

\ms

\no {\bf (HN)} The functions $V^i_j, G \in C^N(\dbR^{d_1}; \dbR)$ with all the derivatives up to the order $N$ bounded.  

\ms

\no For later purpose, we will also need the following stronger version of {\bf (H0)}.

\ms
\no {\bf (H0-N)} Each $K_i: \{(t, r): 0\le r\le t\le T\} \to [0, \infty)$ is infinitely smooth on $\{r<t\}$, and either $K_i \equiv 1$ or $K_i$ has Hurst parameter $H_i>N+{1\over 2}$, in the sense that ${\pa^{\a+\b}\over \pa t^\a \pa r^\b} K_i(t,r) \sim (t-r)^{H_i-\a-\b-{1\over 2}}$  when $t-r>0$ is small, for all integers $\a,\b\ge 0$ such that $\a + \b \le N$.

\begin{rem}
\label{rem-SDE1}
{\rm (i) When $K_i$ has Hurst parameter $H_i>{1\over 2}$, $X^i$ is H\"{o}lder-$(H_i\wedge 1-\e)$ continuous in $t$ for any small $\e>0$. This implies that $V(X^i_t) \circ dB^j_t = V(X^i_t) dB^j_t$ for any smooth function $V$, where $V(X^i_t) dB^j_t$ denotes the It\^{o} integral, and they coincide with the Young's pathwise integral. On the other hand, when $K_i \equiv 1$, $X^i$ is clearly a semimartingale. So, letting $I_0$ denote the set of  $i>0$ such that $K_i \equiv 1$, we may rewrite \reff{volterra sde} in It\^{o}'s form, and in particular they are wellposeded under {\bf (H0)} and {\bf (H2)}:
\bea
\label{volterra sde: stra}
X^i_t = x_i + \sum_{j=0}^d \int_0^t K_i(t,r) V^i_j(X_r)  dB^j_r + {1\over 2} \sum_{k\in I_0} \int_0^t K_i(t,r) \pa_{k} V^i_j(X_r) V^k_j(X_r) dr.
\eea

(ii) Consider a special case: $x_1=0$, $K_1 \equiv 1$, $V^1_0 \equiv 1$, $V^1_j \equiv 0$ for $j\ge 2$. Then we can easily see that $X^1_t =t$. So the system \reff{volterra sde} actually covers the case that the coefficients $V^i_j$ depends on the time variable $t$. 
\qed}
\end{rem}

\begin{rem}
\label{rem-SDE2}
{\rm (i) In fractional stochastic volatility models, where $X^i$ is interpreted as the volatility (or variance) process of certain underlying asset price, the assumption ${1\over 2}<H_i<1$ implies that the volatility has ``long memory", see Comte-Renault \cite{CR}.   We also refer to \cite{Beran, BU,  ElOmari, GVZ, GSK, LHJBLJ,   PHBIB}  for applications and theory when $H_i>1$.

(ii) The case  $H_i < {1\over 2}$, supported by the empirical studies in Gatheral-Jaisson-Rosenbaum \cite{GJR}, has received very strong attention in the mathematical finance literature in recent years. The singularity of $K_i$ in this case will make the theory much more involved, for example one may need to consider the weak solution to \reff{volterra sde}, and consequently the numerical algorithms will be less efficient. We shall leave this important and challenging case to future study.
\qed}
\end{rem}

\subsection{The functional It\^{o} formula}

 Note that $X^i$ is not a semimartingale when $K_i \neq 1$, which prohibits us from applying many stochastic analysis tools such as the It\^{o} formula directly. To get around of this difficulty, in this subsection we introduce a functional It\^{o} formula, which is established in Viens-Zhang \cite{VZ17} but tailored for the purpose of this paper. 
 
Denote $\dbX_t:= C^0([t, T]; \dbR^{d_1})\cap C^1((t, T]; \dbR^{d_1})$ for each $t\in [0, T]$, equipped with the uniform norm. For each $t\in [0, T]$ and $\phi: \dbX_t \to \dbR$, let $\pa_\bx \phi$ denote the Fr\'echet derivative of $u$. That is, $\pa_\bx \phi(\bx): \dbX_t\to\dbR$ is a linear mapping  satisfying:
\bea
\label{paxu}
\phi(\bx + \eta) -  \phi(\bx) = \la \pa_\bx \phi(\bx), \eta\ra + o(\|\eta\|),\q \forall \bx, \eta\in \dbX_t.
\eea
Similarly we may define the second order derivate $\pa_{\bx\bx} \phi(\bx)$ as a bilinear mapping on $\dbX_t \times \dbX_t$:
\bea
\label{paxxu}
\la \pa_\bx \phi(\bx+\eta_2), \eta_1\ra -  \la \pa_\bx \phi(\bx), \eta_1\ra = \la \pa_{\bx\bx} \phi(\bx), (\eta_1, \eta_2)\ra + o(\|\eta_2\|),\q \forall \eta_1, \eta_2\in \dbX_t.
\eea
We may continue to define higher order derivatives $\pa^{(n)}_\bx \phi(\bx)$ in an obvious manner, and let $C^n(\dbX_t)$ denote the set of continuous functions $\phi:\dbX_t\to \dbR$ which has uniformly continuous derivatives up to order $n$.
Moreover, as in \cite{VZ17} (see also an earlier work \cite{Wang}), we introduce a two time variable process $\Theta^s_t = (\Th^{1,s}_t,\cds, \Th^{d_1, s}_t)$ for $0\le t\le s\le T$:
\bea\label{volterra sde theta}
\Theta_t^{i,s}=x_i+\sum_{j=0}^{d}\int_0^t K_i(s,r)V^i_j(X_r)\circ dB_r^j.	
\eea
This process enjoys the following nice properties:

\begin{itemize}

\item For fixed $s$, the process $t\!\in\! [0, s]\!\to\! \Th^s_t$ is an $\dbF$-progressively measurable semimartingale;

\item  For fixed $t$, the process $s\!\in\! [t, T] \!\to\! \Th^s_t$ is $\cF_t$-measurable, continuous on $[t, T]$, infinitely smooth on $(t, T]$, and  with ``initial" condition $\Th^t_t = X_t$. In particular, $\Th_t \in \dbX_t$, a.s.
\end{itemize}

Then we have the following functional It\^{o} formula, which is essentially the same as \cite[Theorem 3.10]{VZ17}, but in Stratonovic form instead of It\^o form.  
\begin{prop}
\label{prop-Ito}
Let {\bf (H0)} and {\bf (H2)} hold and $\phi\in C^2(\dbX_{T'})$ for some $0<T' <T$. Then
\bea
\label{Ito}
d \phi(\Th_t^{[T', T]}) = \sum_{i=1}^{d_1} \sum_{j=0}^d \la \pa_{\bx_i} \phi(\Th_t^{[T', T]}), K^{[T', T]}_{i, t})\ra  V^i_j(X_t) \circ dB^j_t,~ 0\le t\le T'.
\eea
Here $\Th_t^{[T', T]}$ and $K^{[T', T]}_{i, t}$ denote the paths $\Th^s_t, K_i(s, t)$, $s\in [T', T]$, respectively.
\end{prop}

We now turn to the problem \reff{Y0EG}. For any $t\in [0, T]$ and $\th\in \dbX_t$, introduce
\bea
\label{u}
 u(t, \th) = \dbE[G(X^{t,\th}_T)],\q X^{t, \th, i}_s = \th_s^i +\sum_{j=0}^d \int_t^s K_i(s, r) V^i_j(X^{t, \th}_r) \circ dB^j_r,~ i=1,\cds, d_1 .
\eea
Since $\th$ is differentiable, by Remark \ref{rem-SDE1} (i) it is clear that the above Volterra SDE is wellposed. Moreover,
\beaa
u(T, \tilde x) = G(\tilde x),~\forall \tilde x\in \dbX_T = \dbR^{d_1};\q Y_0 = u(0, x),~\mbox{where $x\in \dbX_0$ is a constant path},
\eeaa
and we have the following simple result, whose proof is postponed to Appendix.
 
  \begin{prop}
\label{prop-ureg}
Under {\bf (H0)} and {\bf (HN)}, we have $u(t,\cd)\in C^{N-1}(\dbX_t)$  for any $t\in [0, T]$.  Moreover, all the involved derivatives are bounded by $C_Ne^{C_NT}$, where $C_N$ depends only on the parameters in  {\bf (H0)} and {\bf (HN)}.
\end{prop}

\subsection{The stochastic Taylor expansion}
Fix $M\ge 1$, and set $T_m:= m\d$, $m=0,\cds, M$, where $\d:= \d_M:= {T\over M}$.  In this subsection we fix $m$ and consider the stochastic Taylor expansion of $u(T_{m+1},\cd)$ at $T_m$. We first introduce some notation: for any $n\ge 1$ and $s\in [T_m, T_{m+1}]$,
\bea
\label{Delta}
\dbT^m_n(s) :=\big\{ \vec t= (t_1,\cds, t_n): T_m\le t_n\le \cds\le t_1 \le s\big\},\q t^m_0 := T_{m+1},~ \dbT^m_n:= \dbT_n(t^m_0).
\eea

Assume $u(T_{m+1},\cd)$ is sufficiently smooth, by \reff{Ito} we have
\bea
\label{expansion1}
&&\dis u(T_{m+1},\Th_{T_{m+1}}^{[T_{m+1},T]}) = u(t^m_0,\Th_{t^m_0}^{[t^m_0,T]})\\
&&\dis =u(t^m_0, \Th^{[t^m_0,T]}_{T_m})+\sum_{i_1=1}^{d_1}\sum_{j_1=0}^d\int_{T_m}^{t^m_0}\bla\pa_{\bx_{i_1}}u(t^m_0,\Th^{[t^m_0,T]}_{t_1}),K_{i_1, t_1}^{[t^m_0,T]}  \bra V_{j_1}^{i_1}(\Th^{t_1}_{t_1})\circ dB_{t_1}^{j_1},\nonumber
\eea
where we used the fact that $X_{t_1} = \Th^{t_1}_{t_1}$.
Now fix $t_1\in [T_m, t^m_0]$, note that $t_2 \in [T_m, t_1]\mapsto \Th^{t_1}_{t_2}$ is a semimartingale, and $\th\in \dbX_{t^m_0}\mapsto \la\pa_{\bx_{i_1}}u(t^m_0,\th),K_{i_1, t_1}^{[t^m_0,T]} )  \ra$ is in $C^2(\dbX_{t^m_0})$. Then
\beaa
&&\dis dV_{j_1}^{i_1}(\Th^{t_1}_{t_2}) = \sum_{i_2=1}^{d_1}\sum_{j_2=0}^d \pa_{x_{i_2}}V_{j_1}^{i_1}(\Th^{t_1}_{t_2}) K_{i_2}(t_1,t_2)V_{j_2}^{i_2}(X_{t_2})  \circ    dB_{t_2}^{j_2},\\
&&\dis d\bla\pa_{\bx_{i_1}}u(t^m_0,\Th^{[t^m_0,T]}_{t_2}),K_{i_1, t_1}^{[t^m_0,T]} \bra \\
&&\dis = \sum_{i_2=1}^{d_1}\sum_{j_2=0}^d\bla\pa_{\bx_{i_1}\bx_{i_2}}u(t^m_0,\Th_{t_2}^{[t^m_0,T]}), ~(K_{i_1, t_1}^{[t^m_0,T]},K_{i_2, t_2}^{[t^m_0,T]}) \bra V_{j_2}^{i_2}(\Th^{t_2}_{t_2}) \circ   dB_{t_2}^{j_2}.
\eeaa
Apply It\^{o}'s formula and plug these into \reff{expansion1}, we obtain
\bea
\label{expansion2}
&&\dis \!\!\!\!\!\!\!\! u(t^m_0,\Th_{t^m_0}^{[t^m_0,T]})=u(t^m_0, \Th^{[t^m_0,T]}_{T_m})+\sum_{i_1=1}^{d_1}\sum_{j_1=0}^d\int_{T_m}^{t^m_0}\bla\pa_{\bx_{i_1}}u(t^m_0,\Th^{[t^m_0,T]}_{T_m}),K_{i_1, t_1}^{[t^m_0,T]}  \bra V_{j_1}^{i_1}(\Th_{T_m}^{t_1})\circ dB_{t_1}^{j_1}\nonumber\\
&&\dis \!\!\!\!\!\!\!\!~ +\sum_{i_1,i_2=1}^{d_1}\sum_{j_1,j_2=0}^d\int_{\dbT^m_2}\Big[\big\la\pa_{\bx_{i_1}\bx_{i_2}}u(t^m_0,\Th_{t_2}^{[t^m_0,T]}),(K_{i_1, t_1}^{[t^m_0,T]}, K_{i_2, t_2}^{[t^m_0,T]}) \big\ra~ V_{j_1}^{i_1}(\Th^{t_1}_{t_2})V_{j_2}^{i_2}(\Th^{t_2}_{t_2}) \\
&&\dis\qq +\big\la\pa_{\bx_{i_1}}u(t^m_0,\Th^{[t^m_0,T]}_{t_2}), K_{i_1, t_1}^{[t^m_0,T]} \big \ra ~K_{i_2}(t_1,t_2) \pa_{x_{i_2}} V_{j_1}^{i_1}(\Th^{t_1}_{t_2}) V_{j_2}^{i_2}(\Th^{t_2}_{t_2})\Big] \circ dB_{t_2}^{j_2} \circ   dB_{t_1}^{j_1}.\nonumber
\eea

The formulae \reff{expansion1} and \reff{expansion2} are the first order and second order expansions of $u(t^m_0,\Th_{t^m_0}^{[t^m_0,T]})$. For higher order expansions,  we introduce the following notation. For any $n\ge 1$, denote $\cI_n :=  \{1,\cds, d_1\}^n$ with elements $\vec i = (i_1,\cds, i_n)$, $\cJ_n:= \{0,\cds, d\}^n$ with elements $\vec j=(j_1,\cds, j_n)$,   and introduce a set of mappings for the indices:
  \beaa
 \cS_n := \Big\{\vec\k = (\k_1,\cds, \k_n): \k_l \in \{0, 1, \cds, l-1\}, l=1,\cds, n\Big\}.
 \eeaa
 Given $\vec i \in \cI_n, \vec j \in \cJ_n, \vec t \in \dbT_n, \vec \k\in \cS_n$, $\vec x=(x_1,\cds, x_n)\in (\dbR^{d_1})^{n}$, $\th\in \dbX_{t^m_0}$, and $\f: \dbT^m_n\to \dbR$, $\psi: \dbR^{d_1}\to \dbR$,  denote 
 \bea
  \label{cKcV}
 \left.\ba{c}
 \dis \cN_\a(\vec\k):=\{l\in \{1,\cds, n\}: \k_l =\a\},\q \a = 0,\cds, n,\ms\\
 \dis   \cK(\vec i,  \vec\k; \vec t) :=\prod_{l=1}^n  K_{i_l}(t_{\k_l}, t_l),\q \cK_+(\vec i,  \vec\k; \vec t) :=\prod_{\a=1}^n \prod_{l\in \cN_\a(\vec\k)} K_{i_l}(t_{\k_l}, t_l),\\
 \dis \vec \cK_0(\vec i,  \vec\k; \vec t) := \big( K^{[t^m_0, T]}_{i_l, t_l}\big)_{l\in \cN_0(\vec\k)} ,\ms\\
\dis \pa_{\vec i}^{\vec \k}u(t^m_0,\th):=\pa_{\bx_{i_{l_1}}}\cdots \pa_{\bx_{i_{l_k}}} u(t^m_0,\th),~\mbox{where $\{l_1, \cds, l_k\}= \cN_0(\vec\k)$},\ms \\
 \dis \pa^{\vec \k, \alpha}_{\vec i} \psi( x):= \pa_{x_{i_{l_1}}}\cds \pa_{x_{i_{l_k}}} \psi(x), ~\mbox{where $\{l_1, \cds, l_k\}= \cN_\a(\vec\k)$},\ms\\
\dis  \cV(\vec i, \vec j, \vec \k; \vec t, \vec x, \th):=  \prod_{\alpha = 1}^{n} \pa^{\vec \k, \alpha}_{\vec i}  V^{i_\a}_{j_\a} (x_\a) \cK_+(\vec i,  \vec \k; \vec t) ~\big \la \pa_{\vec i}^{\vec \k}u(t^m_0, \th) , \vec \cK_0(\vec i,  \vec\k; \vec t)  \big\ra ,\ms \\
 \dis\Th^{\vec t}_{s}:= (\Th^{t_1}_{s},\cds, \Th^{t_n}_{s}),~s\le t_n,\q  \f(\vec t)\circ dB^{\vec j}_{\vec t} :=  \f(\vec t)  \circ dB^{j_n}_{t_n} \circ \cds \circ dB^{j_1}_{t_1}.
 \ea\right.
 \eea
 Note that $\vec \cK_0$ and $\cV$ here actually depend on $m$, but we omit this dependence for notational simplicity. 
We then have the following expansion, whose proof is postponed to Appendix.

\begin{prop}
\label{prop-expansion}
For any $N\ge 1$, under {\bf (H0)} and {\bf (H(N+3))}, we have
{\small\bea\label{expansionN}
u(t^m_0, \Th^{[t^m_0,T]}_{t^m_0})&=& u(t^m_0, \Th^{[t^m_0,T]}_{T_m})+\sum_{n=1}^N\sum_{\vec i \in\cI_n, \vec j\in\cJ_n, \vec\k\in \cS_n}\int_{\dbT^m_n} \cV(\vec i,\vec j, \vec\k;\vec t, \Th^{\vec t}_{T_m},\Th^{[t^m_0,T]}_{T_m}) \circ dB_{\vec t}^{\vec j}\nonumber \\
&&+\sum_{\vec i\in\cI_{N+1}, \vec j\in\cJ_{N+1}, \vec \k\in\cS_{N+1}}\int_{\dbT^m_{N+1}} \cV(\vec  i,\vec j, \vec\k;\vec t, \Th^{\vec t}_{ t_{N+1}},\Th^{[t^m_0,T]}_{t_{N+1}}) \circ dB_{\vec t}^{\vec j}.
\eea }
\end{prop}

\subsection{The remainder estimate}
In this subsection we estimate the remainder term in Taylor expansion, which will provide guideline for our numerical algorithm later.  For an appropriate  function $\f: \dbT^m_n \times \O\to \dbR$ and $T_m\le s \le t^m_0$, denote
 \bea
 \label{fnnorm}
 \|\f(\cd)\|_{s, \vec j}^2 :=  \|\f(\cd)\|_{[T_m,s], \vec j}^2  :=  \dbE_m\Big[ \Big|\int_{\dbT^m_n(s)} \f(\vec t~\!) \circ dB^{\vec j}_{\vec t}\Big|^2\Big],\q\mbox{where}\q \dbE_m:= \dbE_{\cF_{T_m}}.
 \eea
 Moreover, for $\vec j\in \cJ_n$, $\vec t\in \dbT_n$, and $1\le l \le n$, denote
 \bea
 \label{jl}
 \vec j_l := (j_1, \cds, j_l),\q \vec j_{-l} := (j_{l+1},\cds, j_n),\q \vec t_{-l} := (t_{l+1},\cds, t_n).
 \eea
 We first have the following simple but crucial lemma, whose proof is postponed to Appendix.
\begin{lem}
\label{lem-EdB} 
Fix $n\ge 2$, $\vec j\in \cJ_n$,  and let $\f: \dbT^m_n \times \O\to \dbR$ be bounded, jointly measurable in all variables, and, for each $\vec t\in \dbT^m_n$,  $\f(\vec t)$ is $\cF_{t_n}$-measurable in $\o$. There exists a universal constant $C>0$ such that, for any $T_m\le s \le t^m_0$,
\bea
\label{EdB}
\dis \dbE_m\Big[ \int_{\dbT^m_n(s)}\!\!\!\!\!\! \f(\vec t~\!) \circ dB^{\vec j}_{\vec t}\Big]  = \left\{\ba{lll}
 \dis \int_{T_m}^s  \dbE_m\Big[\int_{\dbT^m_{n-1}(t_1)}\!\!\!\!\! \f(t_1, \vec t_{-1}) \circ dB^{\vec j_{-1}}_{\vec t_{-1}}\Big] dt_1,\qq j_1=0,\ms\\
 \dis 0,\qq\qq\qq\qq\qq\qq\qq\qq\q~  j_1\neq 0, j_2,\ms\\
 \dis {1\over 2} \int_{T_m}^s \dbE_m\Big[\int_{\dbT^m_{n-2}(t_1)}\!\!\! \!\!\! \f(t_1, t_1, \vec t_{-2}) \circ dB^{\vec j_{-2}}_{\vec t_{-2}}\Big] dt_1,~ j_1=j_2 >0.
 \ea\right.\\
 \label{E|dB|}
 \|\f(\cd)\|_{s, \vec j}^2 \le \left\{\ba{lll}
\dis C\d^2 \esup_{{T_m}\le s'\le s} \|\f(s', \cd)\|_{s', \vec j^{-1}}^2,\qq\qq\qq\qq\qq\qq\q j_1=0,\\
\dis C\d \esup_{{T_m}\le s'\le s} \|\f(s', \cd)\|_{s', \vec j^{-1}}^2,\qq\qq\qq\qq\qq\qq\q~ j_1\neq 0, j_2,\\
\dis C\d\esup_{{T_m}\le s'\le s} \|\f(s', \cd)\|_{s', \vec j^{-1}}^2 +C\d^2 \esup_{{T_m}\le s'\le s} \|\f(s',s',  \cd)\|_{s', \vec j^{-2}}^2,~ j_1= j_2> 0.\\
 \ea\right.
 \eea
\end{lem}

Note that $B^0$ and $\{B^j\}_{j\ge 1}$ contribute differently in \reff{EdB} and \reff{E|dB|}. Alternatively, we note that $B^0_t = t$ is Lipschitz continuous, but $B^j_t$ is H\"{o}lder-$({1\over 2}-\e)$ continuous for $j\ge 1$. To provide a more coherent error estimate, we shall modify \reff{expansionN} slightly. For any $1\le n\le N$ and $p\ge 1$, denote
{\small
\bea
\label{cJkm}
\left.\ba{c}
\dis \cJ_{n,N} :=\Big\{\vec j\in \cJ_n:  \|\vec j \|  \le N\Big\},\q\mbox{where}\q \|\vec j \| := n + \sum_{l=1}^n \1_{\{j_l =0\}},\\
\dis  A^m_N:= \sup_{n\le N} \sup_{\vec i \in \cI_n, \vec j\in \cJ_n,  \|\vec j\|=N, \vec t \in \dbT^m_n} \sup_{\vec x\in (\dbR^{d_1})^{n}, \th\in \dbX_{t^m_0}} \Big|\sum_{\vec\k \in \cS_{N}} \cV(\vec  i,\vec j, \vec\k;\vec t, \vec x, \th)\Big|.
\ea\right.
\eea
}

We then have the following tail estimate, whose proof is postponed to Appendix.

\begin{thm}
\label{thm-Taylor}
Let {\bf (H0)} and {\bf (H(N+3))} hold, and $R^m_N$ be determined by
\bea
\label{Taylor}
\left.\ba{c}
\dis u(t^m_0, \Th^{[t^m_0,T]}_{t^m_0})=  I^m_N + R^m_{N},\q\mbox{where}\\
\dis I^m_N:= u(t^m_0, \Th^{[t^m_0,T]}_{T_m})  + \sum_{n=1}^N  \sum_{\vec i\in \cI_n, \vec j\in \cJ_{n,N}, \vec \k\in \cS_n} \int_{\dbT^m_n} \cV(\vec i,\vec j, \vec\k;\vec t, \Th^{\vec t}_{T_m},\Th^{[t^m_0,T]}_{T_m}) \circ dB^{\vec j}_{\vec t}.
\ea\right.
\eea
Then there exists a constant $C_N>0$, which depends only on $N$ and $d, d_1$, such that
\bea
\label{RNest}
\big|\dbE_m[R^m_N]\big|\le  \Big(\dbE_m[|R^m_N|^2]\Big)^{1\over 2} \le  C_N \Big[ A^m_{N+1}\d^{N+1\over 2} + A^m_{N+2} \d^{N+2\over 2}+A^m_{N+3}\d^{N+3\over 2}\Big].
\eea
\end{thm}

\ms
\begin{rem}
\label{rem-Taylor}
{\rm Clearly, for fixed $N$, the error  in \reff{RNest} will be smaller  when $\d$ is smaller, when $G$ and $V^i_j$ are smoother (so that $u$ is smoother), and when the dimensions $d$ and $d_1$ are smaller (so that $C_N$ is smaller). This is consistent with our numerical results later. 
\qed}
\end{rem}

\section{The cubature formula: the one period case}
\label{sect-cubature}
\setcounter{equation}{0}
Note that \reff{RNest} is effective when $\d$ is small. In this section we consider the case that $T$ is small. Then we may simply set $M=1$ and thus $\d = T$. We shall apply the results in \S \ref{sect-Taylor} with $m=0$. In particular, in this case $\dbE_0 = \dbE$. For notational simplicity, in this section we shall omit the superscript $^0$, e.g. $\dbT_n= \dbT^0_n$, $I_N = I^0_N$, and $R_N= R^0_N$.

\subsection{Simplification of the stochastic Taylor expansion}
In this case we have: denoting $t_0 := T$, 
 \bea
 \label{cV0}
 \left.\ba{c}
\dis  \Th^t_{T_m} = \Th^t_0 = x,\q \dbX_{t^m_0} = \dbX_T = \dbR^{d_1},\q u(T, x) = G(x),~ x\in  \dbR^{d_1},\ms\\
\dis \big \la \pa_{\vec i}^{\vec \k}u(t^m_0, x) , \vec \cK_0(\vec i,  \vec\k; \vec t)  \big\ra = \pa_{\vec i}^{\vec \k, 0} G(x) \prod_{l\in \cN_0(\vec\k)} K_{i_l}(T, t_l),\ms\\
\dis \cV(\vec i, \vec j, \vec\k; \vec t, (x,\cds, x), x)=  \cV_0(\vec i, \vec j, \vec\k; x)\cK(\vec i, \vec\k; \vec t), \\
\dis \mbox{where}\q  \cV_0(\vec i, \vec j, \vec\k; x):= \prod_{\alpha = 1}^{n} \pa^{\vec \k, \alpha}_{\vec i}  V^{i_\a}_{j_\a}(x) \pa_{\vec i}^{\vec \k, 0} G(x).
\ea\right.
 \eea
Thus \reff{Taylor} becomes
 \bea\label{Taylor1}
G(X_T)= I_N+R_N = G(x) +\sum_{n=1}^N\sum_{\vec i \in\cI_n, \vec j\in\cJ_{n, N}, \vec\k\in \cS_n}\!\!\!\!\!\! \cV_0(\vec i, \vec j, \vec\k; x)\int_{\dbT_n} \cK(\vec i, \vec\k; \vec t) \circ dB_{\vec t}^{\vec j}+R_N.
\eea 
Moreover,  by abusing the notation we may modify $A^0_N$ and define $A_N$ as follows:
\bea
\label{AN}
 A_N:= \sup_{n\le N} \sup_{\vec i \in \cI_n, \vec j\in \cJ_n,  \|\vec j\|=N, \vec t \in \dbT^m_n} \sup_{x\in \dbR^{d_1}} \Big|\sum_{\vec\k \in \cS_{N}}  \cV_0(\vec i, \vec j, \vec\k; x)\cK(\vec i, \vec\k; \vec t)\Big|.
\eea

\begin{rem}Motivated from the Taylor expansion \eqref{Taylor1}, the step-N Volterra signature should have the following form in the space $\bigoplus_{n=0}^{N}(\dbR^{d_1+1})^{\otimes n}$:
	\bea\label{sig}
	\sum_{n=0}^{N} \sum_{\vec i\in \cI_n, \vec j\in \cJ_n, \vec\k\in \cS_n} \Big(\int_{\dbT_n} \cK(\vec i,  \vec \k; \vec t~\!) \circ dB^{\vec j}_{\vec t}\Big)(e_{j_1}\otimes \cdots \otimes e_{j_n} ),
	\eea
where $\{e_j\}_{j=0,1,\cdots,d_1}$ denotes the canonical basis of $\dbR^{d_1+1}$.
	Below, we shall focus on the expectation of the  Volterra signature at any step.
\end{rem}

To facilitate the cubature method in the next subsection, we shall rewrite \reff{Taylor1} further slightly.  Note that, for fixed $N$, the mapping $(\vec i, \vec j, \vec\k)\in {\bigcup_{n\le N}} \cI_n\times \cJ_{n,N}\times \cS_n \to \cV_0(\vec i, \vec j, \vec \k; \cd)$  (as a function of $x$) is not one to one, so we may combine the terms with the same $\cV_0(\vec i, \vec j, \vec \k; \cd)$.  That is, we may rewrite \reff{Taylor1} as:
\bea
\label{Taylor1new}
\left.\ba{c}
\dis G(X_T) = G(x) +  \sum_{\phi\in \dbV_{N}} \phi(x) \G^\phi_{N}  + R_N,\q  \mbox{where}\ms\\
\dis \dbV_{N} := \bigcup_{n\le N}\Big\{\cV_0(\vec i, \vec j, \vec \k; \cd): (\vec i, \vec j, \vec\k)\in \cI_n\times \cJ_{n,N}\times \cS_n\Big\} \subset C(\dbR^{d_1}; \dbR),\\
\dis  \G^\phi_{N}:= \sum_{n=1}^N \sum_{\vec i\in \cI_n, \vec j\in \cJ_{n,N}, \vec \k\in \cS_n} \1_{\{\cV_0(\vec i, \vec j, \vec \k; \cd)=\phi\}}\int_{\dbT_n} \cK(\vec i,  \vec \k; \vec t~\!) \circ dB^{\vec j}_{\vec t}.
\ea\right.
\eea
Since it requires rather complicated notations to characterize $\phi\in \dbV_{N}$ precisely in the general case, we leave it to the special cases we will actually compute numerically.

\subsection{The cubature formula}
We now extend the cubature formula for Brownian motion in \cite{LV03} to the Volterra setting, especially for the Taylor expansion \reff{Taylor1new}. From now on, we set $\O:= C([0, T]; \dbR^d)$ the canonical space, $B$ the canonical process, and thus $\dbP$ is the Wiener measure so that $B$ is a $\dbP$-Brownian motion.  
For some $W \ge 1$, $L\ge 1$, we introduce a discrete probability measure $Q$ on $\O$: for some constants $a_{k,l} =(a^1_{k,l},\cds, a^d_{k,l})\in \dbR^{d}$, $k=1,\cds, W$, $l=1,\cds, L$,
\bea
\label{QK}
\left.\ba{c}
\dis Q:= \sum_{k=1}^{2W} \l_k \d_{\o_k},\q\mbox{where $\d_\cd$ denotes the Dirac measure},\q \l_k>0, \q \sum_{k=1}^{2W} \l_k = 1,\\  
\dis   \l_{W+k} = \l_k,\q \o_{W+k}=-\o_k,\q  k=1,\cds, W,\ms\\
\dis  \o_{k,0} = 0,\q \o_{k,t} = \o_{k, s_{l-1}} + {a_{k,l}\over \sqrt{T}} [t-s_{l-1}], ~ t\in (s_{l-1}, s_{l}], ~ s_l := {l\over L}T,~ 0=1,\cds, L.
\ea\right.
\eea
Here, the second line implies that $Q$ is symmetric, since Brownian motion is symmetric. Also, it is ok to consider non-uniform partition $0=s_0<\cds<s_L=T$. Recall \reff{cKcV}, for each piecewise linear $\o = (\o^1, \cds, \o^{d_1})$ as in \reff{QK}, $\vec j \in \cJ_n$, and $\f: \dbT_n \to \dbR$,  denote 
\bea
\label{pats2}
 \int_{\dbT_n} \f(\vec t~\!) ~d\o^{\vec j}_{\vec t} ~:=~  \int_{\dbT_n}  \f(\vec t~\!)  ~d\o^{j_n}_{t_n}  \cds  d\o^{j_1}_{t_1},\q\mbox{where}\q \o^0_t:= t.
 \eea
 Then we have
 \bea
 \label{EQK}
 \dbE^{Q}\Big[\int_{\dbT_n} \cK(\vec i,  \vec \k; \vec t~\!) \circ dB^{\vec j}_{\vec t}\Big] = \sum_{k=1}^{2W}\l_k \int_{\dbT_n} \cK(\vec i,  \vec \k; \vec t~\!) ~d(\o_k)^{\vec j}_{\vec t}.
 \eea

\begin{defn}\label{def: cubature}
Let $N \ge 1$, $W\ge 1$. We say $Q$ defined in \reff{QK} is an N-Volterra cubature formula on $[0, T]$ if, recalling $\dbE=\dbE^\dbP$,   
\bea
\label{cubature}
 \dbE^{Q}[ \G^\phi_{N}] = \dbE[ \G^\phi_{N}] \q\mbox{for all}\q \phi\in \dbV_{N},\qq\mbox{and hence}\q \dbE^{Q}[I_{N}] = \dbE[I_{N}].
\eea
\end{defn}

Recall our goal \reff{volterra sde}-\reff{Y0EG}. Our main idea is the following approximation:
\bea
\label{Y0EG2}
&\dis Y_0:=\dbE\big[G(X_T)\big] \approx Y^{Q}_0 := \sum_{k=1}^{2W} \l_k G(X_{T}(\o_k)),\q \mbox{where} \\
\label{Xomega}
&\dis X^{i}_{t}(\o) =x_i + \sum_{j=0}^d \int_0^t K_i(t,r) V^i_j(X_{r}(\o))  d\o^j_{r},~ i=1,\cds, d_1.
\eea

We now have the main result of this section, whose proof is postponed to Appendix.
\begin{thm}
\label{thm-Qest}
Under {\bf (H0)} and {\bf (H(N+3))}, we have: recalling \reff{AN} and \reff{cV0},
\bea
\label{Qest1}
\left.\ba{c}
\dis |Y_0 - Y^{Q}_0| \le C_N \Big[ A_{N+1}  \big(1+C_{Q}^{N-1}\big)T^{N+1\over 2} + A_{N+2} \big(1+C_{Q}^{N-2}\big) T^{N+2\over 2}+A_{N+3} T^{N+3\over 2} \Big],\ms\\
\dis \mbox{where}\q C_{Q} :=\max_{1\le k\le W, 1\le j \le d, 1\le l\le L} |a_{k,l}^j|.
\ea\right.
\eea
In particular, if each $K_i$ is rescalable, in the sense that there exists an $\a_i \in [0, \infty)$ (not necessarily the same as $H_i-{1\over 2}$) such that 
\bea
\label{rescale}
K_i(ct, cr) = c^{\a_i} K(t,r),\q\mbox{for all}~ 0\le r< t.
\eea
Then all the $a_{k,l}$ and hence $C_Q$ are independent of $T$.
\end{thm}

 \subsection{A simplification of the cubature formula}
Due to the symmetric properties of Brownian motion and $Q$, we may simplify the requirement \reff{cubature}. Recall \reff{cKcV} and abuse the notation, for $\vec j\in \cJ_n$ we denote
\bea
\label{cNj}
\cN_\a(\vec j):=\{l\in \{1,\cds, n\}: j_l =\a\},\q \a = 0,\cds, d.
\eea
 
\begin{lem}
\label{lem-odd}
Let {\bf (H0)} hold and $\vec j\in \cJ_n$ be such that $|\cN_\a(\vec j)|$ is odd for some $\a=1,\cds, d$, in particular if $\|\vec j\|$ is odd, then
\bea
\label{odd}
\dbE\Big[\int_{\dbT_n} \cK(\vec i,  \vec \k; \vec t~\!) \circ dB^{\vec j}_{\vec t}\Big] =0=\dbE^{Q}\Big[\int_{\dbT_n} \cK(\vec i,  \vec \k; \vec t~\!) \circ dB^{\vec j}_{\vec t}\Big].
\eea
\end{lem}
\proof One may easily derive the first equality from \reff{EdB} by induction on $n$. The second equality follows directly from the symmetric properties of $Q$.
\qed

Note further that, when $\vec j = (0,\cds, 0)\in \cJ_n$, we have $d B^{\vec j}_{\vec t} = d \o^{\vec j}_{\vec t} = dt_n\cds dt_1$. This, together with Lemma \ref{lem-odd},  implies the following result immediately. 
\begin{thm}
\label{thm-cubature}
Let {\bf (H0)} and {\bf (H(N+3))} hold and denote 
\bea
\label{hV0}
\left.\ba{lll}
\dis \bar \cJ_{n, N} := \{ \vec j\in \cJ_{n,N}\backslash \{(0,\cds,0)\}: ~|\cN_\a(\vec j)|~\mbox{is even for all $\a=1,\cds, n$}\big\},\ms\\
\dis \bar \dbV_N := \bigcup_{n\le N}\Big\{\cV_0(\vec i, \vec j, \vec \k; \cd): (\vec i, \vec j, \vec\k)\in \cI_n\times \bar \cJ_{n,N}\times \cS_n\Big\} \subset \dbV_N,\\
\dis \bar \G^\phi_{N}:= \sum_{n=1}^N \sum_{\vec i\in \cI_n, \vec j\in \bar\cJ_{n,N}, \vec \k\in \cS_n} \1_{\{\cV_0(\vec i, \vec j, \vec \k; \cd)=\phi\}}\int_{\dbT_n} \cK(\vec i,  \vec \k; \vec t~\!) \circ dB^{\vec j}_{\vec t}.
\ea\right.
\eea
Then $Q$ satisfies \reff{cubature} if and only if: 
\bea
\label{cubature2}
 \dbE^{Q}[ \bar \G^\phi_{N}] = \dbE[ \bar\G^\phi_{N}],\q\mbox{for all}~ \phi \in \bar \dbV_{N}.
\eea
\end{thm}

When $N$ is odd, note that $\bar\cJ_{N, N}=\emptyset$, so we will get the cubature formula for free at the $N$-th order. Therefore, we shall always consider odd $N$.
\begin{eg}
\label{eg-barcJ}
(i) In the case $N=3$, obviously we have 
\bea
\label{cJ3}
\bar \cJ_{1, 3} = \bar \cJ_{3,3} = \emptyset,\q \bar \cJ_{2,3} = \{(j, j): 1\le j\le d\}.
\eea

(ii) In the case $N=5$, we have:
\bea
\label{cJ5}
\left.\ba{c}
\dis \bar \cJ_{1, 5} =  \bar \cJ_{5,5}=\emptyset,\q \bar \cJ_{2,5} = \big\{(j, j): 1\le j\le d\},\q \bar \cJ_{3,5} =\big\{(j, j, 0), (j, 0, j), (0, j, j)\big\},\ms\\
\dis \bar \cJ_{4,5} = \big\{(j, j, j, j), (j, j, \tilde j, \tilde j), (j, \tilde j, j, \tilde j), (j, \tilde j, \tilde j, j): 1\le j\neq \tilde j\le d\big\}.
\ea\right.
\eea
\end{eg}

\section{The cubature formula: the multiple period case}
\label{sect-cubatureM}
\setcounter{equation}{0}
In this case we consider general $T$, and we use the setting in \S\ref{sect-Taylor}, in particular $\d:= {T\over M}$. 

\subsection{The cubature formula on each subinterval $[T_m, T_{m+1}]$}
Recall \reff{Taylor}.  Note that in \reff{Taylor1}  $\cV_0(\vec i, \vec j, \vec\k; x)$ and $\int_{\dbT_n} \cK(\vec i, \vec\k; \vec t) \circ dB_{\vec t}^{\vec j}$ are separated and the cubature measure $ Q$ is determined only by $\int_{\dbT_n} \cK(\vec i, \vec\k; \vec t) \circ dB_{\vec t}^{\vec j}$. In \reff{Taylor}, however, 
\beaa
\cV(\vec i,\vec j, \vec\k;\vec t, \Th^{\vec t}_{T_m},\Th^{[t^m_0,T]}_{T_m}) =  \prod_{\alpha = 1}^{n} \pa^{\vec \k, \alpha}_{\vec i}  V^{i_\a}_{j_\a} (\Th^{t_\a}_{T_m}) \cK_+(\vec i,  \vec \k; \vec t) ~\big \la \pa_{\vec i}^{\vec \k}u(t^m_0, \Th^{[t^m_0,T]}_{T_m}) , \vec \cK_0(\vec i,  \vec\k; \vec t)  \big\ra
\eeaa
and we are not able to move the term $\prod_{\alpha = 1}^{n} \pa^{\vec \k, \alpha}_{\vec i}  V^{i_\a}_{j_\a} (\Th^{t_\a}_{T_m}) $ outside of the stochastic integral, which prohibits us from constructing a desirable $ Q_m$ to match the conditional expectations of $I^m_N$: $\dbE^{ Q_m}_m[I^m_N]=\dbE_m[I^m_N]$. In light of \reff{RNest}, we shall instead content ourselves with
\bea
\label{Qerror}
\Big|\dbE^{ Q_m}_m[I^m_N]-\dbE_m[I^m_N]\Big| \le C \d^{N+1\over 2}.
\eea
We shall remark though, in general conditional expectations are only defined in a.s. sense, which requires specifying the probability on $\cF_{T_m}$. However, here we will construct $Q_m$ only on the paths on $[T_m, T_{m+1}]$. For this purpose, we interpret the conditional expectations in a pathwise sense, as we explain in the remark below, so that \reff{Qerror} could make sense.

\begin{rem}
\label{rem-cond}
Under our conditions, one can easily see that $\dbE_m[I^m_N] = v_m(\Th^{[T_m, T]}_{T_m})$, for a deterministic function $v_m \in C(\dbX_{T_m})$. Similarly, for the $Q_m$ we are going to construct, we will interpret it as a regular conditional probability distribution and thus we also have the structure $\dbE^{ Q_m}_m[I^m_N] = \tilde v_m(\Th^{[T_m, T]}_{T_m})$ for a deterministic function $\tilde v_m \in C(\dbX_{T_m})$. Then by \reff{Qerror} we actually mean a stronger result:
\bea
\label{pointwiseest}
|\tilde v_m(\th) - v_m(\th)|\le C\d^{N+1\over 2},\q\mbox{for all}\q \th\in \dbX_{T_m}.
\eea
We refer to \cite[Chapter 9]{Zhang} for more details of the pathwise stochastic analysis. In this paper, since our main focus is the approximation, to avoid introducing further complicated notations, we abuse the notation slightly and write them as conditional expectations. 
\end{rem}

From now on, we shall assume $K_i$ is sufficiently smooth in $(t, r)$. Then, recalling \reff{volterra sde theta},  the mapping $s\in [T_m, T]\to \Th^s_{T_m}$ is smooth: for any $\a\ge 0$,
\bea
\label{paTh}
{\pa^\a\over \pa s^\a} \Th^{i, s}_{T_m} = \sum_{j=0}^d \int_0^{T_m} {\pa^\a\over \pa s^\a} K_i(s, r) V^i_j(X_r) \circ dB^j_r.
\eea
Then the mapping $\vec t\in \dbT^m_n \to \check \cV(\vec t ):= \cV(\vec i,\vec j, \vec\k;\vec t, \Th^{\vec t}_{T_m},\Th^{[t^m_0,T]}_{T_m})$ is also smooth. Our idea is to introduce further the Taylor expansion of $\check \cV(\vec t)$ at $\vec T_m := (T_m,\cds, T_m)$. For any $\vec \a = (\a_1, \cds, \a_n)$ with $\a_l = 0, 1,\cds$, denote $\|\vec \a\| := \sum_{l=1}^n \a_l$ and $\vec \a! := \prod_{l=1}^n (\a_l!)$. Then by the standard Taylor expansion formula we have, for any $k\ge 0$,
\bea
\label{checkRk}
\left.\ba{c}
\dis \check \cV(\vec t~\!) = \sum_{\|\vec \a\|\le k} {1\over \vec\a!} {\pa^{\|\vec\a\|}\over \pa {\vec t}^{\vec \a}} \check \cV(\vec T_m) \prod_{l=1}^n (t_l-T_m)^{\a_l} + \check R_k(\vec t),~ \mbox{where}~ {\pa^{\|\vec\a\|}\over \pa {\vec t}^{\vec \a} }:= {\pa^{\|\vec\a\|}\over  \pa t_1^{\a_1}\cds \pa t_n^{\a_n}}\\
\dis \mbox{and}\q |\check R_k(\vec t)| \le C_{n,k} \sup_{\|\vec \a\|=k+1} \sup_{s_l\in [T_m, t_l], l=1,\cds, n}  \Big|{\pa^{k+1}\over \pa \vec t^{\vec\a}} \check \cV(s_1, \cds, s_n)\Big| \d^{k+1}.
\ea\right.
\eea
We now extend Theorem \ref{thm-Taylor}. Recall Lemma \ref{lem-odd} and Theorem \ref{thm-cubature}.
 
\begin{thm}
\label{thm-Taylorm}
Let $N$ be odd and  {\bf (H0-${N-1\over 2}$)}, {\bf (H(N+3))} hold.  Let $\check R^m_N$ be determined by: 
 \bea
\label{checkRmN}
&&\dis\qq\qq\qq\qq\qq I^m_N = \check I^m_N  +  \check R^m_N,\q\mbox{where}\\
&&\dis \check I^m_N := u(t^m_0, \Th^{[t^m_0,T]}_{T_m})  + \sum_{n=1}^N  \sum_{\vec i\in \cI_n,  \vec \k\in \cS_n}\Big[\sum_{\vec j\in \cJ_{n,N}\backslash \bar\cJ_{n,N}}\int_{\dbT^m_n} \cV(\vec i,\vec j, \vec\k;\vec t, \Th^{\vec t}_{T_m},\Th^{[t^m_0,T]}_{T_m}) \circ dB^{\vec j}_{\vec t} \nonumber\\
&&\dis   + \sum_{\vec j\in \bar\cJ_{n,N}}\sum_{\vec \a: \|\vec \a\|\le {N-\|\vec j\|\over 2}} {1\over \vec\a!}  {\pa^{\|\vec\a\|}\over \pa \vec t^{\vec \a}}\cV(\vec i,\vec j, \vec\k;\vec T_m, \Th^{\vec T_m}_{T_m},\Th^{[t^m_0,T]}_{T_m}) \int_{\dbT^m_n} \prod_{l=1}^n (t_l-T_m)^{\a_l} \circ dB^{\vec j}_{\vec t}\Big].\nonumber
\eea
Then, there exists a constant $C^m_N$, which depends on $N$, $H_i$, and the upper bounds of $V^i_j$ and their derivatives up to the order $N+2$, such that
\bea
\label{checkRmNest} 
 &\dis \dbE_m[|\check R^m_N|^2] \le C^m_N e^{C^m_NT}\d^{N+1}.
 \eea
 \end{thm}
\proof For each $\vec j\in \bar \cJ_{n,N}$, noting that $\|\vec j\|$ is even and $N$ is odd, set $k := {N-\|\vec j\|-1\over 2}$. Using the notations in \reff{checkRk}, one can see that ${\pa^{k+1}\over \pa \vec t^{\vec\a}} \check \cV$ involves the derivatives of $V^i_j$ and $u(t^m_0,\cd)$ up to the order $n + k+1$, and the derivatives of $K_i$ up to the order $k+1$. Note that $2\le \|\vec j\|\le N-1$, then
\beaa
&\dis n + k+1 \le \|\vec j\| + {N-\|\vec j\|-1\over 2} + 1 = {N+\|\vec j\| +1\over 2} + 1\le N+1,\\
&\dis k+1 = {N-\|\vec j\|-1\over 2} + 1\le  {N-3\over 2} + 1\le {N-1\over 2}.
\eeaa
Recall Proposition \ref{prop-ureg} for the bounds of the derivatives of $u$. Now following the arguments in Theorem \ref{thm-Taylor}, one can easily see that, for some appropriate constant $C$ depending on the parameters specified in this theorem, 
\beaa
\dbE\Big[ \Big|\int_{\dbT^m_n} \check R_k(\vec t) \circ dB^{\vec j}_{\vec t}\Big|^2\Big] \le C  e^{C^m_NT}\d^{2(k+1) + \|\vec j\|}= C e^{C^m_NT}\d^{N+1}.
\eeaa
This implies \reff{checkRmNest} immediately.
 \qed
 
 We next introduce $Q_m$ as in \reff{QK}, but on paths on $[T_m, T_{m+1}]$: 
 \bea
\label{QKm}
\left.\ba{c}
\dis Q_m:= \sum_{k=1}^{2W} \l_k \d_{\o_k},\q \l_k>0,\q \sum_{k=1}^{2W} \l_k = 1,\\
\dis \l_{W+k} = \l_k,\q \o_{W+k}=-\o_k,\q  1\le k\le W,\\
\dis  \o_{k,0} = T_m,~ \o_{k,t} = \o_{k, s_{l-1}} \!\!+ {a_{k,l}\over \sqrt{\d}} [t-s_{l-1}], ~ t\in (s_{l-1}, s_{l}], ~ s_l := T_m+{l\over L}\d,~ 0\le l\le L.
\ea\right.
\eea
 Recall Remark \ref{rem-cond}.
 \begin{defn}\label{def: cubaturem}
Let $N \ge 1$ be odd and fix $m$. We say $Q_m$  defined in \reff{QKm} is a modified $N$-Volterra cubature formula on $[T_m, T_{m+1}]$ if, for all $n\le N$, $\vec j\in \bar\cJ_{n,N}$, and $\|\vec\a\| \le {N-\|\vec j\|-1\over 2}$,  
\bea
\label{cubaturem}
 \dbE^{Q_m}_m\Big[ \int_{\dbT^m_n} \prod_{l=1}^n (t_l-T_m)^{\a_l} \circ dB^{\vec j}_{\vec t}\Big] = \dbE_m\Big[ \int_{\dbT^m_n} \prod_{l=1}^n (t_l-T_m)^{\a_l} \circ dB^{\vec j}_{\vec t}\Big].
  \eea
\end{defn}
 
  \begin{rem}
 \label{rem-CQm}
 (i) The equations in \reff{cubaturem} corresponding to $\|\vec \a\|=0$ exactly characterize the cubature measures for  Brownian motions. That is, our modified $N$-Volterra cubature formula is a cubature formula for the standard one, but not vice versa in general. In the case $N=3$, however, as we will see in Example \ref{eg-QmK} (i) below, the two are equivalent.
 
(ii) The kernel $\prod_{l=1}^n (t_l-T_m)^{\a_l}$ in \reff{cubaturem} is rescalable in the sense of \reff{rescale}. Then by the same arguments as in Theorem \ref{thm-Qest},  $\dis  C_{Q_m} := \max_{1\le k\le W, 1\le j \le d, 1\le l\le L} |a_{k,l}^j|$ is independent of $\d$ (or $M$). Indeed, as in  \reff{omegalinear2} below, $Q_m$ is  a modified N-Volterra cubature formula on $[T_m, T_{m+1}]$ if and only if the following $Q^*_N$ is a  modified N-Volterra cubature formula on $[0,1]$:
 \bea
\label{Q*}
\left.\ba{c}
\dis Q^*_N:= \sum_{k=1}^{2W} \l_k \d_{\o_k},\q \l_k>0,\q \sum_{k=1}^{2W} \l_k = 1,\\
\dis \l_{W+k} = \l_k,\q \o_{W+k}=-\o_k,\q  1\le k\le W,\\
\dis  \o_{k,0} = 0,~ \o_{k,t} = \o_{k, s_{l-1}} + a_{k,l} [t-s^*_{l-1}], ~ t\in (s^*_{l-1}, s^*_{l}], ~ s^*_l := {l\over L},~ 0\le l\le L.
\ea\right.
\eea
We emphasize that $Q^*_N$ is universal in the sense that it depends only on $N$,  the dimensions, and our construction of the cubature measure, but does not depend on $T, M$ or even $K$. In particular, $\dis C_{Q^*_N} := \max_{1\le k\le W, 1\le j \le d, 1\le l\le L} |a_{k,l}^j|$ is independent of $T$, $M$, or $\d$, and $C_{Q_m} = C_{Q^*_N}$.
 \end{rem}
 
 \begin{thm}
\label{thm-Qestm}
Let $N$ be odd and  {\bf (H0-${N-1\over 2}$)}, {\bf (H(N+3))} hold, and $Q_m$ be as in Definition \ref{def: cubaturem}. Then, for the  $C^m_N$ as in Theorem \ref{thm-Taylorm} and $C_{Q^*_N}$ in Remark \ref{rem-CQm}, we have
\bea
\label{Qestm}
\left.\ba{lll}
\dis \Big| \dbE^{Q_m}_m\big[u(T_{m+1}, \Th^{[T_{m+1}, T]}_{T_{m+1}})\big] - \dbE_m\big[u(T_{m+1}, \Th^{[T_{m+1}, T]}_{T_{m+1}})\big]\Big|  \le C^m_N \big(1+ C_{Q^*_N}^{N-1}\big)e^{C^m_NT} \d^{N+1\over 2} \ms\\
\dis\qq + C_N\Big[ A^m_{N+1}  \big(1+ C_{Q^*_N}^{N-1}\big) \d^{N+1\over 2} + A^m_{N+2}  \big(1+ C_{Q^*_N}^{N-2}\big)\d^{N+2\over 2} +A^m_{N+3} \d^{N+3\over 2}\Big].
\ea\right.
\eea
\end{thm}
\no This proof is also postponed to Appendix.

\begin{eg}
\label{eg-QmK}
(i) When $N=3$, recall \reff{cJ3} and note that ${N-\|\vec j\|-1\over 2} =0$ for $\vec j=(j, j)\in \bar\cJ_{2,3}$, we see that \reff{cubaturem} is equivalent to the cubature formula for  standard Brownian motions:
\bea
\label{QmN3}
 \dbE^{Q_m}_m\Big[ \int_{\dbT^m_2}  \circ dB^{(j,j)}_{\vec t}\Big] = \dbE_m\Big[ \int_{\dbT^m_2}    \circ dB^{(j,j)}_{\vec t}\Big] = {\d\over 2}.
  \eea

(ii) In the case $N=5$, recall \reff{cJ5} and note that ${N-\|\vec j\|-1\over 2} =1$ for $\vec j=(j, j)\in \bar\cJ_{2,5}$ and ${N-\|\vec j\|-1\over 2} =0$ for $\vec j\in \bar\cJ_{3,5}\cup \bar \cJ_{4, 5}$, then \reff{cubaturem} is equivalent to: for $1\le j\neq \tilde j \le d$, $l=1,2$,
\bea
\label{QmN5}
\left.\ba{lll}
\dis \dbE^{ Q_m}_m\Big[\int_{\dbT^m_2}  \circ dB^{(j, j)}_{\vec t}\Big]  =\dbE_m\Big[\int_{\dbT^m_2}  \circ dB^{(j, j)}_{\vec t}\Big] = {\d\over 2},\\
\dis \dbE^{Q_m}_m\Big[\int_{\dbT^m_2} (t_l-T_m) \circ dB^{(j, j)}_{\vec t}\Big]= \dbE_m\Big[\int_{\dbT^m_2} (t_l-T_m) \circ dB^{(j, j)}_{\vec t}\Big]={\d^2\over 4},\\
\dis \dbE^{Q_m}_m\Big[\int_{\dbT^m_3} \circ dB^{\vec j}_{\vec t}\Big]=  \dbE_m\Big[\int_{\dbT^m_3} \circ dB^{\vec j}_{\vec t}\Big] = {\d^2\over 4},\q \vec j= ( j, j, 0), (0, j, j),\\
\dis \dbE^{Q_m}_m\Big[\int_{\dbT^m_3} \circ dB^{( j, 0, j)}_{\vec t}\Big]=  \dbE_m\Big[\int_{\dbT^m_3} \circ dB^{(j, 0, j)}_{\vec t}\Big] = 0,\\
\dis \dbE^{Q_m}_m\Big[\int_{\dbT^m_4} \circ dB^{\vec j}_{\vec t}\Big]=   \dbE_m\Big[\int_{\dbT^m_4} \circ dB^{\vec j}_{\vec t}\Big] = {\d^2\over 8},\q \vec j=(j,j, j, j), (j, j, \tilde j, \tilde j),\\
\dis \dbE^{Q_m}_m\Big[\int_{\dbT^m_4} \circ dB^{\vec j}_{\vec t}\Big]=   \dbE_m\Big[\int_{\dbT^m_4} \circ dB^{\vec j}_{\vec t}\Big] = 0, \q \vec j=(j,\tilde j, j, \tilde j), (j, \tilde j, \tilde j,  j).
\ea\right.
\eea
\end{eg}

\subsection{The cubature formula on the whole interval $[0, T]$}
Recall $Q_m$ is defined on $C([T_m, T_{m+1}]; \dbR^{d_1})$. We shall now compose all the $Q_m$:
\bea
\label{Qmultiple}
Q := Q_0 \otimes \cds \otimes Q_{M-1}.
\eea
Here $\otimes$ refers to independent composition. Then $Q$ is a probability measure on $\O = C([0, T]; \dbR^{d_1})$.  Similarly let $\dbP_m$ denote the Wiener measure on $C([T_m, T_{m+1}]; \dbR^{d_1})$, then $\dbP = \dbP_0 \otimes\cds\otimes \dbP_{M-1}$.  The following result extends \cite[Theorem 3.3]{LV03} to our setting.

\begin{thm}
\label{thm-QestMultiple}
Let $N$ be odd and  {\bf (H0-${N-1\over 2}$)}, {\bf (H(N+4))} hold, and $Q$ is defined by \reff{Qmultiple} with each $Q_m$ as in Definition \ref{def: cubaturem}. Then, for the $C^m_N$ in Theorem \ref{thm-Qestm} and $C_{Q^*_N}$ in Remark \ref{rem-CQm}, we have
\bea
\label{QestMultiple}
\left.\ba{lll}
\dis \Big| \dbE^{Q}\big[G(X_T)\big] - \dbE\big[G(X_T)\big]\Big|\le \sum_{m=0}^{M-1}\Big[C^m_N  \big(1+ C_{Q^*_N}^{N-1}\big) e^{C^m_NT}\d^{N+1\over 2}\\
\dis + C_N \big[ A^m_{N+1}  \big(1+ C_{Q^*_N}^{N-1}\big) \d^{N+1\over 2}+ A^m_{N+2} \d^{N+2\over 2} \big(1+ [C_{Q^*_N}\sqrt{\d}]^{N-2}\big)+A^m_{N+3} \d^{N+3\over 2}\big]\Big].
\ea\right.
\eea
Moreover, for a possibly larger $C_N$ which may depend on the bounds of the derivatives of $V^i_j$ up to the order $N+4$, and the $C_{Q^*_N}$, but not on $M$, we have\footnote{ The constant $e^{C_N T}$ below is due to the estimate for the derivatives of $u$ in Proposition \ref{prop-ureg}. If one can improve this estimate, under certain technical conditions, then one can replace $e^{C_N T}$ with the new bound for the derivatives of $u$ up to the order $N+3$. This comment is valid for the estimates in \reff{checkRmNest} , \reff{Qestm}, \reff{QestMultiple} as well.}
\bea
\label{QestMultiple2}
\dis \Big| \dbE^{Q}\big[G(X_T)\big] - \dbE\big[G(X_T)\big]\Big|\le C_N M e^{C_NT}\d^{N+1\over 2} = C_Ne^{C_NT}{T^{N+1\over 2}\over M^{N-1\over 2}},
\eea
which converges to $0$ as $M\to \infty$.
\end{thm}
\proof Note that, recalling \reff{u}.
\beaa
&&\dis \Big| \dbE^{Q}\big[G(X_T)\big] - \dbE\big[G(X_T)\big]\Big| = \Big| \dbE^{Q_0\otimes \cds \otimes Q_{M-1}}\big[G(X_T)\big] - \dbE^{\dbP_0\otimes \cds\otimes \dbP_{M-1}}\big[G(X_T)\big]\Big|\\
&&\dis \le \sum_{m=0}^{M-1} \Big| \dbE^{Q_0\otimes \cds \otimes Q_m \otimes \dbP_{m+1}\otimes \cds\otimes \dbP_{M-1}}\big[G(X_T)\big] - \dbE^{Q_0\otimes \cds Q_{m-1} \otimes \dbP_{m}\otimes \cds\otimes \dbP_{M-1}}\big[G(X_T)\big]\Big|\\
&&\dis \le \sum_{m=0}^{M-1} \Big| \dbE^{Q_0\otimes \cds \otimes Q_m} \big[u(T_{m+1}, \Th^{[T_{m+1}, T]}_{T_{m+1}})\big] - \dbE^{Q_0\otimes \cds Q_{m-1} \otimes \dbP_m}\big[u(T_{m+1}, \Th^{[T_{m+1}, T]}_{T_{m+1}})\big]\Big|\\
&&\dis \le \sum_{m=0}^{M-1} \dbE^{Q_0\otimes \cds \otimes Q_{m-1}} \Big[\Big| \dbE^{Q_m}_m \big[u(T_{m+1}, \Th^{[T_{m+1}, T]}_{T_{m+1}})\big] - \dbE^{\dbP_m}_m\big[u(T_{m+1}, \Th^{[T_{m+1}, T]}_{T_{m+1}})\big]\Big|\Big].
\eeaa
Recall Remark \ref{rem-cond} and note that $\dbE^{\dbP_m}_m = \dbE_m$, then by \reff{pointwiseest} we see that \reff{QestMultiple} follows directly from \reff{Qestm}. Finally, by \reff{cKcV}, \reff{cJkm}, and Proposition \ref{prop-ureg} we obtain \reff{QestMultiple2}. 
\qed

\begin{rem}
\label{rem-cubaturem}
(i) Compared to Theorem \ref{thm-Qest}, the above Theorem \ref{thm-QestMultiple} allows us to deal with large $T$. Moreover, 
compared with the $Q$ in \reff{cubature2},  it is easier to construct the cubature measure $Q_m$ in \reff{cubaturem} and the $Q$ in \reff{Qmultiple}. The price to pay, however, is that \reff{cubaturem} requires higher regularity of $K_i$ in order to have the desired convergence rate.

(ii) Provided sufficient regularity {\bf (H0-${N-1\over 2}$)} on $K_i$ $($and {\bf (H(N+4))} on $V^i_j$ and $G$ $)$, we have the convergence and its rate in \reff{QestMultiple2} as $M\to \infty$, which is very desirable in theory. However, by \reff{QKm} and \reff{Qmultiple}, we see that each $Q_m$ will involve $2W$ paths, and thus the independent composition $Q$ will involve $(2W)^M$ paths. Therefore,  practically we still don't want to make $M$ too large, which in turn means that $T$ cannot be too large.  We note that the same difficulty arises in the  Brownian setting, and there have been various ideas on improving the efficiency, for example the recombination schemes in \cite{HayakawaOberhauserLyons, LL12}. It will be very interesting to explore these ideas and see if they can be extended to the Volterra setting.

(iii) Note that the choice of $Q_m$ is not unique. In particular, \reff{cubaturem} involves certain number of equations. To make it solvable, we need to allow for sufficient number of parameters $\l_k$, $a_{k, l}$, $1\le k\le W$, $0\le l\le L-1$. As mentioned in (ii), the complexity of our cubature algorithm increases dramatically for large $W$, but is much less sensitive to the value of $L$. So whenever possible, we would prefer a small $W$ while allowing for a reasonably large $L$. We shall remark that, when the dimension $d_1$ is large, typically we need a large $W$. This consideration is not serious for the one period case, which however requires $T$ to be small.

(iv) Clearly we have a better rate for a larger $N$ (again, provided sufficient regularity). However, a larger $N$ implies more equations in \reff{cubaturem}, which in turn requires larger values of $W$ and/or $L$. In the meantime, a larger $N$ implies a larger $C_N$ in \reff{QestMultiple2}. So the algorithm may not be always more efficient for a larger $N$. 
\end{rem}

\section{A one dimensional model}
\label{sect-1d}
\setcounter{equation}{0}
In this section we focus on the following model with $d=d_1=1$: 
\bea
\label{X1}
X_t = x_0 + \int_0^t K(t, r) V(X_r) \circ dB_r,\q K(t, r) = (t-r)^{H-{1\over 2}},
\eea
where the Brownian motion $B$ is one dimensional and the Hurst parameter $H>{1\over 2}$.  We investigate a few cases in details and compute the desired $Q$. We shall illustrate the efficiency of our algorithm in these cases by several numerical examples in \S\ref{sect-numerical} below.

Note that in this case $V^1_0\equiv 0$, then there is no need to consider $j=0$. So, for $n\le N$, by abusing the notation we may view $\cI_n = \{(1,\cds, 1)\}$, $\cJ_n = \cJ_{n,N} = \{(1,\cds, 1)\}$. We may  omit $\vec i = (1,\cds, 1)$, $\vec j=(1,\cds,1)$ inside $\cK$ and $\cV$ in \reff{cKcV}.

\subsection{The multiple period case with order $N=3$}
\label{sect-N3-mul}

We shall  construct $Q_m$ as in \reff{QKm}  for a fixed $m$. In the numerical examples in \S\ref{sect-numerical}, we may simply compose these $Q_m$ independently as in \reff{Qmultiple}. 

In this case \reff{QmN3} consists of only one equation:
\bea
\label{Qm3-1}
\dis \dbE^{ Q_m}\Big[\int_{T_m}^{T_{m+1}}\int_{T_m}^{t_1} dB_{t_2} \circ dB_{t_1}\Big] = {\d\over 2}.
\eea
To construct $Q_m$, we set $W=1$  and $L=1$  in \reff{QKm}:
\bea
\label{Qm3-2}
\l_1 =  {1\over 2},\q d \o_{1, t} = {a\over \sqrt{\d}} dt, ~t\in [T_m, T_{m+1}].
\eea
Then \reff{Qm3-1} becomes
\beaa
{\d\over 2} = \int_{T_m}^{T_{m+1}}\int_{T_m}^{t_1} d\o_{t_2}  d\o_{t_1} = {a^2\over \d} \int_{T_m}^{T_{m+1}}\int_{T_m}^{t_1} dt_2  dt_1 = {a^2 \d\over 2}.
\eeaa
Thus
\bea
\label{Qm3-3}
a = 1 \q\mbox{and hence}\q \o_{1, t} = {t- T_m\over \sqrt{\d}}, ~t\in [T_m, T_{m+1}].
\eea
We remark that the above computation does not involve $H$, in fact, as we saw in Example \ref{eg-QmK} (i), the cubature measure in this case coincides with that of standard Brownian motion. However, in order to have the desired error estimate, by Theorem \ref{thm-Qestm} we need $H>{3\over 2}$.

\subsection{The multiple period case with order $N=5$} \label{sect-N5-mul}
While we may apply \reff{QmN5} directly, in this one dimensional case actually we may simplify the problem further. Note that the corresponding term which requires the further expansion \reff{checkRk} or \reff{checkRmN} is: recalling \reff{expansion2} and abusing the notation $\check \cV(\vec t)$,
\beaa
&\dis \sum_{\vec\k \in \cS_2}\int_{\dbT^m_2} \cV(\vec\k;\vec t, \Th^{\vec t}_{T_m},\Th^{[t^m_0,T]}_{T_m}) \circ dB^{(1,1)}_{\vec t}  = \int_{\dbT^m_2} \check \cV(\vec t)  \circ dB^{(1,1)}_{\vec t}\\
&\dis\mbox{where}\q \check\cV(\vec t) := \big\la\pa_{\bx\bx}u(t^m_0,\Th_{T_m}^{[t^m_0,T]}),(K_{t_1}^{[t^m_0,T]}, K_{t_2}^{[t^m_0,T]}) \big\ra~ V(\Th^{t_1}_{T_m})V(\Th^{t_2}_{T_m}) \\
&\dis\qq\qq +\big\la\pa_{\bx}u(t^m_0,\Th^{[t^m_0,T]}_{T_m}), K_{t_1}^{[t^m_0,T]} \big \ra ~K(t_1,t_2) \pa_{x} V(\Th^{t_1}_{T_m}) V(\Th^{t_2}_{T_m}).
\eeaa
For $N=5$, we shall assume {\bf (H0-2)}, namely  $H>{5\over 2}$, then $K(T_m, T_m) = \pa_t K(T_m, T_m) = \pa_r K(T_m, T_m) =0$. Thus
\beaa
&& \pa_{t_1} \check\cV(\vec t)\big|_{\vec t=(T_m, T_m)}=\big\la\pa_{\bx\bx}u(t^m_0,\Th_{T_m}^{[t^m_0,T]}),(\pa_r K_{r}^{[t^m_0,T]}\big|_{r=T_m}, K_{T_m}^{[t^m_0,T]}) \big\ra~ V(X_{T_m})V(X_{T_m}) \\
 &&\qq\qq +\big\la\pa_{\bx\bx}u(t^m_0,\Th_{T_m}^{[t^m_0,T]}),(K_{T_m}^{[t^m_0,T]}, K_{T_m}^{[t^m_0,T]}) \big\ra~ V'(X_{T_m})V(X_{T_m}) \pa_s\Th^s_{T_m}\big|_{s=T_m},\\
 && \pa_{t_2} \check\cV(\vec t)\big|_{\vec t=(T_m, T_m)}=\big\la\pa_{\bx\bx}u(t^m_0,\Th_{T_m}^{[t^m_0,T]}),(K_{T_m}^{[t^m_0,T]}, \pa_r K_{r}^{[t^m_0,T]}\big|_{r=T_m}) \big\ra~ V(X_{T_m})V(X_{T_m}) \\
 &&\qq\qq +\big\la\pa_{\bx\bx}u(t^m_0,\Th_{T_m}^{[t^m_0,T]}),(K_{T_m}^{[t^m_0,T]}, K_{T_m}^{[t^m_0,T]}) \big\ra~ V'(X_{T_m})V(X_{T_m})\pa_s\Th^s_{T_m}\big|_{s=T_m}.
  \eeaa
 Note that $ \bla \pa_{xx} u, (\eta_1, \eta_2)\bra = \bla \pa_{xx} u, (\eta_2, \eta_1)\bra$, then  
 \beaa
 \pa_{t_1} \check\cV(\vec t)\big|_{\vec t=(T_m, T_m)} = \pa_{t_2} \check\cV(\vec t)\big|_{\vec t=(T_m, T_m)}.
 \eeaa
 This leads to the following  expansion:
 \beaa
 &&\dis \int_{\dbT^m_2} \check \cV(\vec t)  \circ dB^{(1,1)}_{\vec t} = \check \cV(T_m, T_m)  \int_{\dbT^m_2}   dB^{(1,1)}_{\vec t} \\
 &&\dis\qq +  \pa_{t_1} \check\cV(\vec t)\big|_{\vec t=(T_m, T_m)}\int_{\dbT^m_2} [(t_1-T_m) + (t_2-T_m)] \circ dB^{(1,1)}_{\vec t} + \check R(\vec t),
 \eeaa
where $\check R(\vec t)$ satisfies the desired estimate. Consequently, we may merge the 2 equations in the second line of \reff{QmN5} into one equation, in the same spirit of \reff{Taylor1new}, by considering only their sum. Therefore, in this case \reff{QmN5} reduces to three equations:
\bea
\label{Qm5-1}
\left.\ba{lll}
\dis \dbE^{ Q_m}\Big[\int_{T_m}^{T_{m+1}}\int_{T_m}^{t_1} dB_{t_2} \circ dB_{t_1}\Big]  = {\d\over 2},\ms\\
\dis \dbE^{ Q_m}\Big[\int_{T_m}^{T_{m+1}}\int_{T_m}^{t_1}\big[(t_1-T_m) + (t_2-T_m)\big] dB_{t_2} \circ dB_{t_1}\Big] = {\d^2\over 2},\ms\\
\dis \dbE^{ Q_m}\Big[\int_{T_m}^{T_{m+1}}\int_{T_m}^{t_1} \int_{T_m}^{t_2}\int_{T_m}^{t_3} dB_{t_4} \circ dB_{t_3} \circ dB_{t_2} \circ dB_{t_1}\Big] = {\d^2\over 8}.
\ea\right.
\eea

To construct $Q_m$, we set $W=2$   and $L=1$ in \reff{QKm}:
\bea
\label{Qm5-2}
\l_1 + \l_2 = {1\over 2},\qq d \o_{k, t} = {a_{k}\over \sqrt{\d}} dt, ~k=1,2.
\eea
Note that, in light of Remark \ref{rem-cubaturem}, we would prefer a small $W$. However, if we set $W=1$ here, the cubature measure $Q_m$ does not exist for any value $L$. 
By straightforward calculation we see that \reff{Qm5-1} becomes
\bea
\label{Qm5-3}
\sum_{k=1}^2 \l_k a_k^2 \d= {\d\over 2};\qq \sum_{k=1}^2 \l_k a_k^2 \d^2 = {\d^2\over 2};\qq \sum_{k=1}^2 \l_k {a_k^4\over 12} \d^2= {\d^2\over 8}.
\eea
In particular, the first two equations coincide, and we obtain
\bea
\label{Qm5-4}
a_1^2 = 4\big[1+ \sqrt{2\l_2\over \l_1}\big],\qq a_2^2 = 4\big[1- \sqrt{2\l_1\over \l_2}\big].
\eea
This requires $\l_1\le {1\over 6}$ so that $\sqrt{2\l_1\over \l_2} \le 1$. Then for any $0\le \l_1\le {1\over 6}$, we would obtain a solution by \reff{Qm5-4}. One particular solution is:
\bea
\label{Qm5-5}
\l_1 := {1\over 6},\q \l_2 := {1\over 3},\q a_1 = \sqrt{3},\q a_2 :=0.
\eea
We remark that, in this case $-\o_2 = \o_2=0$, so we actually have a total of three paths, instead of four paths: by abusing the notation $\l_2$,
\bea
\label{Qm5-6}
\l_1 = {1\over 6},\q d \o_{1,t} = \sqrt{3\over \d} dt,\q \l_2 = {2\over 3},\q {d\over dt} \o_{2,t}  =0,\q \l_3 = {1\over 6},\q d \o_{3,t} = -\sqrt{3\over \d}dt.
\eea

\subsection{The one period case with order $N=3$}\label{sect-N3}
Recall \reff{hV0} and \reff{cJ3}, in particular we shall only consider $n=2$. Then one can verify straightforwardly that, 
\bea
\label{hV1}
\left.\ba{c}
\dis \cS_2 = \{ (0, 0); (0, 1)\}, \q \bar \dbV_3 = \Big\{\cV_0( \vec\k; \cd): \vec\k\in \cS_2\Big\} = \big\{G'' V^2,~ G' V' V\big\},\\
\dis \bar \G^{G''V^2}_3 =  \int_{\dbT_2} \cK((0,0); \vec t) \circ dB^{(1,1)}_{\vec t} = \int_0^T \int_0^{t_1} [(T-t_1)(T-t_2)]^{H-{1\over 2}}  dB_{t_2}\circ dB_{t_1},\\
\dis \bar \G^{G'V'V}_3 =  \int_{\dbT_2} \cK((0,1); \vec t) \circ dB^{(1,1)}_{\vec t} = \int_0^T \int_0^{t_1} [(T-t_1)(t_1-t_2)]^{H-{1\over 2}}  dB_{t_2}\circ dB_{t_1}.
\ea\right.
\eea
By \reff{EdB}, one can easily compute that
\bea
\label{EQK1}
\dis \dbE\big[\bar \G^{G''V^2}_3 \big] = {T^{2H}\over 4H},\q \dbE\Big[\bar \G^{G'V'V}_3 \Big]  = 0.
\eea

To construct $Q$, we set $W=1$ and $L=2$ in \reff{QK} :
\bea
\label{QK1}
\l_1  = {1\over 2},\q d\o_{1, t} = \Big[{a_1\over \sqrt{T}}  \1_{[0, {T\over 2}]}(t)  + {a_2\over \sqrt{T}}  \1_{({T\over 2}, T]}(t)\Big]dt.
\eea
 For notational simplicity, we introduce 
\bea
\label{Htilde}
H_- := H-{1\over 2}>0,\q  H_+ := H+{1\over 2} >1.
\eea
Then one may compute:
\bea
\label{cKdo1}
\left.\ba{c}
\dis  \dbE^Q\big[\bar \G^{G''V^2}_3 \big] =\!\!\! \int_0^T \!\!\! \int_0^{t_1}\!\!\! [(T-t_1)(T-t_2)]^{H_-}  d\o_{1,t_2} d\o_{1, t_1} = {T^{2H}\over 2H_+^2} \big[(1-{1\over 2^{H_+}}) a_1  + {a_2\over 2^{H_+}}\big]^2,\ms\\
\dis \dbE^Q\big[\bar \G^{G'V'V}_3 \big]= \!\!\!\int_0^T \!\!\! \int_0^{t_1} [(T-t_1)(t_1-t_2)]^{H_-}  d\o_{1,t_2} d\o_{1, t_1}  =  T^{2H}\big[c_1 a_1^2 +  c_2 a_1a_2 + c_3 a_2^2\big],
\ea\right.
\eea
where
\bea
\label{c123}
\left.\ba{c}
\dis  c_1 := \int_0^{1\over 2} \int_0^{t_1} [(1-t_1) (t_1-t_2)]^{H_-} dt_2 dt_1 = {1\over H_+}\int_0^{1\over 2} (1-t)^{H_-} t^{H_+} dt,\ms\\
\dis  c_2 := \int_{1\over 2}^1\! \int_0^{1\over 2} [(1-t_1) (t_1-t_2)]^{H_-} dt_2 dt_1 ={1\over H_+} \int_{1\over 2}^1\! (1-t)^{H_-} \big[t^{H_+}  - (t-{1\over 2})^{H_+}\big] dt,\ms\\
\dis  c_3 := \int_{1\over 2}^1 \int_{1\over 2}^{t_1} [(1-t_1) (t_1-t_2)]^{H_-} dt_2 dt_1 = {1\over H_+} \int_{1\over 2}^1 (1-t)^{H_-} (t-{1\over 2})^{H_+} dt.
\ea\right.
\eea
Combining \reff{EQK1} and \reff{cKdo1}, we obtain from \reff{cubature} that
\bea
\label{Equation1}
 \Big[(2^{H_+}-1) a_1  + a_2\Big]^2 =  {H_+^2 2^{2H_+}\over 2H},\q c_1 a_1^2 +  c_2 a_1a_2 + c_3 a_2^2=0.
\eea
First by the second equation we obtain  (we may use the other one as well):
\bea
\label{c4}
{a_2\over a_1} = c_4:= {-c_2 + \sqrt{c_2^2 - 4 c_1c_3}\over 2c_3}.
\eea
Plug this into the first equation in \reff{Equation1} we obtain one solution:
\bea
\label{Solution1}
a_1 = {H_+ 2^{H_+} \over \sqrt{2H} [2^{H_+}+c_4-1]},\q a_2 = {H_+ 2^{H_+} c_4\over \sqrt{2H} [2^{H_+}+c_4-1]}.
\eea

\begin{eg}
\label{eg1}
Set $H={3\over 2}$ at above, then one may compute straightforwardly that
\beaa
&\dis c_1 = {5\over 384},\q c_2 = {10\over 384},\q c_3 =  {1\over 384},\q c_4 = -5+2\sqrt{5},\\
&\dis a_1 ={\sqrt{5} + 1\over \sqrt{3}},\q a_2 = {5-3\sqrt{5}\over \sqrt{3}},
\eeaa
and we obtain $Q$ through \reff{QK1}.
\end{eg}

\subsection{The one period case with order $N=5$}\label{sect-N5}
Recall \reff{hV0} and \reff{cJ5}, in particular we shall only consider $n=2$ and $n=4$. Clearly $\bar \dbV_3 \subset \bar \dbV_5$, for the $\bar \dbV_3$ in \reff{hV1}. Moreover, again omitting $\vec i = (1,\cds, 1), \vec j=(1,\cds, 1)$,
\beaa
\bar\dbV_5 \backslash \bar \dbV_3 = \Big\{\cV_0(\vec \k; \cd): \vec\k \in \cS_4\Big\} 
\eeaa
Recall  \reff{Htilde} and the Gamma function $\G(\a,\b):= \int_0^1 (1-t)^{\a-1} t^{\b-1} dt$, and denote by $\f^{(k)}$ the $k$-th derivative of $\f$. We then have the following result.

\begin{lem}
\label{lem-hV2}
For the above model, we have $\bar\dbV_5 \backslash \bar \dbV_3  = \{\phi_\a\}_{1\le \a\le 7}$, where
 \bea
\label{hV2}
\left.\ba{c}
\dis \phi_1 = G^{(4)} V^4,\q \phi_2 = G^{(3)} V' V^3,\q \phi_3 = G'' V'' V^3,\q \phi_4 = G'' (V')^2 V^2, \\
\dis \phi_5 = G' V^{(3)} V^3,\q \phi_6 = G' V'' V' V^2,\q  \phi_7 = G' (V')^3 V.
\ea\right.
\eea
Moreover, denote $\g^{\phi_\a}_4 := \dbE[\bar\G^{\phi_\a}_{4}]$, we have
\bea
\label{EQK2}
\left.\ba{c}
\dis \g^{\phi_1} = {T^{4H}\over 32H^2},\q  \g^{\phi_2}_{4} = {T^{4H}\over 4H}\G(2H, H_+),\ms\\
\dis  \g^{\phi_3}_{4} =\g^{\phi_4}_{4} = {T^{4H}\over 4H} \G(2H, 2H_+),\q  \g^{\phi_5}_{4}=\g^{\phi_6}_{4}=\g^{\phi_7}_{4}=0.
\ea\right.
\eea
\end{lem}
\proof Note that $\cS_4 =  \bigcup_{\a=1}^7 \cS_{4, \a}$,
\bea
\label{S4}
\left.\ba{lll}
\dis  \cS_{4,1} := \big\{(0,0,0,0)\big\},\q \cS_{4, 5} := \big\{(0, 1,1,1)\},\q \cS_{4, 7} := \big\{(0, 1,2,3)\}, \\
 \dis\cS_{4, 2} := \big\{(0,0,0,1), (0,0, 1, 0), (0,1, 0, 0), (0,0, 0, 2), (0,0,2, 0), (0,0,0, 3)\big\},\\
 \dis\cS_{4,3} := \big\{(0,0, 1,1), (0,1,0,1), (0,1,1,0), (0,0, 2,2)\big\},\\
 \dis\cS_{4,4}:= \big\{(0,0, 1, 2), (0,0, 2,1), (0,0, 1, 3), (0,0, 2, 3), (0,1, 0, 2), (0, 1, 2, 0), (0,1,0, 3)\big\},\\
 \dis\cS_{4, 6} := \big\{(0, 1,1,2), (0,1,2,1), (0,1,2,2), (0,1,1, 3)\}.
\ea\right.
\eea
By \reff{cV0} one can check that $\cV_0(\vec\k; \cd) = \phi_\a$, $\forall \vec\k\in \cS_{4, \a}$, $\a=1,\cds, 7$. Then, by \reff{EdB},
\beaa
&& \dbE\Big[\int_{\dbT_4} \cK(\vec\k; \vec t~\!)\circ dB^{(1,1,1,1)}_{\vec t}\Big] = {1\over 4} \int_0^T \int_0^{t_1}  \cK(\vec\k; (t_1, t_1, t_3, t_3)) dt_3 dt_1\\
&&={T^{4H}\over 4} \int_0^1 \int_0^{t_1}  \big[(1-t_1) (t_{\k_2}-t_1) (t_{\k_3}-t_3)(t_{\k_4}-t_3)\big]^{H_-} dt_3 dt_1 .
\eeaa
Now the expectations in \reff{EQK2} follow from straightforward computation.
\qed

We next construct a desired $Q$. Note that  \reff{cubature2} consists of  $2$ equations for $n=2$ and $7$ equations for $n=4$. To allow for sufficient flexibility, we set $W=2$ and $L=4$ in \reff{QK}: 
\bea
\label{QK2}
\left.\ba{c}
\dis \l_1+\l_2 = {1\over 2},\q d \o_{k, t} = \sum_{l=1}^4 {a_{k, l}\over \sqrt{T}} \1_{ [s_{l-1}, s_l)}(t)dt, \q  k=1,2.
\ea\right.
\eea
Similar to \reff{cKdo1}, the following result is obvious.
\begin{lem}
\label{lem-cub2}
For the above model \reff{QK2}, we have, for $k=1,2$,
 \bea
\label{cub2}
\left.\ba{lll}
\dis \int_{\dbT_2} \cK((0,0); (t_1, t_2)) d (\o_k)^{(1,1)}_{\vec t} ={1\over 2H_+^2 T }\Big[\sum_{l=1}^4 [s_l^{H_+}-s_{l-1}^{H_+}] a_{k,l}\Big]^2,\\
\dis \int_{\dbT_2} \cK((0,1); (t_1, t_2)) d (\o_k)^{(1,1)}_{\vec t} ={1\over T}\sum_{1\le l_2\le  l_1\le  4} c(l_1, l_2) a_{k,l_1} a_{k, l_2},\\
\dis \int_{\dbT_4} \cK(\vec\k; \vec t~\!) d (\o_k)^{(1,1,1,1)}_{\vec t} = {1\over T^2}\sum_{1\le l_4\le l_3 \le l_2\le l_1\le 4} c(\vec\k, \vec l) a_{k, l_1} a_{k, l_2} a_{k, l_3} a_{k,l_4}
\ea\right.
\eea
where, 
\bea
\label{ckl}
\left.\ba{lll}
\dis c(l, l) := \!\int_{s_{l-1}}^{s_l}\!\! \int_{s_{l-1}}^{t_1} \!\!\!\! [(T-t_1) (t_1-t_2)]^{H_-}dt_2 dt_1={1\over H_+} \int_{s_{l-1}}^{s_l} \!\!\!\! (T-t)^{H_-} (t-s_{l-1})^{H_+} dt,\ms\\
\dis c(l_1, l_2) := \int_{s_{l_1-1}}^{s_{l_1}}\int_{s_{l_2-1}}^{s_{l_2}} [(T-t_1) (t_1-t_2)]^{H_-}dt_2 dt_1\ms\\
\dis\qq ={1\over H_+} \int_{s_{l_1-1}}^{s_{l_1}} (T-t)^{H_-} [(t-s_{l_2-1})^{H_+}-(t-s_{l_2})^{H_+}] dt,\q l_1>l_2,\ms\\
\dis c(\vec\k, \vec l) :=  \int_{s_{l_1-1}}^{s_{l_1}} \int_{s_{l_2-1}}^{s_{l_2}\wedge t_1} \int_{s_{l_3-1}}^{s_{l_3}\wedge t_2}\int_{s_{l_4-1}}^{s_{l_4}\wedge t_3}   \cK(\vec\k; \vec t~\!)dt_4dt_3dt_2dt_1. 
\ea\right.
\eea
\end{lem}

Combine \reff{EQK1}, \reff{EQK2}, and \reff{cub2}, we have the following rsult.
\begin{thm}
\label{thm-equation2} \reff{cubature} is equivalent to the following equations:
\bea
\label{Equation2}
\left.\ba{c}
\dis \l_1 + \l_2 = {1\over 2},\\
\dis {1\over H_+^2 T}\sum_{k=1}^2\l_k\Big[\sum_{l=1}^4 [s_l^{H_+}-s_{l-1}^{H_+}] a_{k,l}\Big]^2={T^{2H}\over 4H},\ms\\
\dis \sum_{k=1}^2\l_k\sum_{1\le l_2\le  l_1\le  4} c(l_1, l_2) a_{k,l_1} a_{k, l_2} =0,\ms\\
\dis {2\over 4! T^2} \sum_{\vec l \in \{1,\cds, 4\}^4} \sum_{\vec\k \in \cS_{4,\a}} \sum_{k=1}^2 \l_k c(\vec\k, \vec l~\!)  a_{k, l_1} a_{k, l_2} a_{k, l_3} a_{k,l_4} = \g^{\phi_\a}_4,\q \a=1,\cds, 7.
\ea\right.
\eea
\end{thm}

We remark that \reff{Equation2} consists of $10$ equations, with $10$ unkowns: $\l_k, a_{k,l}$, $l=1,2,3,4$, $k=1,2$.  Since these equations are nonlinear, in particular they involve $4$-th order polynomials of $a_{k,l}$, in general we are not able to derive explicit solutions as in \reff{Solution1}. Indeed, even the existence of solutions is not automatically guaranteed, and in that case we can actually increase $W$ and/or $L$ in \reff{QK2} to allow for more unknowns. Nevertheless, we can solve \reff{Equation2} numerically, and our numerical examples in the next section show that the numerical solutions of \reff{Equation2} serve for our purpose well.

\section{A fractional stochastic volatility model}
\label{sect-rough}
\setcounter{equation}{0}
Consider a financial market where $S_t$ denotes the underlying asset price and $U_t$ is the volatility process:
\bea
\label{Heston}
\left.\ba{c}
\dis S_t= S_0 + \int_0^t b_1(r,S_r,U_r)dr+\int_0^t \si_1(r,S_r,U_r)\circ dB^1_r,\\
\dis U_t=U_0+\int_0^t K(t,r)b_2(r,U_r)ds+\int_0^t K(t,r)\si_2(r,U_r)\circ dB^2_r.
\ea\right.
\eea
Here $B^1, B^2$ are correlated Brownian motions with constant correlation $\rho\in [-1, 1]$,  $K(t,r) = (t-r)^{H-{1\over 2}}$ with Hurst parameter $H>{1\over 2}$. Assume for simplicity the interest rate is $0$, and our goal is to  compute the option price $\dbE[G(S_T)]$.

Note that \reff{Heston} involves the time variable $t$, so we are in the situation with $d=2$ and $d_1 =3$ in \reff{volterra sde}. Indeed, denote $X_t  = (X^0_t, X^1_t, X^2_t) := (t, S_t, U_t)$, then we have 
\bea
\label{Heston-setting}
\left.\ba{c}
\dis x = (0, S_0, U_0),\q K_0=K_1 = 1, ~ K_2 = K,\\
\dis V^0 = (1, 0, 0),\q  V^1 = (b_1, \si_1, 0),\q V^2 = (b_2, 0, \si_2).
\ea\right.
\eea
Here, for notational simplicity, we use indices $(0,1,2)$ instead of $(1, 2, 3)$ for $X$. We shall emphasize that, although $B^1, B^2$ are correlated here, the Taylor expansions \reff{expansionN} and \reff{Taylor} will remain the same, but the expectations in Lemma \ref{lem-EdB} need to be modified in an obvious way. In particular, \reff{odd} will hold true only when $\|\vec j\|$ is odd. Therefore, in this section we modify the $\bar \cJ_{n,N}$ in \reff{hV0}: still denoted as $\bar \cJ_{n,N}$ by abusing the notation,
\bea
\label{barcJnN2}
\bar \cJ_{n, N} := \{ \vec j\in \cJ_{n,N}\backslash \{(0,\cds,0)\}: ~\|\vec j\|~\mbox{is even}\big\}.
\eea
Then all the results in the previous sections remain true. 
Alternatively, we may express $B^1, B^2$ as linear combinations of independent Brownian motions. However, this will make  $V^1, V^2$ more complicated and does not really simplify the analysis below.

\subsection{The multiple period case with $N=3$}
\label{sect-N3-2d-mul}
By \reff{barcJnN2}, we see that
\beaa
\bar\cJ_{1,3} = \bar \cJ_{3,3} = \emptyset,\q \bar\cJ_{2,3} = \{(j_1, j_2): j_1, j_2>0\}.
\eeaa
Thus the cubature measure $Q_m$ should satisfy: for $\vec j \in \bar\cJ_{2,3}$,
\bea
\label{2dQm}
\dbE^{Q_m}_m\big[\int_{\dbT^m_n} \circ dB^{\vec j}_{\vec t}\big] = \dbE_m\big[\int_{\dbT^m_n} \circ dB^{\vec j}_{\vec t}\big]. 
\eea 
This is the same as the Brownian motion case. More precisely,
\bea
\label{2dQm2}
\left.\ba{lll}
\dis \dbE^{Q_m}_m\big[\int_{\dbT^m_2} \circ dB^{(1,1)}_{\vec t}\big] = \dbE^{Q_m}_m\big[\int_{\dbT^m_2} \circ dB^{(2,2)}_{\vec t}\big] = {\d\over 2},\\
\dis \dbE^{Q_m}_m\big[\int_{\dbT^m_2} \circ dB^{(1,2)}_{\vec t}\big] = \dbE^{Q_m}_m\big[\int_{\dbT^m_2} \circ dB^{(2,1)}_{\vec t}\big] = {\rho \d\over 2}.
\ea\right. 
\eea 
To construct $Q_m$, we set $W=2$   and $L=1$ in \reff{QKm}: noting that $\o$ is $2$-dimensional,
\bea
\label{2dQm3}
\l_1 = \l_2 = {1\over 4},\q d \o_{k, t} = d(\o^1_{k,t}, \o^2_{k,t})= \big( {a^1_k\over \sqrt{\d}},~ {a^2_k\over \sqrt{\d}}\big)dt,  ~k=1,2.
\eea
Plug these into \reff{2dQm2}, we have
\bea
\label{2dQm4}
{\d\over 4}[|a_1^1|^2 + |a_2^1|^2] =  {\d\over 4}[|a_1^2|^2 + |a_2^2|^2] = {\d\over 2},\q {\d\over 4}[a_1^1 a_1^2 + a_2^1a_2^2] =   {\rho\d\over 2}
\eea
One can easily solve the above equations:
\bea
\label{2dQm5}
\left.\ba{c}
\dis a^1_1 = \sqrt{2} \sin(\th_1),\q a^1_2 = \sqrt{2} \cos(\th_1),\q a^2_1 = \sqrt{2} \sin(\th_2),\q a^2_2 = \sqrt{2} \cos(\th_2),\ms\\
\dis \mbox{for any $\th_1, \th_2$ satisfying $\cos(\th_1 - \th_2) = \rho$}.
\ea\right.
\eea

\subsection{The one period case with $N=3$}
\label{sect-N3-2d}
We first note that, due to the multiple dimensionality here, the system corresponding to \reff{cubature2} will be pretty large, especially when $N=5$ in the next subsction.  However, since $(X^0_t, X^1_t) = (t, S_t)$ are not of Volterra type, the system can be simplified significantly. Furthermore, we shall modify the cubature method slightly as follows. 

\begin{rem}
\label{rem-groupI}
Recall that in \reff{Taylor1new} and \reff{hV0} we group the terms with the same $\cV_0(\vec i, \vec j, \vec \k; \cd)$. Note that the mapping $(\vec i, \vec j, \vec\k) \to \int_{\dbT_n} \cK(\vec i,  \vec \k; \vec t~\!) \circ dB^{\vec j}_{\vec t}$  is also not one to one, and since $(X^0_t, X^1_t) = (t, S_t)$, many terms do not appear in the Taylor expansion (or say the corresponding $\cV_0(\vec i, \vec j, \vec \k; \cd)=0$). It turns out that it will be more convenient to group the terms based on  $\int_{\dbT_n} \cK(\vec i,  \vec \k; \vec t~\!) \circ dB^{\vec j}_{\vec t}$ for this model, as we will do in this and the next subsections. To be precise, let $\tilde \dbV^0_N$ denote the terms $\int_{\dbT_n} \cK(\vec i,  \vec \k; \vec t~\!) \circ dB^{\vec j}_{\vec t}$ with $(\vec i, \vec j, \vec\k)\in \cI_n\times (\cJ_{n,N}\backslash \{(0,\cds,0)\})\times \cS_n$ appearing in the Taylor expansion, and  $\D\tilde \dbV^0_N:= \tilde \dbV^0_N\backslash\tilde \dbV^0_{N-1}$. We emphasize that, unlike in \reff{Taylor1new}, the elements here are not $\cV_0$. We then modify Definition \ref{def: cubature} by replacing \reff{cubature} with:
\bea
\label{cubatureI}
\dbE^Q\Big[\int_{\dbT_n} \cK(\vec i,  \vec \k; \vec t~\!) \circ dB^{\vec j}_{\vec t}\Big]=\dbE\Big[\int_{\dbT_n} \cK(\vec i,  \vec \k; \vec t~\!) \circ dB^{\vec j}_{\vec t}\Big], 
\q\mbox{for all the terms in $\tilde \dbV^0_N$.}
\eea 
 We remark that, if $\int_{\dbT_n} \cK(\vec i',  \vec \k'; \vec t~\!) \circ dB^{\vec j'}_{\vec t}=\int_{\dbT_n} \cK(\vec i,  \vec \k; \vec t~\!) \circ dB^{\vec j}_{\vec t}$ (as random variables), then automatically we have $\int_{\dbT_n} \cK(\vec i',  \vec \k'; \vec t~\!) ~d\o^{\vec j'}_{\vec t}=\int_{\dbT_n} \cK(\vec i,  \vec \k; \vec t~\!) ~d\o^{\vec j}_{\vec t}$ for all $\o$ in \reff{QK}.
\end{rem}

Recall \reff{barcJnN2} that Lemma \ref{lem-odd} remains true when $\|\vec j\|$ is odd,  we shall only find $\D\tilde\dbV^0_{2}$. Instead of applying the Taylor expansion \reff{Taylor1} directly on \reff{Heston}, we first expand $\dbE[G(S_T)]$ as in \reff{Taylor1}.  Indeed, 
note that $\Th^s_t=(t, S_t, \Th^{s,2}_t)$, where
\bea
\label{HestonTh}
\Th^{s,2}_t=U_0+\int_0^t K(s,r)b_2(r,U_r)ds+\int_0^t K(s,r)\si_2(r,U_r)\circ dB^2_r.
\eea
Then, applying the chain rule we have: denoting $\tilde b_1 := G' b_1$ and $\tilde \si_1 := G' \si_1$,
\bea
\label{chainrule}
\left.\ba{c}
\dis G(S_T) = G(S_0) + \int_0^T \Big[\tilde b_1(\Th^{t_1}_{t_1})dt_1+\tilde \si_1(\Th^{t_1}_{t_1})\circ dB^1_{t_1}\Big]\\
\dis = R_{\neq 2} +\int_{\dbT_1}\tilde b_1(\Th^{t_1}_{t_1})dt_1+ \int_{\dbT_2} \sum_{i=1}^2 \pa_i \tilde \si_1(\Th^{t_1}_{t_2}) \si_i(\Th^{t_2}_{t_2}) K_i(t_1, t_2)\circ dB^i_{t_2}\circ dB^1_{t_1},
\ea\right.
\eea
Here $R_{\neq n}$  is a generic term whose order is not equal to $n$. Then,  recalling Remark \ref{rem-groupI},   by \reff{Heston-setting} we can easily obtain
\bea
\label{2dhV2}
\left.\ba{c}
\dis \D\tilde\dbV^0_{2} := \{\tilde \G_2^{(1,1)}, ~\tilde \G_2^{(1,2)} \},\q\mbox{where}\q  \tilde \G_2^{(1,i)}:=  \int_{\dbT_2}K_i(t_1, t_2)d B^i_{t_2} \circ d B^1_{t_1},\ms\\
\dis \dbE\big[\tilde\G^{(1,i)}_{2}\big]={1\over 2}\int_{\dbT_1} K_i(t_1, t_1)d\la B^i, B^1\ra_{t_1}= \left\{\ba{lll} {T\over 2},\q i=1,\\ 0, \q  i=2.\ea\right.
\ea\right.
\eea

To construct $Q$, we set $W=1$ and $L=1$  in \reff{QK} : $d \o_{1,t} = ( {a_1\over \sqrt{T}},~ {a_2\over \sqrt{T}})dt$. Then
\beaa
\dbE^Q[\tilde \G_2^{(1,1)}]= {|a_1|^2\over T} \int_{\dbT_2} \!\! d\vec t = {T\over 2} |a_1|^2,\q \dbE^Q[\tilde \G_2^{(1,2)}]= {a_1a_2\over T} \int_{\dbT_2} \!\! K_2(t_1, t_2) d\vec t= {T^{H_+}\over H_+(H_++1)}a_1a_2.
\eeaa
 By \reff{2dhV2} we have
\beaa
{T\over 2}|a_1|^2  ={T\over 2},\qq {T^{H_+}\over H_+(H_++1)}a_1a_2 =0,\q\mbox{implying:}\q  a_1 = 1,\q a_2 =0.
\eeaa

We remark that this solution is independent of $\rho$. In fact, numerical results (which are not reported in the paper) show that this does not provide a good approximation, even when $T$ is small. So we shall move to the order $N=5$ in the next subsection, although it becomes much more involved to find the cubature measure.

\subsection{The one period case with $N=5$}
\label{sect-N5-2d}
In this case, besides the $\D\tilde\dbV^0_2$ in \reff{2dhV2}, we also need  $\D\tilde\dbV^0_4$. For this purpose, we need the $2$-nd order expansion of $\tilde b_1(\Th^{t_1}_{t_1})$ and the $3$-rd order expansion of $\tilde \si_1(\Th^{t_1}_{t_1})$. These derivations are straightforward, but rather tedious. We thus postpone them to Appendix and turn to numerical examples first.

\section{Numerical examples}
\label{sect-numerical}
\setcounter{equation}{0}

\subsection{The algorithm}
\label{sect-algorithm}

Our numerical algorithm consists of the following five steps. We are illustrating only the algorithm in \S\ref{sect-N5}. The algorithms in the other subsections, especially   those in \S\ref{sect-rough}, need to be modified slightly in the obvious manner.

{\it Step 1.} Compute $\dbE\Big[\int_{\dbT_n} \cK(\vec i,  \vec \k; \vec t~\!) \circ dB^{\vec j}_{\vec t}\Big]$  for each $n$ and $(\vec i, \vec j, \vec \k)$  by using \reff{EdB}, and then compute $ \dbE[ \bar\G^\phi_{N}]$ in \reff{cubature2} for each $\phi\in \bar \dbV_{N}$. 

{\it Step 2.} Compute $\dbE^{Q}\Big[\int_{\dbT_n} \cK(\vec i,  \vec \k; \vec t~\!) \circ dB^{\vec j}_{\vec t}\Big]$ by \reff{EQK} and then compute $ \dbE^Q[ \bar\G^\phi_{N}]$ in \reff{cubature2} for each $\phi\in \bar \dbV_{N}$.

{\it Step 3.} Establish the equations \reff{cubature2} with unknowns $\l_k, a_{k, l}$, $k=1,\cds, W$, $l=1,\cds, L$, from \reff{QK}, then solve these equations  to obtain a desired $Q$. We may in general use numerical methods to solve these (polynomial) equations, when explicit solutions are not available, see Remark \ref{rem-optimization} below.

{\it Step 4.} For each $\o_k$ obtained in Step 3, solve the (deterministic) ODE \reff{Xomega}, by discretizing $[0, T]$ equally into $D$ pieces. That is, denoting $h:= {T\over D}$ and $\o_{-h}:=0$, 
\bea
\label{XD}
X^{D,i}_{lh}(\o) =x_i + \sum_{j=0}^d \sum_{\a=0}^{l-1} K_i(lh, \a h) V^i_j(X^D_{\a h}(\o))  {\o^j_{(\a+1)h}-\o^j_{(\a-1) h}\over 2},\q l=0,\cds, D,
\eea
For convenience, we typically set $D$ as a multiple of $ML$  for the $M, L$ in \reff{cubature}.

{\it Step 5.} We obtain the approximation by \reff{Y0EG2}: $\dis Y^{Q}_0 \approx Y^{Q,D}_0:=\sum_{k=1}^{2W} \l_k G(X^D_{T}(\o_k))$.

We note that, provided the conditions in Remark \ref{rem-Taylor}, our algorithm is deterministic and is much more efficient than the probabilistic methods, e.g. the Euler scheme in Zhang \cite{ZhangXicheng}.

\begin{rem}
\label{rem-optimization} 
(i) In Subsections \ref{sect-N3-mul}, \ref{sect-N5-mul}, \ref{sect-N3}, \ref{sect-N3-2d-mul}, \ref{sect-N3-2d}, we have obtained the cubature paths, then we can move to Step 4 directly.

(ii) We remark that Steps 1-3 depend only on the model, more specifically only on $K_i$, but not on the specific forms of $V^i_j$ or $G$. So, given the model, we may compute the desired $Q$ offline, and then for each $V$ and $G$, we only need to complete Steps 4 and 5. 

(iii) To illustrate the idea for Step 3, we consider the equations in \reff{Equation2}. We shall use the steepest decent method to minimize the following weighted sum:
\bea
\label{optimization}
&\dis \inf_{\l_k\ge 0, a_{k,l}\in \dbR, k=1,2, l= 1,2,3,4}  \Big[ \b_1 \big|\l_1 + \l_2 - {1\over 2}\big|^2 + \b_2\big| \sum_{k=1}^2\l_k\sum_{1\le l_2\le  l_1\le  4} c(l_1, l_2) a_{k,l_1} a_{k, l_2}\big|^2\nonumber\\
&\dis + \b_3\big| {1\over H_+^2 }\sum_{k=1}^2\l_k\big[\sum_{l=1}^4 [s_l^{H_+}-s_{l-1}^{H_+}] a_{k,l}\big]^2-{T^{2H}\over 4H}\big|^2 \nonumber\\
&\dis\q + \sum_{\a=1}^7 \bar \b_\a \big|
 {2\over 4!} \sum_{\vec l \in \{1,\cds, 4\}^4} \sum_{\vec\k \in \dbS_{4,\a}} \sum_{k=1}^2 \l_k c(\vec\k, \vec l)  a_{k, l_1} a_{k, l_2} a_{k, l_3} a_{k,l_4} - \G^{\phi_\a}_4\big|^2\Big],\nonumber
\eea
where $\b_i, \bar \b_\a>0$, $i=1,2,3$ and $\a=1,\cds, 7$, are some appropriate weights.

 (iv) For (iii) above, we may replace it with any efficient solver for equations \reff{Equation2}.

(v) While our algorithm is more sensitive to $W$ than to $L$, as pointed out in Remark \ref{rem-cubaturem} (iii), a large $L$ will increase the difficulty to solve the equations like \reff{Equation2}. However, since this can be done offline, the impact of $L$ is less serious. 
\end{rem}

\begin{rem}
\label{rem-D}
In this paper we focus on the impacts of $M$ and $N$, but do not analyze rigorously the impact of  $D$ (or $h$) in Step 4, which can be chosen much larger than $ML$. We shall only comment on it heuristically in this remark.

(i) First, by \reff{XD} and in particular due to the path dependence, it is clear the running time of the algorithm grows quadratically (rather than linearly) in $D$. 

(ii) By standard arguments, under mild regularity conditions one can easily show that 
\beaa
 \big|\dbE^Q[G(X_T)] - Y^{Q,D}_0\big| \le C \sup_{k, l} {|a_{k,l}|\over \sqrt{\d}} h = C C_{Q^*_N}  \sqrt{M\over T} {T\over D}= C C_{Q^*_N} {\sqrt{MT}\over D}. 
 \eeaa
 So there is a balance between the quadratic running cost and this error estimate. Theoretically, given an error level $\e$, we shall choose the parameters $M, D$, and $Q$ which satisfy $C C_{Q^*_N} {\sqrt{MT}\over D}\le \e$ and minimize the computational cost. Since the algorithm is much more sensitive to $M, Q$ than to $D$,   we content ourselves in this paper to choose a reasonably large $D$ and we  identify  $Y^{Q,D}_0$ with $Y^Q_0$ to emphasize the dependence on $Q$. Indeed, our numerical results show that the total error is not sensitive to $D$, see Example \ref{eg0-D} below.
\end{rem}

In the rest of this section, all the numerics are based on the use of Python  3.7.6 under Quad-Core Intel Core i5 CPU (3.4 GHz). For the running time, we use  $s$ and $ms$ to denote second and millisecond,  respectively.

\subsection{An illustrative one dimensional linear model}
\label{sect-eg0}
In this subsection we present a one dimensional numerical example:
\bea
\label{eg0-X}
X_t=x_0+\int_0^t (t-r)^{H-{1\over 2}}  dB_r,\q Y_0 = \dbE[G(X_T)].
\eea
In this case, $X_T \sim$ Normal$(x_0, {T^{2H}\over 2H})$, so essentially we can compute the true value of $Y_0$:
\bea
\label{eg0-Y0}
Y^{true}_0 =  {1\over \sqrt{2\pi}}\int_\dbR G(x_0 + {\sqrt{T^{2H}\over 2H}} x) e^{-{x^2\over 2}}dx.
\eea 
We use $Y^{cub}_0$, $Y^{mul, M}_0$, and $Y^{Euler}_0$  to denote the values computed by using the one period cubature formula,  the  multiple period cubature formula with $M$ periods, and the Euler scheme, respectively.
We shall compare our numerical results with this true value. In particular, since the cubature method is deterministic, while $Y^{Euler}_0$ is random, we shall explain how we compare the numerical results of the cubature methods with $Y^{Euler}_0$.

Our first example shows that, when $T$ is small, the one period method in \reff{cubature2} is more efficient than the multiple period method in \reff{cubaturem} with $M=1$, especially when $H$ is small. We remark that, although $M=1$,  \reff{cubature2} and  \reff{cubaturem} have different kernels and thus have different cubature paths. 

\begin{eg}
\label{eg0-1-M}
Consider \reff{eg0-X} with $G(x)=x^2$ and $x_0=0$. Then, with order $N=3$, 
\bea
\label{1-M}
|Y^{cub}_0 - Y^{true}_0| =  0 ~ <~  |Y^{mul,1}_0 - Y^{true}_0| = {H_-^2\over 2H H_+^2}  T^{2H}.
\eea
We see that the last error is small when $T$ is small or $H$ is large. However, it is still larger than the error of $Y^{cub}_0$ which is $0$ in this case.
\end{eg}
\proof  First, by \reff{eg0-Y0} it is clear that  $Y^{true}_0 = {T^{2H}\over 2H}$.

For the one period cubature method with $N=3$, by \reff{Xomega} and \reff{QK1} we have
\beaa
X_T(\o_1) &=& \int_0^{T\over 2} (T-t)^{H_-} {a_1\over \sqrt{T}} dt + \int_{T\over 2}^T (T-t)^{H_-} {a_2\over \sqrt{T}} dt \\
&=& {a_1\over \sqrt{T}} {1\over H_+}[T^{H_+} - ({T\over 2})^{H_+}] + {a_2\over \sqrt{T}} {1\over H_+} ({T\over 2})^{H_+} = \Big[a_1 ( 2^{H_+}-1) + a_2\Big] {T^H\over H_+ 2^{H_+}},\\
X_T(\o_2) &=& - X_T(\o_1).
\eeaa
Then, by \reff{Y0EG2} and \reff{Equation1},
\beaa
Y^{cub}_0 = {1\over 2}\sum_{k=1}^2|X_T(\o_k)|^2  =  \big[a_1 ( 2^{H_+}-1) + a_2\big]^2 {T^{2H}\over H_+^2 2^{2H_+}} = {H_+^2 2^{2H_+}\over 2H} {T^{2H}\over H_+^2 2^{2H_+}} ={T^{2H}\over 2H}= Y^{true}_0.
\eeaa

However, for the multiple period cubature method with $M=1$ and $N=3$, by \reff{Qm3-2} and \reff{Qm3-3} we have: by abusing the notations $\o_k$,
\beaa
&\dis X_T(\o_1) = \int_0^T (T-t)^{H_-} {1\over \sqrt{T}} dt = {T^H\over H_+},\q X_T(\o_2) = -X_T(\o_1),\\
&\dis Y^{mul,1}_0 = {1\over 2}\Big[|X_T(\o_1)|^2 + |X_T(\o_2) |^2\Big] =  {T^{2H}\over H_+^2},\\
&\dis |Y^{mul,1}_0 - Y^{true}_0| = |{1\over H_+^2} - {1\over 2H}| T^{2H} = {H_-^2\over 2H H_+^2}  T^{2H}.
\eeaa

\vspace{-8mm}
\qed

\ms

Our next example shows that, when $T$ is large, the one period algorithm fails, but the multiple period algorithm does converge when $M$ becomes large, as we expect.

\begin{eg}\label{eg0-M}
Consider \reff{eg0-X} with $T=3$, $H=5/2$, $G(x)=(x-1/2)^+$, $x_0=0.56$.  We compute the value $Y^{cub}_0$  from the one period model with $N=3$ and $Y_0^{mul, M}$ from the  multiple period model with $N=3$ and $M$ from $1$ to $5$. The cubature paths are constructed as in Example \ref{eg0-1-M}, with $D= 300 $ in Step 4.  Numerical results are reported in Table \ref{Table-eg0-M}. 
\begin{table}[t!]
\centering
    \begin{tabular}{|c | c | c | c | c| c| c|}
     \hline   
$Y^{true}_0$& $Y_0^{cub}$ & $Y^{mul, 1}_0$ & $Y^{mul, 2}_0$& $Y^{mul, 3}_0$& $Y^{mul, 4}_0$&  $Y^{mul, 5}_0$ \\ \hline
$2.8112$& $3.5157$   & $2.6281$ & $3.2450$& $3.1967$& $3.0340$&  $2.8883$ \\ \hline
      \end{tabular}
\caption{The numerical results for Example \ref{eg0-M}} 
     \label{Table-eg0-M} 
\end{table}
\end{eg}
\no We remark that, here the function $G$ is not as smooth as required in Theorems \ref{thm-Qest} and \ref{thm-QestMultiple}. Nevertheless we see from the numerical results that the cubature method still works well. 

We now present an example to compare the accuracy between $N=3$ and $N=5$.
\begin{eg}\label{eg0-N}
Consider \reff{eg0-X} with $T=0.3$,  $H=3/2$, and we try three different $G$ with corresponding initial value $x_0$.  We choose $D= 30$ in Step 4, and compute $Y_0^{mul, 2}$ with $N=3$ and $N=5$, respectively. The numerical results are reported in Table \ref{Tabel-eg0-N}. 
\end{eg}
\begin{table}[t!]
\centering
    \begin{tabular}{|c | c | c | c |}
     \hline   
$G(x)~\slash~ x_0$  & $\cos(x)~\slash~  1$  & $x^2~\slash~  1$  & $(x-1/2)^+~\slash~  0.56$\\ \hline
$Y^{true}_0$ & $0.5378641$   & $1.0090376$ & $0.0751964$  \\ \hline
$Y^{mul, 2}_0(N=3)$ &  $ 0.5380251$ &  $1.0084375$ & $0.0740474$    \\ \hline
$Y^{mul, 2}_0(N=5)$ &  $ 0.5380277$  &  $1.0084375$ & $0.0751558$   \\ \hline
      \end{tabular}
\caption{The numerical results for Example \ref{eg0-N}} 
     \label{Tabel-eg0-N} 
\end{table}

We note that, although we have a better rate of convergence in Theorem \ref{thm-QestMultiple} when $N=5$, the numerical results in this example do not show such improvement. One explanation is that the constant $C_N$ in \reff{QestMultiple2} becomes larger when we increase $N$, see Remark \ref{rem-cubaturem} (iv). The numerical results for the one period cubature method does not show significant improvement either when we increase the order from $N=3$ to $N=5$. However, for the  fractional stochastic volatility model  \reff{Heston}, as we saw in \S \ref{sect-N3-2d}, the one period cubature method with $N=3$ does not depend on $\rho$ at all, and thus it is clearly not as good as the  one period cubature method with $N=5$.

Our next example compares the one period cubature method with the Euler scheme. More examples concerning this comparison will be presented in the next two subsections.
\begin{eg}
\label{eg0-Euler}
Consider \reff{eg0-X} with $H={3\over 2}$, $T=0.2$, $G(x) = \cos(x)$, and $x_0=1$. We choose $D= 12$ in Step 4, and $Y^{cub}_0$ is computed with $N=5$. The numerical results are shown in Table \ref{Table-eg0-Euler}. 
\end{eg}

\begin{figure}
  \begin{minipage}[b]{.45\linewidth}
    \centering
    \begin{tabular}{| c | c | c | c |}
     \hline   
k  & 1& 2   \\ \hline

$\l_k$ &$ 0.15332891$ & $0.34667109$  \\ \hline

$a_{k,1}$   & $-3.04533315$  & $1.57981296$   \\ \hline
   
$a_{k,2}$   & $0.71729258$  & $-2.08974376$  \\ \hline

$a_{k,3}$   & $-0.60085202$  & $2.33258457 $  \\ \hline

$a_{k,4}$   & $0.12029985$  & $-4.5060389$   \\ \hline
      \end{tabular}
     \captionof{table}{The approximate cubature paths  in  Example \ref{eg0-Euler}} 
     \label{Table-paths} 
  \end{minipage}\hfill
  \begin{minipage}[b]{.45\linewidth}
\centering
\includegraphics[scale=0.45]{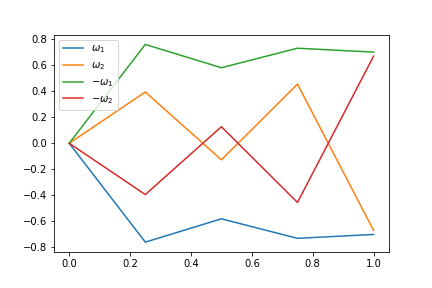}
 \captionof{figure}{The approximate cubature paths  in  Example \ref{eg0-Euler} (rescaling to $T=1$)} 
\label{Fig-paths}
  \end{minipage}
\end{figure}

We now explain the numerical results in Table \ref{Table-eg0-Euler}.  First, in this case $Y^{true}_0=0.53959$. Next, by solving \reff{optimization} numerically, we report the approximate solution of \reff{Equation2} in Table \ref{Table-paths}, see also Figure \ref{Fig-paths} for the plot of the  four approximate cubature paths (after rescaling to $T=1$). We then obtain $Y_0^{cub}= 0.54005$.  In Table \ref{Table-eg0-Euler}, the cubature error $e^{cub}:=|Y^{cub}_0-Y^{true}_0| = 0.00046$. The reported running time is for Steps 4 and 5 only, since Steps 1-3 can be completed offline once for all. 

\begin{table}[t!]
\centering
    \begin{tabular}{| c | c | c |c| c |}
     \hline   
$\widetilde M_{Euler}$& $100$  & $500$& $1000$   \\ \hline

$e^{cub}$ (cubature time)& $0.00046$ ($1.86$ms)& $0.00046$ ($1.86$ms) & $0.00046$ ($1.86$ms) \\ \hline

$e^{Euler}_{mean}$ (Euler time)& $0.00338$ ($14.2$ms)  & 0.00156  ($70$ms)  & $0.00114$ ($141$ms)   \\ \hline
   
$SD^{Euler}$ (percentile)& 0.0026 ($19.4\%$) & 0.00119 ($24\%$)  & 0.00081 ($27.5\%$) \\ \hline
      \end{tabular}
\caption{The numerical results for Example \ref{eg0-Euler}}
     \label{Table-eg0-Euler} 
\end{table}

For the Euler scheme, we also set time discretization step $D=12$. Let $\widetilde M_{Euler}$ denote the sample size in Euler scheme, namely the number of simulated paths of the Brownian motion. Clearly both the approximate value $Y^{Euler}_0$ and the running time depend on $\widetilde M_{Euler}$, in particular, the latter is proportional to $\widetilde M_{Euler}$.  Note that $Y^{Euler}_0$ is random.  We repeat the Euler scheme $1000$ times, each with sample size $\widetilde M_{Euler}$, and obtain $Y^{Euler,i}_0$ and the corresponding Euler scheme errors $e^{Euler}_i:=|Y^{Euler,i}_0-Y^{true}_0|$, $1\le i\le 1000$. We shall use the sample median $e^{Euler}_{median}$ of $\{e^{Euler}_i\}_{1\le i\le 1000}$ to measure the accuracy of the Euler scheme. We can then compare $e^{cub}$ and $e^{Euler}_{median}$, and to have a more precise comparison, we will actually compute the percentile of the cubature error $e^{cub}$ among the Euler scheme errors $\{e^{Euler}_i\}_{1\le i\le 1000}$: $\a$-percentile means about  $1000\times \a\%=10\a$ of   $\{e^{Euler}_i\}_{1\le i\le 1000}$ are smaller than $e^{cub}$. So $50\%$ roughly means $e^{cub}= e^{Euler}_{median}$  and the two methods have the same accuracy; while $\a\%<50\%$ means that $e^{cub} < e^{Euler}_{median}$ and the cubature method has better accuracy: the smaller $\a\%$ is, the better the cubature method outperforms.  Moreover, since $\{e^{Euler}_i\}_{1\le i\le 1000}$ are i.i.d., we may use the Normal approximation to compute the percentile. We will report their mean $e^{Euler}_{mean} := {1\over 1000}\sum_{i=1}^{1000} e^{Euler}_i \approx e^{Euler}_{median}$ and  standard deviation $SD^{Euler}$, and then the percentile $\a\% \approx \Phi({e^{cub}-e^{Euler}_{mean}\over SD^{Euler}})$, where $\Phi$ is the cdf of the standard normal.

For the above example, we test three cases: $\widetilde M_{Euler}= 100,  500, 1000$, and the numerical results are reported in Table \ref{Table-eg0-Euler}. As we see, when $\widetilde M_{Euler}= 500$, the Euler schemes takes  $70$ milliseconds (for each run, not for $1000$ runs), which is about $38$ times slower than the $1.86$ milliseconds used for the cubature method, and the percentile of the cubature method is $24\%$. So the cubature method outperforms the Euler scheme both in running time and in accuracy. When we increase the sample size $\widetilde M_{Euler}$ to $1000$, the percentile increases to $27\%$, so the cubature method still outperforms in accuracy. In this case the running time of the Euler scheme increases to $141$ milliseconds, which is about $76$ times slower than the cubature method.  On the other hand, if we  decrease the sample size $\widetilde M_{Euler}$ to $100$, the running time of the Euler scheme drops to $14.2$ milliseconds, which is still $7.6$ times slower than the cubature method, but the accuracy deteriorates further with a percentile $19.4\%$. So in all three cases, the cubature method outperforms the Euler scheme significantly both in running time and in accuracy. 

We conclude this subsection with an example concerning the impact of $D$.
\begin{eg}
\label{eg0-D}
Consider the same setting as in Example \ref{eg0-Euler}, but we try three different $D$'s. The numerical results are shown in Table \ref{Table-eg0-D}.\footnote{The running time for the Euler scheme grows quadratically in $D$, as expected. However, the running time for the cubature method grows slower than quadratically, especially when $D$ is not that large. This is possibly because in our code, for the sake of readability, there is a relatively time consuming step whose cost grows linearly in $D$. The efficiency of our cubature method could be improved slightly further, when $D$ is small, if we write the code in a more straightforward way.}
\begin{table}[t!]
\centering
    \begin{tabular}{| c | c | c |c| c |}
     \hline   
$D$& $12$  & $60$& $120$   \\ \hline

$e^{cub}$ (cubature time)& $0.00046$ ($1.86$ms)& $0.00046$ ($15.3$ms) & $0.00046$ ($54.7$ms) \\ \hline

$e^{Euler}_{mean}$ (time, $\tilde M=100$)& $0.00338$ ($14.2$ms)  & $0.0034$ ($314$ms)   & $0.00336$ ($1.24$s)   \\ \hline
   
$e^{Euler}_{mean}$ (time, $\tilde M=500$)&  $0.00156$  ($70$ms)  & $0.00161$ ($1.54$s)   & $0.00108$  ($5.9$s)  \\ \hline
$e^{Euler}_{mean}$ (time, $\tilde M=1000$)& $0.00114$ ($141$ms) & $0.00113$ ($3.11$s)  & $0.00107$ ($11.5$s)\\ \hline
      \end{tabular}
\caption{The numerical results for Example \ref{eg0-D}}
     \label{Table-eg0-D} 
\end{table}
\end{eg}
As we can see, the cubature method is not sensitive to $D$. The Euler scheme does rely on our choices of $D$ and $\tilde M_{Euler}$. However, in all the above choices, the cubature method  outperforms both in running time and in accuracy.

\subsection{A one dimensional nonlinear model}
We now consider the following nonlinear model, but still in $1$-dimensional setting:
\bea
\label{eg1-X}
X_t=x_0+\int_0^t (t-r)^{H-{1\over 2}}  \cos(X_s) dB_r,\q Y_0 = \dbE[G(X_T)].
\eea
We shall use one period cubature method with $N=5$ when $T$ is small, and multiple period cubature method with $N=3$ and appropriate $M$ when $T$ is large. Our main purpose is to compare the efficiency of the cubature method with that of the Euler scheme. 

We first note that, by Remark \ref{rem-optimization}  (ii), the cubature paths for \reff{eg1-X} is the same as those for \reff{eg0-X}. In particular, for the one period method with $N=5$, we may continue to use the paths in Table \ref{Table-paths}.  For comparison purpose, we will use the same $D$ for the cubature method and the Euler scheme. For the $\widetilde M_{Euler}$ in Euler scheme, there is an obvious tradeoff between the running time and the accuracy. While one may try to find an ``optimal" $\tilde M_{Euler}$ for a given $D$, such an analysis relies on a precise idea on the constants involved in the error estimates, as in Remark \ref{rem-D} (ii). Since our main focus is the cubature method and since our examples show that the cubature method outperforms significantly (under our strong conditions), we do not go through that analysis. Instead,  unless stated otherwise, for simplicity in the rest of this section we shall always set 
\beaa
\widetilde M_{Euler} = 500, \q\mbox{and we repeat the Euler scheme for $1000$ times},
\eeaa
 each time with  $\widetilde M_{Euler}$ simulation paths. We use $\{Y^{Euler,i}_0\}_{1\le i\le 1000}$ and $Y^{true}_0$ to compute $e^{Euler}_{mean}$ and $SD^{Euler}$. However, in this case we are not able to compute the exact value of $Y^{true}_0$ as in \reff{eg0-Y0}.  Since the convergence of the Euler scheme approximations is  well understood, for comparison purpose we shall set the true value as the sample mean of the Euler scheme approximations:
\bea
\label{eg1-Y0}
Y^{true}_0 = {1\over 1000} \sum_{i=1}^{1000} Y^{Euler,i}_0.
\eea

In the first example we show the impact of the regularity of  $G$.

\begin{eg}
\label{eg1-Greg}
Consider \reff{eg1-X} with $H={3\over 2}$,  $T=0.2$, and we consider three different $G$'s with corresponding $x_0$. We choose $D=12$, and compare the one period cubature method with $N=5$ with the Euler scheme. The numerical results are reported in Table \ref{Table-eg1-Greg}.
\begin{table}[t!]
\centering
    \begin{tabular}{| c | c | c | c |}
     \hline   
$G(x)$~\slash~ $x_0$  & $\cos(x)$~\slash~ 1  & $x^2$~\slash~ 1& $(x-0.5)^+$ ~\slash~ 0.56  \\ \hline
$Y^{true}_0$ & 0.5401  & 1.00073 & 0.0617  \\ \hline
$Y^{cub}_0$ &  0.5402 &1.00041  &    0.0601  \\ \hline
   $e^{cub}$ (time) & 0.0001(1.75ms) &0.00032(1.66ms)  &0.0016 (1.63ms)    \\ \hline
  $e^{Euler}_{mean}$ (time)  & 0.0008(80.4ms) & 0.00200(82ms )  & 0.00147(75ms)  \\ \hline
 $SD^{Euler}$(percentile) & 0.00064(18\%)  & 0.00149(20.1\%)&  0.0011(54.3\%) \\ \hline
      \end{tabular}
\caption{The numerical results for Example \ref{eg1-Greg}} 
     \label{Table-eg1-Greg}
\end{table}
\end{eg}
As we can see, the cubature method outperforms the Euler scheme in all these three cases. However, when $G$ becomes less smooth,  the advantage of the cubature method fades away, which is consistent with the theoretical observation in Remark \ref{rem-Taylor}.

The next example illustrates the impact of $T$.
\begin{eg} 
\label{eg1-T}
Consider \reff{eg1-X} with $H={3\over 2}$, $G(x) = x^2$, $x_0=1$,  and three different values of $T$: $0.2, ~0.5, ~0.8$. We choose $D=12, ~30,~ 48$, respectively, and compare the one period cubature method with $N=5$ with the Euler scheme.   The numerical results are reported in Table \ref{Table-eg1-T}.
\end{eg}
\begin{table}[t!]
\centering
    \begin{tabular}{|c| c | c | c |}
     \hline   
$T$  & 0.2  & 0.5& 0.8   \\ \hline
$Y^{true}_0$ & 1.00073  & 1.0109 & 1.045  \\ \hline
$Y^{cub}_0$ &  1.00041 &1.0056  &   1.022   \\ \hline
   $e^{cub}$ (time) & 0.00032(1.66ms)  &0.0054(5.25ms)  &0.0234(11ms)    \\ \hline
  $e^{Euler}_{mean}$ (time)  & 0.00200 (82ms)  & 0.0332(445ms) &0.0152(1.16s)  \\ \hline
 $SD^{Euler}$(percentile) & 0.00149(20.1\%)  & 0.0098(57.4\%)&  0.0115(69.4\%) \\ \hline
      \end{tabular}
\caption{The numerical results for Example \ref{eg1-T}} 
     \label{Table-eg1-T} 
\end{table}
Again consistent with our theoretical result, the performance of the cubature method decays when $T$ gets large. In the above example, the cubature method obviously outperforms the Euler scheme when $T=0.2$, and still works better when $T=0.5$, but in the case $T=0.8$, there is a tradeoff between the speed and the accuracy and it is hard to claim the cubature method is more efficient. In the last case, we shall use the  multiple period cubature method, as we do in the next example, and we can easily see that the cubature method outperforms the Euler scheme again.

\begin{eg} 
\label{eg1-Tbig}
Consider \reff{eg1-X} with $H={3\over 2}$,  $T=1$, and three different $G$'s with corresponding $x_0$.  We choose $D=100$, and compare the  multiple period cubature method with $M=5$ and $N=3$ with the Euler scheme.   The numerical results are reported in Table \ref{Table-eg1-Tbig}. 
\begin{table}[t!]
\centering
    \begin{tabular}{| c | c | c | c |}
     \hline   
$G(x)$~\slash~ $x_0$  & $\cos(x)$~\slash~ 1  & $x^2$~\slash~ 1& $(x-0.5)^+$ ~\slash~ 0.56  \\ \hline
$Y^{true}_0$ & $0.5136$  & $1.098$ & $0.2275$  \\ \hline
$Y^{mul, 5}_0$ &  $ 0.5186$ & $1.084$ &  $0.2297$    \\ \hline
$e^{mul, 5}$ (time) &  $ 0.005(300ms)$ & $0.014(302ms)$ &  $0.0022(302ms)$    \\ \hline
  $e^{Euler}_{mean}$ (time)  &0.0074(4.2s)  & 0.023 ($4.36$s)  & 0.011 ($4.41$s)  \\ \hline
 $SD^{Euler}$(percentile) &  $0.0056$ (37.37\%) & $0.017$ ($35.1\%$) & $0.0089$ ($21\%$)  \\ \hline
      \end{tabular}
\caption{The numerical results for Example \ref{eg1-Tbig}} 
     \label{Table-eg1-Tbig} 
\end{table}
\end{eg}

\subsection{A fractional stochastic volatility model}
In this section we consider the following special case of \reff{Heston} with $H={3\over 2}$ and $\rho={1\over 2}$:
\bea
\label{RVM}
\left.\ba{c}
\dis dS_t= S_t b_1(U_t)dt+ S_t \sigma_1(U_t)\circ dB_t^1, \ms\\
\dis	U_t= 1 +\int_0^t [t-s][{1\over 2}-{1\over 3}U_s] ds+\int_0^t [t-s] \si_2(U_s)\circ dB_t^2,\\
\ea\right.
\eea
Again we will use one period cubature method with $N=5$ when $T$ is small,  multiple period cubature method with $N=3$ and appropriate $M$ when $T$ is large, and we shall compare the efficiency between the  cubature method and the Euler scheme. 

\begin{eg}
\label{eg-rvm1}
Consider \reff{RVM} with $T=0.1$,  $b_1(U)=U$, $\si_1(U)=\si_2(U)=\cos(U)$, and we consider three different $G$'s with corresponding $S_0$.  We choose $D=12$, and compare the one period cubature method with  $N=5$ with the Euler scheme. The numerical results are reported in Table \ref{Table-rvm1}.
\end{eg}
\begin{table}[t!]
\centering
    \begin{tabular}{|c | c | c | c |}
     \hline   
$G(x)$~\slash~ $S_0$  & $\cos(x)$~\slash~ 1  & $x^2$~\slash~ 1& $(x-0.5)^+$ ~\slash~ 0.56  \\ \hline
$Y^{true}_0$ & $ 0.4270$  & $1.2947$ & $ 0.1320$  \\ \hline
$Y^{cub}_0$ &  $ 0.4257$ & $1.2967$ &  $0.1283$    \\ \hline
   $e^{cub}$ (time) &  $ 0.0013(8.53ms)$ & $0.0020(8.9ms)$ &  $0.0037(8.68ms)$    \\ \hline
  $e^{Euler}_{mean}$ (time)  & 0.0063 (273ms)  & 0.0157(283ms) &0.0037(272ms)    \\ \hline
 $SD^{Euler}$(percentile) &  $0.0047$ (21.3\%) & $0.0119$ ($19.2\%$) & $0.0028$ ($50\%$)  \\ \hline
      \end{tabular}
      \centering
\caption{The numerical results for Example \ref{eg-rvm1}} 
     \label{Table-rvm1} 
\end{table}

\begin{figure}[t!]
\centering
\begin{subfigure}{.45\textwidth}
\centering	
\includegraphics[scale=0.5]{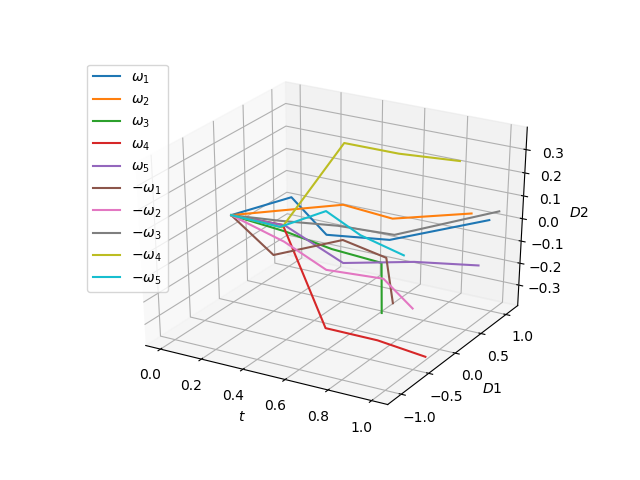}
\end{subfigure}
\caption{The cubature paths for the model \reff{RVM} (rescaling to $T=1$)} 
\label{Fig-rvm}
\end{figure}

For the cubature method, we first compute the cubature paths following the same idea as in Remark \ref{rem-optimization}. By \S\ref{sect-N5-2d2} below we choose $W=5$ and $L=4$. Then  we obtain
\beaa
&\l_1=0.0247245002,\q \l_2=0.0561159547, \q \l_3= 0.417734596e, \\
&\l_4=0.00142494883,\q \l_5= 4.44061201e-17
 \eeaa
and the $10$ paths are plotted in Figure \ref{Fig-rvm} (after rescaling to $T=1$). 

\begin{eg}
\label{eg-rvm2}
Consider the same setting as in Example \ref{eg-rvm1}, except that $b_1(U)=\si_1(U)=\si_2(U)=\sqrt{U}$. The numerical results are reported in Table \ref{Table-rvm2}.
\end{eg}

 \begin{table}
 \begin{tabular}{|c | c | c | c |}
     \hline   
$G(x)$~\slash~ $S_0$  & $\cos(x)$~\slash~ 1  & $x^2$~\slash~ 1& $(x-0.5)^+$ ~\slash~ 0.56  \\ \hline
$Y^{true}_0$ & 0.37897   & 1.4932 &  0.17098  \\ \hline
$Y^{cub}_0$ & 0.3698  & 1.4887  & 0.17797 \\ \hline
   $e^{cub}$ (time) &  $ 0.00917(8.06ms)$ & $0.0045(8.39ms)$ &  $0.00699(8.07ms)$    \\ \hline
  $e^{Euler}_{mean}$ (time)  & 0.0119(306ms)  & 0.0361(311ms) &0.007(304ms)   \\ \hline
 $SD^{Euler}$(percentile) &  0.0092(40.8\%)  & 0.0270(19.1\%) &  0.0052(50\%)  \\ \hline
      \end{tabular}
      \centering
\caption{The numerical results for Example \ref{eg-rvm2}} 
     \label{Table-rvm2} 
\end{table}
We remark that Example \ref{eg-rvm2} uses the same cubature paths as in Example \ref{eg-rvm1}.

\begin{eg}
\label{eg-rvm1-H}
Consider \reff{RVM} with $T=1$, $b_1(U)=U$, $\si_1(U)=\si_2(U)=\cos(U)$, $G(x)=(x-1/2)^+$, $S_0=0.56$, and consider three different $H$. We shall compare the efficiency of the multiple period cubature method with $M=3$, $N=3$ and the Euler scheme. We set $D= 100$. Recalling $\rho = {1\over 2}$, for the cubature method we use $\theta_1=\frac{\pi}{6}$ and $\theta_2=-\frac{\pi}{6}$  in \reff{2dQm5}. We repeat the Euler scheme $100$ times (instead of $1000$ times).
The numerical results are reported in Table \ref{Table-rvm1-H}. 
\end{eg}
\begin{table}[t!]
\centering
    \begin{tabular}{|c | c | c | c |}
     \hline   
$H$  & $1$  & $3/2$& $5/2$   \\ \hline
$Y^{true}_0$ & $ 1.36$  & $1.33$ & $ 1.2957$  \\ \hline
$Y^{mul, 3}_0$ &  $ 1.299$ & $1.302$ &  $1.286$    \\ \hline
$e^{mul, 3}$(time) &  $ 0.061(2.23s)$ & $0.028(2.26s)$ &  $0.0097(2.42s)$    \\ \hline
   $e^{Euler}_{mean}$(time) & 0.033 (15.8s)  & 0.036(15.8s) &0.037(15.5s)    \\ \hline
 $SD^{Euler}$(percentile) &  $0.024$ (81.1\%) & $0.028$ ($41.7\%$) & $0.0260$ ($21\%$)  \\ \hline
      \end{tabular}
      \centering
\caption{The numerical results for Example \ref{eg-rvm1-H}} 
     \label{Table-rvm1-H} 
\end{table}
We see that the cubature method outperforms the Euler scheme when $H={3\over 2}$ and $H={5\over 2}$, especially in the latter case, but does not seem to work well when $H=1$.  This is consistent with our theoretical result. 

\subsection{Some concluding remarks}
\label{sect-conclude}
We first note that, our theoretical convergence analysis for the cubature method, namely Theorems \ref{thm-Qest} and \ref{thm-QestMultiple}, is complete, provided sufficient regularities on $K$ and $(V, G)$ (corresponding to $N$). In particular, it holds true for arbitrary dimensions and arbitrarily large $T$ (for Theorem \ref{thm-QestMultiple}).

For the numerical efficiency, in its realm, the cubature method has clear advantages over the Euler scheme. In light of Remark \ref{rem-cubaturem}, the cubature method requires the following three conditions though: (i) sufficient regularity, so as to obtain the desired error estimate; (ii) low dimension (and relatively small $N$), so that $W$ can be relatively small; and (iii) not too large $T$, so that $M$ can be relatively small. We remark that the standard cubature method in \cite{LV03} for diffusions also requires these conditions. However, the constraints are more severe in the Volterra setting here, for example, the constant $C_N$ in  \reff{QestMultiple2} is larger here, and for given dimensions and $N$, there are more equations required in \reff{cubaturem} and hence we may need a larger $W$. When the number of cubature paths $(2W)^M$ is large (recalling again Remark \ref{rem-cubaturem} (ii)), it will be interesting to explore if the approach in \cite{CM12, CM14} could help reduce the complexity, which we shall leave for future research. The less smooth case, especially when $H<{1\over 2}$, requires novel idea to extend our approach. 

We shall also remark that the parameters $M, N$ in the cubature method cannot be too big.  For the Euler scheme, by increasing the sample size $\widetilde M_{Euler}$  gradually one may improve the accuracy ``continuously" at the price of sacrificing the speed. For the cubature method we have only limited choices on $M, N$ and thus lose the flexibility of improving its accuracy ``continuously". Consequently, the cubature method is more appropriate in the situations where one has strong requirement on the speed but is less stringent on the accuracy.

 \section{Appendix}
\label{sect-appendix}
\setcounter{equation}{0}

\no{\bf Proof of Proposition \ref{prop-ureg}. Recall \eqref{u}. Following rather standard arguments, we see that $\pa_\bx u(t,\cd)$ exists and has the following representation:
\bea
\label{tdX}
\left.\ba{c}
\dis \la \pa_\bx u(t, \th), \eta\ra = \dbE\big[\pa_x G(X^{t,\th}_T) \cd \td_\eta X^{t,\th}_T\big],\q \th, \eta\in \dbX_t, \ms \\
\dis\mbox{where}\q \td_\eta X^{t,\th,i}_s = \eta^i_s + \sum_{j=0}^d \int_t^s K_i(s, r) \pa_x V^i_j (X^{t, \th}_r) \td_\eta X^{t,\th,j}_r \circ dB^j_r,\q i=1,\cds, d_1.
\ea\right.
\eea
We note that  at above we need the second derivative of $V^i_j$ exists so that the Stratonovich integration $\circ$ makes sense. Then we see immediately $\|\pa_\bx u(t, \th)\|\le C e^{C^m_N(T-t)}$. By similar arguments we can prove the results for higher order derivatives. In particular, we note that the $(N-1)$-th derivative of $u$ would involve the $N$-th derivative of $V^i_j$. 
\qed
}

\bs
\no{\bf Proof of Proposition \ref{prop-expansion}.}  We first note that, by Proposition \ref{prop-ureg} $u(t^m_0,\cd)\in C^{N+2}(\dbX_{t^m_0})$, so the right side of \reff{expansionN} makes sense.

When $N=1, 2$, one may verify easily that \reff{expansionN} reduces to \reff{expansion1}, \reff{expansion2}, respectively. Assume \reff{expansionN} holds  for $N-1$. For $\vec i\in\cI_{N}$, $\vec j\in\cJ_{N}$, $\vec\k\in \cS_N$, $\vec t \in \dbT_{N}$, we have
\beaa
&&\dis \D_{N+1}(\vec i, \vec j, \vec \k; \vec t):= \cV(\vec  i,\vec j, \vec\k;\vec t, \Th^{\vec t}_{ t_{N}},\Th^{[t^m_0,T]}_{t_{N}}) - \cV(\vec  i,\vec j, \vec\k;\vec t, \Th^{\vec t}_{T_m},\Th^{[t^m_0,T]}_{T_m})\\
&&\dis =\prod_{\alpha = 1}^{N} \pa^{\vec \k, \alpha}_{\vec i}  V^{i_\a}_{j_\a} (\Th^{t_\a}_{ t_{N}}) \cK_+(\vec i,  \vec \k; \vec t) ~\big \la \pa_{\vec i}^{\vec \k}u(t^m_0, \Th^{[t^m_0,T]}_{t_{N}}) , \vec \cK_0(\vec i,  \vec\k; \vec t)  \big\ra \\
&&\dis \q- \prod_{\alpha = 1}^{N} \pa^{\vec \k, \alpha}_{\vec i}  V^{i_\a}_{j_\a} (\Th^{t_\a}_{T_m})\cK_+(\vec i,  \vec \k; \vec t) ~\big \la \pa_{\vec i}^{\vec \k}u(t^m_0, \Th^{[t^m_0,T]}_{T_m}) , \vec \cK_0(\vec i,  \vec\k; \vec t)  \big\ra.
\eeaa
Note that, for $s\in [T_m, t_N]$, by It\^o's formula and Proposition \ref{prop-Ito} we have
\bea
\label{dVu}
\left.\ba{lll}
\dis d \pa^{\vec \k, \alpha}_{\vec i}  V^{i_\a}_{j_\a} (\Th^{t_\a}_s) = \sum_{\tilde i =1}^{d_1} \sum_{\tilde j=0}^d \pa_{x_{\tilde i}} [\pa^{\vec \k, \alpha}_{\vec i}  V^{i_{\a}}_{j_{\a}}] (\Th^{t_{\a}}_s) K_{\tilde i}(t_{\a}, \tilde t) V^{\tilde i}_{\tilde j}(X_s) \circ dB^{\tilde j}_s,\ms\\
\dis d\big \la \pa_{\vec i}^{\vec \k}u(t^m_0, \Th^{[t^m_0,T]}_s) , \vec \cK_0(\vec i,  \vec\k; \vec t)  \big\ra \ms\\
\dis =  \sum_{\tilde i =1}^{d_1} \sum_{\tilde j=0}^d \big \la \pa_{\bx_{\tilde i}}\pa_{\vec i}^{\vec \k}u(t^m_0, \Th^{[t^m_0,T]}_s) , (\vec \cK_0(\vec i,  \vec\k; \vec t), K_{\tilde i, s}^{[t^m_0, T]})  \big\ra V^{\tilde i}_{\tilde j}(X_s) \circ dB^{\tilde j}_s.
\ea\right.
\eea
Then, by It\^o's formula we have, for given $\vec i\in\cI_{N}$, $\vec j\in\cJ_{N}$, $\vec t \in \dbT_{N}$,
\bea
\label{XiN}
&&\dis \D_{N+1}(\vec i, \vec j, \vec\k; \vec t)= \sum_{\tilde i =1}^{d_1} \sum_{\tilde j=0}^d  \int_{T_m}^{t_N}  \Xi_{N+1}\big((\vec i, \tilde i), (\vec j, \tilde j), \vec\k; (\vec t, s)\big)  \circ dB^{\tilde j}_s,\q\mbox{where}\nonumber\\
 &&\dis \Xi_{N+1}\big((\vec i, \tilde i), (\vec j, \tilde j), \vec\k; (\vec t, s)\big) := \sum_{\tilde\a=1}^{N} \prod_{\a \in\{1,\cds, N\}\backslash\{\tilde \a\}} \pa^{\vec \k, \alpha}_{\vec i}  V^{i_\a}_{j_\a} (\Th^{t_\a}_{s})\times\nonumber\\
&&\dis \qq \pa_{x_{\tilde i}} [\pa^{\vec \k, \tilde \alpha}_{\vec i}V^{i_{\tilde \a}}_{j_{\tilde \a}}](\Th^{t_{\tilde \a}}_{s}) V^{\tilde i}_{\tilde j}(\Th^s_s) \cK(\vec i,  \vec \k; \vec t) K_{\tilde i}(t_{\tilde \a}, s)  ~\big \la \pa_{\vec i}^{\vec \k}u(t^m_0, \Th^{[t^m_0,T]}_s) , \vec \cK_0(\vec i,  \vec\k; \vec t)  \big\ra \\
&&\dis \q+  \prod_{\a=1}^{N} \pa^{\vec \k, \alpha}_{\vec i}  V^{i_\a}_{j_\a} (\Th^{t_\a}_{s})V^{\tilde i}_{\tilde j}(\Th^s_s)\cK(\vec i,  \vec \k; \vec t)~\big \la \pa_{\bx_{\tilde i}}\pa_{\vec i}^{\vec \k}u(t^m_0, \Th^{[t^m_0,T]}_s) , (\vec \cK_0(\vec i,  \vec\k; \vec t), K_{\tilde i, s}^{[t^m_0, T]})  \big\ra.\nonumber
\eea
Thus, since \reff{expansionN} holds  for $N-1$ by induction assumption, we have
\beaa
&&\dis u(t^m_0, \Th^{[t^m_0,T]}_{t^m_0}) - u(t^m_0, \Th^{[t^m_0,T]}_{T_m})\\
&&\dis = \sum_{n=1}^{N-1}\sum_{\vec i \in\cI_n, \vec j\in\cJ_n, \vec\k\in \cS_n}\int_{\dbT^m_n} \cV(\vec i,\vec j, \vec\k;\vec t, \Th^{\vec t}_{T_m},\Th^{[t^m_0,T]}_{T_m}) \circ dB_{\vec t}^{\vec j} \\
&&\q +\sum_{\vec i\in\cI_{N}, \vec j\in\cJ_{N}, \vec\k\in \cS_N}\int_{\dbT^m_N} \cV(\vec  i,\vec j, \vec\k;\vec t, \Th^{\vec t}_{ t_{N}},\Th^{[t^m_0,T]}_{t_{N}}) \circ dB_{\vec t}^{\vec j}\\
&&=\sum_{n=1}^{N}\sum_{\vec i \in\cI_n, \vec j\in\cJ_n, \vec\k\in \cS_n}\int_{\dbT^m_n} \cV(\vec i,\vec j, \vec\k;\vec t, \Th^{\vec t}_{T_m},\Th^{[t^m_0,T]}_{T_m}) \circ dB_{\vec t}^{\vec j} \\
&&\dis\q +\sum_{(\vec i, \tilde i)\in\cI_{N+1}, (\vec j, \tilde j)\in\cJ_{N+1}, \vec\k\in \cS_N}\int_{\dbT^m_{N+1}} \Xi_{N+1}\big((\vec i, \tilde i), (\vec j, \tilde j), \vec \k; (\vec t, s)\big) \circ dB_{(\vec t, s)}^{(\vec j, \tilde j)}.
\eeaa 
Then it suffices to show that, 
\bea
\label{Xi=cV}
\Xi_{N+1}\big((\vec i, \tilde i), (\vec j, \tilde j), \vec\k; (\vec t, s)\big)  =  \sum_{\tilde\a=0}^N \cV((\vec i, \tilde i), (\vec j, \tilde j), (\vec\k, \tilde \a); (\vec t, s), \Th^{(\vec t, s)}_{s},\Th^{[t^m_0,T]}_s). 
\eea

Indeed, recall \reff{cKcV}, one can verify that
\beaa
 &&\dis \Xi_{N+1}\big((\vec i, \tilde i), (\vec j, \tilde j), \vec\k; (\vec t, s)\big) \\
&&\dis =  \sum_{\tilde\a=1}^{N} \prod_{\a=1}^N  \pa^{(\vec \k, \tilde \a), \alpha}_{(\vec i, \tilde i)}  V^{i_\a}_{j_\a} (\Th^{t_\a}_{s})V^{\tilde i}_{\tilde j}(\Th^s_s)\cK((\vec i, \tilde i),  (\vec \k, \tilde \a); (\vec t, s)) \times\\
&&\dis \qq\qq\qq   \big \la \pa_{(\vec i, \tilde i)}^{(\vec \k, \tilde \a)}u(t^m_0, \Th^{[t^m_0,T]}_s) , \vec \cK_0((\vec i, \tilde i),  (\vec\k, \tilde \a); (\vec t, s))  \big\ra \\
&&\dis \q+  \prod_{\a=1}^{N} \pa^{(\vec \k, 0), \alpha}_{(\vec i, \tilde i)}  V^{i_\a}_{j_\a}  (\Th^{t_\a}_{s}) V^{\tilde i}_{\tilde j}(\Th^s_s)\cK((\vec i, \tilde i),  (\vec \k, 0); (\vec t, s))\times\\
&&\dis \qq\qq\qq \big \la \pa_{(\vec i, \tilde i)}^{(\vec \k, 0)}u(t^m_0, \Th^{[t^m_0,T]}_s) , (\vec \cK_0((\vec i, \tilde i),  (\vec\k, 0); (\vec t,s))  \big\ra\\
&&\dis =  \sum_{\tilde\a=0}^{N} \prod_{\a=1}^{N+1}  \pa^{(\vec \k, \tilde \a), \alpha}_{(\vec i, \tilde i)}  V^{i_\a}_{j_\a} (\Th^{t_\a}_{s})\cK((\vec i, \tilde i),  (\vec \k, \tilde \a); (\vec t, s)) \times\\
&&\dis \qq\qq\qq   \big \la \pa_{(\vec i, \tilde i)}^{(\vec \k, \tilde \a)}u(t^m_0, \Th^{[t^m_0,T]}_s) , \vec \cK_0((\vec i, \tilde i),  (\vec\k, \tilde \a); (\vec t, s))  \big\ra \\
&&\dis =  \sum_{\tilde\a=0}^N \cV((\vec i, \tilde i), (\vec j, \tilde j), (\vec\k, \tilde \a); (\vec t, s), \Th^{(\vec t, s)}_{s},\Th^{[t^m_0,T]}_s).
\eeaa
This proves \reff{Xi=cV}, hence \reff{expansionN}  for $N$. 
\qed

\bs
\no{\bf Proof of Lemma \ref{lem-EdB}.}  Recall that $B^0_t=t$, the case $j_1=0$ is  obvious. We now assume $j_1>0$.  Fix ${T_m}\le s\le t^m_0$ and denote 
\beaa
\psi(t_1, t_2):= \int_{\dbT^m_{n-2}(t_2)} \f(t_1, t_2, \vec t_{-2}) \circ dB^{\vec j_{-2}}_{\vec t_{-2}}.
\eeaa
 Then, when $j_2\neq j_1$, we have
\beaa
\int_{\dbT^m_n(s)}\!\!\! \f(\vec t~\!) \circ dB^{\vec j}_{\vec t} = \int_{T_m}^{s}\Big[\int_{T_m}^{t_1}\psi(t_1, t_2)\circ dB^{j_2}_{t_2} \Big]\circ dB^{j_1}_{t_1} =\int_{T_m}^{s} \Big[\int_{T_m}^{t_1}\psi(t_1, t_2)\circ dB^{j_2}_{t_2} \Big] dB^{j_1}_{t_1},
\eeaa
coinciding with the Ito integral, and when $j_1=j_2$, 
\beaa
\int_{\dbT^m_n(s)}\!\!\! \f(\vec t~\!) \circ dB^{\vec j}_{\vec t} = \int_{T_m}^{s} \Big[\int_{T_m}^{t_1}\psi(t_1, t_2)\circ dB^{j_2}_{t_2} \Big] dB^{j_1}_{t_1} +{1\over 2}\int_{T_m}^{s} \psi(t_1, t_1) dt_1.
\eeaa
Then one may verify \reff{EdB} and \reff{E|dB|} straightforwardly.
\qed

\bs
\no{\bf Proof of Theorem \ref{thm-Taylor}.} First, similar to \reff{expansion2} and \reff{expansionN}, one can verify that 
\bea
\label{RNrep}
 \left.\ba{c}
 \dis R^m_N =  \sum_{n=1}^{N+1}  \sum_{\vec i\in \cI_n, \vec j\in \cJ_n}\big[ \1_{\{\|\vec j\|=N+1\}}  +  \1_{\{\|\vec j_{-1}\|=N, j_1=0\}}\big]\times\\
 \dis \int_{\dbT^m_n}  \sum_{\vec\k \in \cS_{n}} \cV(\vec  i,\vec j, \vec\k;\vec t, \Th^{\vec t}_{ t_{n}},\Th^{[t^m_0,T]}_{t_{n}}) \circ dB^{\vec j}_{\vec t}.
\ea\right.
\eea
Fix $n$, $\vec i\in \cI_n$, $\vec j\in \cJ_n$, and denote $\f(\vec t) :=  \sum_{\vec\k \in \cS_{n}} \cV(\vec  i,\vec j, \vec\k;\vec t, \Th^{\vec t}_{ t_{n}},\Th^{[t^m_0,T]}_{t_{n}})$ for $\vec t\in \dbT^m_n$. By \reff{E|dB|} and \reff{cJkm} we may prove by induction that, for $s\in [T_m, t^m_0]$, $l=1,\cds, n-1$, and for some constant $C_n$ which may depend on $n$, 
\beaa
\|\f(\cd)\|_{s, \vec j}^2 &\le& C_n\d^{\|\vec j_l\|}\esup_{T_m\le s_l\le \cds\le s_1\le s} \|\f(s_1,\cds, s_l,\cd)\|_{\vec j_{-l}, s_l}^2 \\
&&+ C_n\d^{\|\vec j_{l+1}\|}\esup_{T_m\le s_{l+1}\le \cds\le s_1\le s} \|\f(s_1,\cds, s_{l+1},\cd)\|_{\vec j_{-l-1}, s_{l+1}}^2
\eeaa
In particular, by setting $l=n-1$ we have
\bea
\label{ERNest1}
&\dis \dbE_m\Big[\Big|\int_{\dbT^m_n}\f(\vec t) \circ d B^{\vec j}_{\vec t}\Big|^2 \Big] \le C_n[\L_{n-1} +  \L_n],\q \mbox{where}\\
&\dis  \L_{n-1}:= \d^{\|\vec j_{n-1}\|} \esup_{\vec s\in \dbT^m_{n-1}} \dbE_m\Big[\Big|\int_{T_m}^{s_{n-1}} \f(\vec s, t_n) \circ dB^{j_n}_{t_n}\Big|^2\Big],\q \L_n := \d^{\|\vec j\|}\esup_{\vec s\in \dbT^m_n} \dbE_m\Big[\big|\f(\vec s)\big|^2\Big].\nonumber
\eea
By \reff{cJkm} we have $\L_n \le |A^m_{\|\vec j\|}|^2\d^{\|\vec j\|}$. Moreover, fix $\vec s\in \dbT_{n-1}$. When $j_n=0$, we have
\beaa
&&\dis \dbE_m\Big[\Big|\int_{T_m}^{s_{n-1}} \f(\vec s, t_n) \circ dB^{j_n}_{t_n}\Big|^2\Big] = \dbE_m\Big[\Big|\int_{T_m}^{s_{n-1}} \f(\vec s, t_n) d{t_n}\Big|^2\Big]\\
&&\dis \le \d\dbE_m\Big[\int_{T_m}^{s_{n-1}} |\f(\vec s, t_n)|^2 d{t_n}\Big] \le \d^2 \esup_{T_m\le t_n \le s_{n-1}} \dbE_m\Big[|\f(\vec s, t_n)|^2\Big]\le \d^2|A^m_{\|\vec j\|}|^2.
\eeaa
Then 
\beaa
\L_{n-1} \le \d^{\|\vec j_{n-1}\|} \d^2|A^m_{\|\vec j\|}|^2 = |A^m_{\|\vec j\|}|^2\d^{\|\vec j\|}.
\eeaa
When $j_n >0$, recalling \reff{cKcV} and by \reff{dVu} we have: denoting $s_n:= t_n$,
  \beaa
&&\dis \int_{T_m}^{s_{n-1}} \f(\vec s, t_n) \circ dB^{j_n}_{t_n}  \\
&&\dis = \int_{T_m}^{s_{n-1}} \sum_{\vec\k \in \cS_{n}}\prod_{\alpha = 1}^{n} \pa^{\vec \k, \alpha}_{\vec i}  V^{i_\a}_{j_\a} (\Th^{s_\a}_{t_n}) \cK_+(\vec i,  \vec \k; (\vec s, t_n)) ~\big \la \pa_{\vec i}^{\vec \k}u(t^m_0, \Th^{[t^m_0, T]}_{t_n}) , \vec \cK_0(\vec i,  \vec\k; (\vec s, t_n))  \big\ra\circ dB^{j_n}_{t_n}\\
&&\dis = \int_{T_m}^{s_{n-1}} \sum_{\vec\k \in \cS_{n}}\prod_{\alpha = 1}^{n} \pa^{\vec \k, \alpha}_{\vec i}  V^{i_\a}_{j_\a} (\Th^{s_\a}_{t_n}) \cK_+(\vec i,  \vec \k; (\vec s, t_n)) ~\big \la \pa_{\vec i}^{\vec \k}u(t^m_0, \Th^{[t^m_0, T]}_{t_n}) , \vec \cK_0(\vec i,  \vec\k; (\vec s, t_n))  \big\ra dB^{j_n}_{t_n}\\
&&\dis  + {1\over 2} \int_{T_m}^{s_{n-1}} \sum_{\vec\k \in \cS_{n}} \Big[\sum_{\tilde \a = 1}^n \prod_{\alpha \neq \tilde \a} \pa^{\vec \k, \alpha}_{\vec i}  V^{i_\a}_{j_\a} (\Th^{s_\a}_{t_n}) \cK_+(\vec i,  \vec \k; (\vec s, t_n)) ~\big \la \pa_{\vec i}^{\vec \k}u(t^m_0, \Th^{[t^m_0, T]}_{t_n}) , \vec \cK_0(\vec i,  \vec\k; (\vec s, t_n))  \big\ra \times\\
&&\dis \qq\qq\qq \sum_{\tilde i=1}^{d_1} \big[\pa_{x_{\tilde i}} \pa^{\vec \k, \tilde \alpha}_{\vec i}  V^{i_{\tilde \a}}_{j_{\tilde \a}}(\Th^{s_{\tilde \a}}_{t_n}) K_{\tilde i}(s_{\tilde \a}, t_n) V^{\tilde i}_{j_n}(X_{t_n})]\\
&&\dis \qq\qq +\prod_{\alpha = 1}^{n} \pa^{\vec \k, \alpha}_{\vec i}  V^{i_\a}_{j_\a}(\Th^{s_\a}_{t_n})\cK_+(\vec i,  \vec \k; (\vec s, t_n)) \times\\
&&\dis\qq\qq\qq \sum_{\tilde i=1}^{d_1} ~\big \la \pa_{\bx_{\tilde i}} \pa_{\vec i}^{\vec \k}u(t^m_0, \Th^{[t^m_0, T]}_{t_n}) , (\vec \cK_0(\vec i,  \vec\k; (\vec s, t_n)), K_{\tilde i, t_n}^{[t^m_0, T]})  \big\ra V^{\tilde i}_{j_n}(X_{t_n}) \Big] dt_n\\
&&\dis = \int_{T_m}^{s_{n-1}} \sum_{\vec\k \in \cS_{n}}\prod_{\alpha = 1}^{n} \pa^{\vec \k, \alpha}_{\vec i}  V^{i_\a}_{j_\a} (\Th^{s_\a}_{t_n}) \cK_+(\vec i,  \vec \k; (\vec s, t_n)) ~\big \la \pa_{\vec i}^{\vec \k}u(t^m_0, \Th^{[t^m_0, T]}_{t_n}) , \vec \cK_0(\vec i,  \vec\k; (\vec s, t_n))  \big\ra dB^{j_n}_{t_n}\\
&&\dis \q + {1\over 2} \sum_{\tilde i=1}^{d_1}  \int_{T_m}^{s_{n-1}}\sum_{\vec\k \in \cS_{n}}\sum_{\tilde \a =0}^n \cV\big((\vec i, \tilde i), (\vec j, j_n), (\vec\k, \tilde \a); (\vec s, t_n, t_n), \Th^{(\vec s, t_n, t_n)}_{t_n}, \Th^{[t^m_0,T]}_{t_n}\big) dt_n,
\eeaa
 where, similarly to \reff{XiN}, the last equality can be verified straightforwardly. Then
\beaa
&&\dis \dbE_m\Big[\Big|\int_{T_m}^{s_{n-1}} \f(\vec s, t_n) \circ dB^{j_n}_{t_n}\Big|^2\Big] \le C \d |A^m_{\|\vec j\|}|^2 + C\d^2 |A^m_{\|(\vec j, j_n)\|}|^2.
\eeaa
Thus
\beaa
\L_{n-1} \le C\d^{\|\vec j_{n-1}\|} \big[\d |A^m_{\|\vec j\|}|^2 + \d^2 |A^m_{\|(\vec j, j_n)\|}|^2\big] = C\big[ |A^m_{\|\vec j\|}|^2 \d^{\|\vec j\|}+  |A^m_{\|\vec j\|+1}|^2\d^{\|\vec j\|+1}\big].
\eeaa
So in all the cases, by \reff{ERNest1} we have
\beaa
\dbE_m\Big[\Big|\int_{\dbT^m_n}\f(\vec t) \circ d B^{\vec j}_{\vec t}\Big|^2 \Big] \le  C_n\Big[ |A^m_{\|\vec j\|}|^2\d^{\|\vec j\|} +  |A^m_{\|\vec j\|+1}|^2\d^{\|\vec j\|+1}\Big].
\eeaa
Then, by \reff{RNrep} we obtain
\beaa
\dbE_m[|R^m_N|^2] \le  C_N \sum_{n=1}^{N+1}  \sum_{\vec i\in \cI_n, \vec j\in \cJ_n}  \big[ \1_{\{\|\vec j\|=N+1\}}  +  \1_{\{\|\vec j_{-1}\|=N, j_1=0\}}\big] \Big[ |A^m_{\|\vec j\|}|^2\d^{\|\vec j\|} +  |A^m_{\|\vec j\|+1}|^2\d^{\|\vec j\|+1}\Big].
\eeaa
This implies \reff{RNest} immediately.
\qed

\bs
\no{\bf Proof of Theorem \ref{thm-Qest}.} It is clear that the functional It\^{o} formula \reff{Ito} holds true under $Q$ as well, then \reff{Taylor1new} and \reff{RNrep} also hold true under $Q$. Thus, by \reff{cubature},
\beaa
\big|Y_0 - Y^{Q}_0\big| = \Big| \dbE[R_N] - \dbE^{Q}[R_N]\Big|.
\eeaa
Therefore, by Theorem \ref{thm-Taylor}, it suffices to show that
\bea
\label{Qest2}
\big|\dbE^{Q}[R_N]\big|\le C_N \Big[  A_{N+1} T^{N+1\over 2} \big(1+ C_{Q}^{N-1}\big) + A_{N+2} T^{N+2\over 2}\big(1+  C_{Q}^{N-2}\big)\Big].
\eea

Now for each $n$ and $\vec i\in \cI_n$, $\vec j\in \cJ_n$ as in \reff{RNrep}, note that 
\beaa
(\o_{k, t}^0)'=1, \q |\{l: j_l =0\}| = \|\vec j\|-n,\q |\{l: j_l > 0\}| =2n- \|\vec j\|,
\eeaa
where $\o'$ denotes the time derivative of $\o$. Then 
\bea
\label{Qest3}
&&\dis \Big|\dbE^{Q}\Big[\int_{\dbT_n}  \sum_{\vec\k \in \cS_{n}} \cV(\vec  i,\vec j, \vec\k;\vec t, \Th^{\vec t}_{ t_{n}},\Th^T_{t_{n}}) \circ dB^{\vec j}_{\vec t}\Big]\Big|\le A_{\|\vec j\|} \sum_{k=1}^{2W} \l_k \int_{\dbT_n}  \prod_{l=1}^n |(\o^{j_l}_{k, t_l})'| dt_n \cds dt_1\Big|\nonumber\\
&&\dis \le C_n A_{\|\vec j\|} \sum_{k=1}^{2W} \l_k ({C_Q\over \sqrt{T}})|^{2n-\|\vec j\|} T^n = C_n A_{\|\vec j\|} C_{Q}^{2n-\|\vec j\|} T^{\|\vec j\|}. 
\eea
If $\|\vec j\|=N+1$, we have $2n \ge N+1$. Denote $k:= 2n-(N+1)$, then $0\le k\le N-1$, and
\beaa
A_{\|\vec j\|} |C_{Q}|^{2n-\|\vec j\|} T^{\|\vec j\|} = A_{N+1} T^{N+1\over 2} C_{Q}^k.
\eeaa
If $\|\vec j_{-1}\|=N$, $j_1=0$, then $2n \ge N+2$. Denote $k:= 2n-(N+2)$, then $0\le k\le N-2$, and
\beaa
A_{\|\vec j\|} C_{Q}^{2n-\|\vec j\|} T^{\|\vec j\|} = A_{N+2} T^{N+2\over 2} C_{Q}^k.
\eeaa
Plug these into \reff{Qest3}, we obtain \reff{Qest2} immediately, and hence \reff{Qest1} holds true.

Finally, when \reff{rescale} holds true,  by \reff{cubature} one can easily see that $Q$  in \reff{QK}  is an N-Volterra cubature formula on $[0, T]$ if and only if the following rescaled one $\tilde Q$ is an N-Volterra cubature formula on $[0, 1]$:
\bea
\label{omegalinear2}
\tilde Q = \sum_{k=1}^{2W} \l_k \d_{\tilde \o_t},\qq d \tilde \o_{k,t} =  a_{k,l} dt, ~ t\in (\tilde s_{l-1}, \tilde s_{l}],\qq \tilde s_l := {l\over L},\q l=0,\cds, L.
\eea 
In particular, this implies that $C_Q = C_{\tilde Q}$ is independent of $T$, 
\qed

\bs

\no{\bf Proof of Theorem \ref{thm-Qestm}.} Note that
\beaa
&&\dis \Big| \dbE^{Q_m}_m\big[u(T_{m+1}, \Th^{[T_{m+1}, T]}_{T_{m+1}}\big] - \dbE_m\big[u(T_{m+1}, \Th^{[T_{m+1}, T]}_{T_{m+1}}\big]\Big|\\
&&\dis \le \Big| \dbE^{Q_m}_m[I^m_N] - \dbE_m[I^m_N]\Big| + \big|\dbE^{Q_m}_m[R^m_N]\big| +  \big|\dbE_m[R^m_N]\big|.
\eeaa
Then, following the same arguments as in Theorem \ref{thm-Qest}, it suffices to provide desired estimate for $\Big| \dbE^{Q_m}_m[I^m_N] - \dbE_m[I^m_N]\Big|$. Similar to Lemma \ref{lem-odd}, by the desired symmetric properties, for any $\vec j\in \cJ_{n,N}\backslash \bar \cJ_{n, N}$ we have
\beaa
\dbE^{Q_m}_m\Big[\int_{\dbT^m_n} \cV(\vec i,\vec j, \vec\k;\vec t, \Th^{\vec t}_{T_m},\Th^{[t^m_0,T]}_{T_m}) \circ dB^{\vec j}_{\vec t} \Big] =0 = \dbE_m\Big[\int_{\dbT^m_n} \cV(\vec i,\vec j, \vec\k;\vec t, \Th^{\vec t}_{T_m},\Th^{[t^m_0,T]}_{T_m}) \circ dB^{\vec j}_{\vec t} \Big].
\eeaa
Then, by \reff{cubaturem} we have
\beaa
\Big|\dbE^{Q_m}_m[I^m_N] - \dbE_m[I^m_N] \Big|=\Big| \dbE^{Q_m}_m[\check R^m_N] - \dbE_m[\check R^m_N]\Big|\le \Big| \dbE^{Q_m}_m[\check R^m_N]\Big| + \Big|\dbE_m[\check R^m_N]\Big|. 
\eeaa
The estimate for $\Big|\dbE_m[\check R^m_N]\Big|$ is implied by \reff{checkRmNest}. Moreover, for each $n\le N$ and $\vec j\in \bar\cJ_{n, N}$, again set $k := {N-\|\vec j\|-1\over 2}$ as in the proof of Theorem \ref{thm-Taylorm}. Then, by \reff{checkRk} and similar to \reff{Qest3},
\beaa
\Big| \dbE^{Q_m}_m\Big[\int_{\dbT^m_n} \check R^m_N(\vec t) \circ dB^{\vec j}_{\vec t} \Big] \Big|\le C^m_N C_{Q_m}^{2n-\|\vec j\|} e^{C^m_NT} \d^{N+1\over 2} .
\eeaa
Thus, since $ 0\le 2n-\|\vec j\|\le N-1$,
\beaa
\Big| \dbE^{Q_m}_m[\check R^m_N]\Big| \le \sum_{n=1}^N \sum_{\vec j\in \bar\cJ_{n, N}} C^m_N C_{Q_m}^{2n-\|\vec j\|} e^{C^m_NT}\d^{N+1\over 2} \le C^m_N\big(1+C_{Q_m}^{N-1}\big) e^{C^m_NT}\d^{N+1\over 2} .
\eeaa
Recall $C_{Q_m} = C_{Q^*_N}$, this is the desired estimate.
\qed

\subsection{The one period cubature formula for \reff{Heston} with $N=5$}
\label{sect-N5-2d2}
Denote
\bea
\label{g12H}
\g_1 :={T^{H_++1}\over H_+ (H_++1)},\q \g_2 :={T^{2H_+}\over 8HH_+}.
\eea
First, similar to \reff{Taylor1} we have
\bea
\label{b1Taylor}
\dis \tilde b_1(\Th^{t_1}_{t_1}) 
\dis &=& R_{\neq 2} + \int_0^{t_1} \Big[\pa_0 \tilde b_1(\Th^{t_1}_{t_2})  +\sum_{i=1}^2 \pa_i \tilde b_1(\Th^{t_1}_{t_2}) b_i(\Th^{t_2}_{t_2}) K_i(t_1, t_2)\Big] dt_2 \nonumber\\
&&\dis +  \int_0^{t_1}\int_0^{t_2} \sum_{i,j=1}^2 \Big[\pa_{ji} \tilde b_1(\Th^{t_1}_{t_3}) \si_i(\Th^{t_2}_{t_3})K_j(t_1, t_3) \\
&&\dis \q +  \pa_{i} \tilde b_1(\Th^{t_1}_{t_3}) \pa_j \si_i(\Th^{t_2}_{t_3}) K_j(t_2, t_3)\Big] \si_j(\Th^{t_3}_{t_3}) K_i(t_1, t_2) \circ dB^j_{t_3}\circ dB^i_{t_2}. \nonumber
\eea
We remark that $\pa_i \tilde b_1 = G'' b_1 + G'\pa_i b_1$, however, thanks to the new convention in Remark \ref{rem-groupI}, we do not need to consider $G'' b_1$ and $G'\pa_i b_1$ separately. Then we can easily see that  the terms in $\D\tilde\dbV^0_4$ derived from $\tilde b_1$ consist of the following stochastic integrals: 
\bea
\label{2dhV4b}
\left.\ba{c}
\dis \D\tilde\dbV^0_{4,b} := \Big\{ \tilde\G_{4,b}^{(i,j,k)}: (i,j,k) = (1,1,1), (1,2,1), (1,2,2), (2,1,1), (2,2,1), (2,2,2)\Big\},\\
\dis \mbox{where}\q \tilde \G_{4,b}^{(i,j,k)}:=  \int_{\dbT_3} K_i(t_1, t_2) K_j(t_k, t_3) \circ d B^j_{t_3} \circ d B^i_{t_2}dt_1.
\ea\right.
\eea
Note that, for $i,j, k = 1, 2$,
{\small\bea
\label{2dE2}
 \dis   \dbE\big[\tilde \G_{4,b}^{(i,j,k)}\big] = {1\over 2}\int_{\dbT_2}  K_i(t_1, t_2) K_j(t_k, t_2) d\la B^i, B^j\ra_{t_2} dt_1 =\left\{\ba{lll} {T^2\over 4},& (i, j, k) = (1,1,1),\\ 
 {\rho \g_1\over 2},& (i, j, k) = (1,2,1), (2,1,1), \\
 \g_2,& (i, j, k) = (2,2,1),\\
 0,& (i, j, k) = (1,2,2), (2,2,2).
\ea\right.
\eea}

Similarly, we may have the expansion of $ \tilde \si_1(\Th^{t_1}_{t_1})$ as in \reff{b1Taylor}:
\bea
\label{si1Taylor1}
&&\dis \tilde \si_1(\Th^{t_1}_{t_1}) = R_{\neq 3} +\sum_{i_2=1}^2 \int_0^{t_1}  \int_0^{t_2}  \pa_{i_20} \tilde \si_1(\Th^{t_1}_{t_3}) \si_{i_2}(\Th^{t_3}_{t_3}) K_{i_2}(t_1, t_3) \circ dB^{i_2}_{t_3} dt_2\nonumber\\
&&\dis  +\sum_{i_1, i_2=1}^2 \int_0^{t_1}  \int_0^{t_2} \Big[\pa_{i_2i_1} \tilde \si_1(\Th^{t_1}_{t_3}) b_{i_1}(\Th^{t_2}_{t_3}) K_{i_2}(t_1, t_3) \nonumber\\
&&\dis\qq  +\pa_{i_1} \tilde \si_1(\Th^{t_1}_{t_3}) \pa_{i_2} b_{i_1}(\Th^{t_2}_{t_3}) K_{i_2}(t_2, t_3) \Big] \si_{i_2}( \Th^{t_3}_{t_3}) K_{i_1}(t_1, t_2) \circ dB^{i_2}_{t_3}dt_2 \\
&&\dis + \sum_{i_1=1}^2\int_0^{t_1} \int_0^{t_2} \Big[\pa_{0i_1} \tilde \si_1(\Th^{t_1}_{t_3}) \si_{i_1}(\Th^{t_2}_{t_3})+\pa_{i_1} \tilde \si_1(\Th^{t_1}_{t_3}) \pa_0\si_{i_1}(\Th^{t_2}_{t_3})\Big]  K_{i_1}(t_1, t_2) dt_3\circ dB^{i_1}_{t_2}\nonumber\\
&&\dis + \sum_{i_1,i_2=1}^2\int_0^{t_1} \int_0^{t_2} \Big[\pa_{i_2i_1} \tilde \si_1(\Th^{t_1}_{t_3})b_{i_2}(\Th_{t_3}^{t_3}) \si_{i_1}(\Th^{t_2}_{t_3})K_{i_2}(t_1, t_3)\nonumber \\
&&\dis\qq +\pa_{i_1} \tilde \si_1(\Th^{t_1}_{t_3}) \pa_{i_2}\si_{i_1}(\Th^{t_2}_{t_3})b_{i_2}(\Th_{t_3}^{t_3})K_{i_2}(t_2, t_3) \Big]  K_{i_1}(t_1, t_2) dt_3\circ dB^{i_1}_{t_2}+\xi, \nonumber
\eea
where the terms presented above involve $\pa_t$, and $\xi$ contains the terms without  $\pa_t$:
\bea
\label{si1Taylor2}
&&\dis\xi:= \sum_{i_1, i_2, i_3=1}^2  \int_0^{t_1}\int_0^{t_2} \int_0^{t_3}\Big[ \pa_{i_3i_2i_1}  \tilde \si_1(\Th^{t_1}_{t_4}) \si_{i_1}(\Th^{t_2}_{t_4}) \si_{i_2}(\Th^{t_3}_{t_4})K_{i_2}(t_1, t_3)K_{i_3}(t_1, t_4)\nonumber\\
&&\dis \qq+ \pa_{i_2i_1}  \tilde \si_1(\Th^{t_1}_{t_4}) \pa_{i_3}\si_{i_1}(\Th^{t_2}_{t_4}) \si_{i_2}(\Th^{t_3}_{t_4})K_{i_2}(t_1, t_3)K_{i_3}(t_2, t_4)\nonumber\\
&&\dis \qq+ \pa_{i_2i_1}  \tilde \si_1(\Th^{t_1}_{t_4}) \si_{i_1}(\Th^{t_2}_{t_4}) \pa_{i_3} \si_{i_2}(\Th^{t_3}_{t_4})K_{i_2}(t_1, t_3)K_{i_3}(t_3, t_4)\nonumber\\
&&\dis \qq+ \pa_{i_3i_1}  \tilde \si_1(\Th^{t_1}_{t_4}) \pa_{i_2}\si_{i_1}(\Th^{t_2}_{t_4}) \si_{i_2}(\Th^{t_3}_{t_4})K_{i_2}(t_2, t_3)K_{i_3}(t_1, t_4)\\
&&\dis \qq+ \pa_{i_1}  \tilde \si_1(\Th^{t_1}_{t_4}) \pa_{i_3i_2}\si_{i_1}(\Th^{t_2}_{t_4}) \si_{i_2}(\Th^{t_3}_{t_4})K_{i_2}(t_2, t_3)K_{i_3}(t_2, t_4)\nonumber\\
&&\dis \qq+ \pa_{i_1}  \tilde \si_1(\Th^{t_1}_{t_4}) \pa_{i_2}\si_{i_1}(\Th^{t_2}_{t_4}) \pa_{i_3}\si_{i_2}(\Th^{t_3}_{t_4})K_{i_2}(t_2, t_3)K_{i_3}(t_3, t_4)\Big] \times\nonumber\\
&&\dis\qq \si_{i_3}(\Th^{t_4}_{t_4}) K_{i_1}(t_1, t_2)\circ dB^{i_3}_{t_4}\circ dB^{i_2}_{t_3}\circ dB^{i_1}_{t_2}.\nonumber
\eea
Then we see that  the terms in $\D\tilde\dbV^0_4$ derived from $\tilde \si_1$ are:
\bea
\label{2dhVsi1-1}
\left.\ba{lll} 
\dis \D\tilde\dbV^0_{4,\si}  =\D\tilde\dbV^0_{4,\si, 0}\cup \D\tilde\dbV^0_{4,\si, 1} \cup \D\tilde\dbV^0_{4,\si,2} \cup \D\tilde\dbV^0_{4,\si,3},\q\mbox{where, for}~l=1,2,\ms\\
\dis \D\tilde\dbV^0_{4,\si,0}  := \Big\{\tilde \G^{(1)}_{4,\si, 0}, \tilde \G^{(2)}_{4,\si, 0}\Big\},\ms\\ 
\dis \D\tilde \dbV^0_{4,\si,l} := \Big\{ \tilde \G^{(i_1, i_2, \k_2)}_{4,\si, l}: (i_1, i_2, \k_2)=(1,1,1),(1,2,1),(1,2,2),(2,1,1),(2,2,1),(2,2,2)\Big\},\ms\\
\ea\right.
\eea
and
\bea
\label{2dhVsi1-2}
\left.\ba{lll} 
\dis \D\tilde\dbV^0_{4,\si,3} := \Big\{\tilde \G^{(i_1, i_2, i_3, \k_2, \k_3)}_{4,\si, 3}: (i_1, i_2, i_3, \k_2, \k_3)=(1,1,1,1,1),(1,1,2,1,1), \\
\dis\qq (1,1,2,1,2),(1,1,2,1,3),(1,2,1,1,1),(1,2,1,2,1), (1,2,2,1,1),(1,2,2,1,2), \\
\dis\qq (1,2,2,1,3),(1,2,2,2,1),(1,2,2,2,2),(1,2,2,2,3), (2,1,1,1,1),(2,1,2,1,1), \\
\dis\qq (2,1,2,1,2),(2,1,2,1,3),(2,2,1,1,1),(2,2,1,2,1) ,(2,2,2,1,1),(2,2,2,1,2), \\
\dis\qq (2,2,2,1,3), (2,2,2,2,1),(2,2,2,2,2),(2,2,2,2,3)  \Big\},
\ea\right.
\eea
where
\beaa
\left.\ba{lll} 
\dis \tilde \G^{(i_2)}_{4,\si, 0} := \int_{\dbT_3} K_{i_2}(t_1, t_3) \circ dB^{i_2}_{t_3} dt_2 \circ dB^1_{t_1},\ms\\
\dis \tilde \G^{(i_1, i_2, \k_2)}_{4,\si, 1} := \int_{\dbT_3} K_{i_1}(t_1, t_2)K_{i_2}(t_{\k_2}, t_3) \circ dB^{i_2}_{t_3} dt_2 \circ dB^1_{t_1},\ms\\
\dis \tilde \G^{(i_1, i_2, \k_2)}_{4,\si, 2} := \int_{\dbT_3} K_{i_1}(t_1, t_2) K_{i_2}(t_{\kappa_2}, t_3)dt_3 \circ dB^{i_1}_{t_2} \circ dB^1_{t_1},\ms\\
\dis \tilde \G^{(i_1, i_2, i_3, \k_2, \k_3)}_{4,\si, 3} := \int_{\dbT_4}K_{i_1}(t_1, t_2) K_{i_2}(t_{\k_2}, t_3) K_{i_3}(t_{\k_3}, t_4)dB^{i_3}_{t_4} \circ dB^{i_2}_{t_3} \circ dB^{i_1}_{t_2} \circ dB^1_{t_1}.
\ea\right.
\eeaa
Note that, for $\vec i = (i_1, i_2, i_3)$ and $\vec \k = (\k_1, \k_2, \k_3)$ with $\k_1=1$, $\k_2=1,2$ and $\k_3 = 1,2,3$, we have: recalling \reff{g12H}, 
\bea
\label{2dE3}
\left.\ba{lll}
\dis \dbE\big[ \tilde \G^{(i_2)}_{4,\si, 0}\big] = \dbE\big[ \tilde \G^{(i_1, i_2, \k_2)}_{4,\si, 1}\big]=0,\ms\\
\dis \dbE\big[\tilde \G^{(i_1, i_2, \k_2)}_{4,\si, 2}\big] = {1\over 2} \int_{\dbT_1} K_{i_1}(t_1, t_1)\int_0^{t_1}K_{i_2}(t_{1},t_3)d t_3 d\la B^{i_1},B^1\ra_{t_1},\\
\dis\q = \left\{\ba{lll} {T^2\over 4}, \quad (i_1, i_2, \k_2)= (1,1,1),\\ {\rho\gamma_1\over 2} ,\quad (i_1, i_2, \k_2)=(1,2,1), \\ 0, \quad (i_1, i_2, \k_2)= (2,1,1),  (2,2,1),(2,2,2), (1,2,2),\ea\right.\ms\\
\dis \dbE\big[\tilde \G^{(i_1, i_2, i_3, \k_2, \k_3)}_{4,\si, 3}\big]\\
\dis \q= {1\over 4} \int_0^T \int_0^{t_1} K_{i_1}(t_1, t_1) K_{i_2}(t_{\k_2}, t_3) K_{i_3}(t_{\k_3}, t_3)\big|_{t_2=t_1} d\la B^{i_3}, B^{i_2}\ra_{t_3} d\la B^{i_1}, B^{1}\ra_{t_1}\ms\\
\dis\q = \left\{\ba{lll} {T^2\over 8},~ \vec i = (1,1,1),\ms\\ {\rho \g_1\over 4 },~ \vec i=(1,2,1), ~\mbox{or}~ \vec i = (1, 1, 2)~\mbox{and}~ \k_3 = 1,2,\ms\\ {\g_2\over 2}, ~ \vec i = (1, 2, 2), ~\k_3 = 1, 2,\ms\\0,~ i_1 = 2 ~\mbox{or}~(i_3, \k_3) = (2, 3).
\ea\right.
\ea\right.
\eea
We note that, by \reff{2dhV2}, \reff{2dhV4b}, and \reff{2dhVsi1-1}-\reff{2dhVsi1-2},
\bea
\label{DhVsize}
|\D\tilde\dbV^0_2| = 2,\quad \text{and}\quad|\D\tilde \dbV^0_4| = |\D\tilde\dbV^0_{4, b}| + |\D\tilde \dbV^0_{4, \si}|=6+2+6+6+24=44.
\eea

We now construct a desired $Q$ as in \reff{QK}. Recall $\o^0_t = t$, and $\o_{k,t} = (\o^1_{k,t}, \o^2_{k, t})$.   We remark that, the correlation between $B^1, B^2$ affects only the expectations $\dbE[\tilde\G]$ in \reff{2dE2} and  \reff{2dE3}, but the expectations under $Q$ in \reff{EQK} remain the same. The integrals against $d\o$ corresponding to the stochastic integrals in \reff{2dhV2}, \reff{2dhV4b}, and \reff{2dhVsi1-1}-\reff{2dhVsi1-2} are: for each $k=1,\cds, W$ and for appropriate functions $\f$,
\bea
\label{2ddo}
\left.\ba{lll}
\dis \int_{\dbT_2} \f(\vec t~\!)d \o^{i_1}_{k, t_2} d \o^1_{k, t_1} = \sum_{1\le l_2 \le l_1 \le m}  {a_{k, l_1}^1 a_{k, l_2}^{i_1}\over T} \int_{s_{l_1-1}}^{s_{l_1}}\int_{s_{l_2-1}}^{s_{l_2}\wedge t_1}  \f(\vec t~\!) dt_2 dt_1,\ms\\
\dis \int_{\dbT_3}  \f(\vec t~\!) d \o^{i_2}_{k, t_3}  d \o^{i_1}_{k, t_2}dt_1=  \sum_{1\le l_2 \le l_1 \le m}  {a_{k, l_1}^{i_1} a_{k, l_2}^{i_2} \over T}\int_{s_{l_1-1}}^{s_{l_1}}\int_{s_{l_2-1}}^{s_{l_2}\wedge t_2} \int_{t_2}^T   \f(\vec t~\!) dt_1 dt_3 dt_2,\ms\\
\dis \int_{\dbT_3} \f(\vec t~\!) d\o^{i_2}_{k, t_3} dt_2  d\o^1_{k,t_1}=\sum_{1\le l_2 \le l_1 \le m}  {a_{k, l_1}^1 a_{k, l_2}^{i_2}\over T} \int_{s_{l_1-1}}^{s_{l_1}}\int_{s_{l_2-1}}^{s_{l_2}\wedge t_1} \int_{t_3}^{t_1}  \f(\vec t~\!)dt_2dt_3 dt_1,\ms\\ 
\dis\int_{\dbT_3} \f(\vec t~\!) dt_3  d\o^{i_1}_{k,t_2} d\o^1_{k, t_1} = \sum_{1\le l_2 \le l_1 \le m}  {a_{k, l_1}^1 a_{k, l_2}^{i_1}\over T} \int_{s_{l_1-1}}^{s_{l_1}}\int_{s_{l_2-1}}^{s_{l_2}\wedge t_1}\int_0^{t_2} \f(\vec t~\!) dt_3dt_2 dt_1,\ms\\
\dis \int_{\dbT_4}\f(\vec t~\!)d\o^{i_3}_{k, t_4} d\o^{i_2}_{k, t_3}  d\o^{i_1}_{k, t_2} d\o^1_{k, t_1} =\sum_{1\le l_4\le l_3 \le l_2\le l_1\le 4}  {a_{k, l_1}^{i_1} a_{k, l_2}^{i_2}a_{k, l_3}^{i_3}a_{k, l_3}^{i_4} \over T^2}\times\ms\\
\dis\qq\qq\qq\qq\qq \int_{s_{l_1-1}}^{s_{l_1}} \int_{s_{l_2-1}}^{s_{l_2}\wedge t_1} \int_{s_{l_3-1}}^{s_{l_3}\wedge t_2}\int_{s_{l_4-1}}^{s_{l_4}\wedge t_3}  \f(\vec t~\!)dt_4dt_3dt_2dt_1.
\ea\right.
\eea

\ms

For $N=5$, by \reff{DhVsize}, together with the constraint on $\l_k$, there will be in total $47$ equations. In our example \reff{RVM}, however,  the coefficients are homogeneous, i.e. independent of $t$. Then in \reff{si1Taylor1} the terms $\pa_{i_20}\tilde \si_1 = 0$ and thus the two terms  $\tilde \G^{(i_2)}_{4,\si, 0} \in \D\tilde\dbV^0_{4,\si,0}$, $i_2=1,2$, in \reff{2dhVsi1-1} are not needed (the terms $\int_{\dbT_3} K_{i_1}(t_1, t_2) dt_3\circ dB^{i_1}_{t_2}\circ dB^1_{t_1}$ are stilled needed, even though $\pa_{0i_1}\tilde \si_1 =\pa_{0}\tilde \si_1 = 0$ in \reff{si1Taylor1}, because they appear in the last line of \reff{si1Taylor1} as well when $i_2=1$). Therefore, we have a total of $45$ equations. We thus set $W=5$ and $L=4$. Since each $2$-dimensional path involves $8$ parameters, plus the parameters $\l_k$, there will be $45$ unknowns. For each $\tilde\G \in \D\tilde\dbV_2^0\cup (\D\tilde\dbV_4^0\backslash \D\tilde\dbV^0_{4,\si,0})$, by \reff{chainrule}, \reff{b1Taylor}, \reff{si1Taylor1}, \reff{si1Taylor2}, one can easily derive from \reff{2ddo} the right side of  \reff{cubatureI},  as $2$-nd order or $4$-th order polynomials of $a_{k,l}^i$. This, together with  \reff{2dhV2}, \reff{2dE2}, \reff{2dE3}, as well as $\l_1 +\cds+\l_5 = {1\over 2}$, leads to the required $45$ equations, in exactly the same manner as in \reff{Equation2}.  Again, we are not able to solve these equations explicitly, so we will solve them numerically.

\end{document}